\documentclass[11pt,reqno]{amsart}
\usepackage[margin=1.2in]{geometry}
\usepackage[utf8]{inputenc}
\usepackage{amsmath,amssymb,amsthm,mathrsfs}
\usepackage{enumitem}
\usepackage{tikz,tikz-cd}
\usepackage[colorlinks,linkcolor=blue,citecolor=blue,urlcolor=blue]{hyperref}

\theoremstyle{plain}
\newtheorem{theorem}{Theorem}[section]
\newtheorem{lemma}{Lemma}[section]
\newtheorem{proposition}{Proposition}[section]
\newtheorem{corollary}{Corollary}[section]
\newtheorem{conjecture}{Conjecture}[section]
\theoremstyle{definition}
\newtheorem{definition}{Definition}
\newtheorem{example}{Example}[section]
\theoremstyle{remark}
\newtheorem*{remark}{Remark}
\newtheorem*{remarks}{Remarks}
\numberwithin{equation}{section}
\setlist[enumerate]{itemsep=2pt,topsep=4pt}
\newcommand{\cD}{\mathcal{D}}
\newcommand{\cO}{\mathcal{O}}
\newcommand{\cC}{\mathcal{C}}
\newcommand{\cL}{\mathcal{L}}
\newcommand{\cR}{\mathcal{R}}
\newcommand{\cT}{\mathcal{T}}
\newcommand{\C}{\mathbb{C}}
\newcommand{\Z}{\mathbb{Z}}

\newcommand{\A}{\mathbb{A}}
\newcommand{\del}{\partial}
\newcommand{\fm}{\mathfrak{m}}
\newcommand{\fn}{\mathfrak{n}}
\newcommand{\ml}{\mathfrak{m}_\lambda}
\newcommand{\nl}{\mathfrak{n}_\lambda}
\newcommand{\hO}{\widehat{\cO}}
\newcommand{\hD}{\widehat{\cD}}
\newcommand{\hn}{\widehat{\fn}}
\newcommand{\bD}{\overline{\cD}}
\newcommand{\bM}{\overline{M}}
\newcommand{\Dl}{\Delta_\lambda}
\newcommand{\Gra}{\mathrm{Gr}^{\mathrm{ad}}}

\newcommand{\ff}[2]{[#1]_{#2}}          
\newcommand{\lan}{\langle}
\newcommand{\ran}{\rangle}
\DeclareMathOperator{\Ann}{Ann}
\DeclareMathOperator{\cl}{cl}
\DeclareMathOperator{\Pic}{Pic}
\DeclareMathOperator{\Hom}{Hom}
\DeclareMathOperator{\End}{End}
\DeclareMathOperator{\Ext}{Ext}
\DeclareMathOperator{\Der}{Der}
\DeclareMathOperator{\Aut}{Aut}
\DeclareMathOperator{\GL}{GL}
\DeclareMathOperator{\gr}{gr}
\DeclareMathOperator{\Ker}{Ker}
\DeclareMathOperator{\Coker}{Coker}
\DeclareMathOperator{\im}{Im}
\DeclareMathOperator{\ord}{ord}
\DeclareMathOperator{\codim}{codim}
\DeclareMathOperator{\Wr}{Wr}
\DeclareMathOperator{\res}{res}
\DeclareMathOperator{\Res}{Res}
\DeclareMathOperator{\Supp}{Supp}
\DeclareMathOperator{\spn}{span}
\DeclareMathOperator{\Stairs}{Stairs}
\DeclareMathOperator{\inn}{in}

\begin{document}

\title{Strongly holonomic functions on curves}

\author{Yuri Berest}
\address{Department of Mathematics, Cornell University, Ithaca, NY 14853-4201, USA}
\email{berest@math.cornell.edu}

\author{Isaac Goldberg}
\address{Department of Mathematics, Cornell University, Ithaca, NY 14853-4201, USA}
\email{isg26@cornell.edu}


\begin{abstract}
Let $X$ be a smooth affine irreducible curve over $\C$, and let $\cD = \cD(X)$ be the ring of
global differential operators on $X$. To every nonzero operator $L \in \cD$ we attach a
non-negative integer $n(L)$, the \emph{holonomicity index} of $L$, which is defined in an elementary way, by counting solutions of $L$ in the truncated local rings at the singular points of $L$. Our main theorem shows that $n(L)$ coincides with the Calogero--Moser index of the Weyl closure $\cl(L) = \cD(K) L \cap \cD$ of $L$, i.e.\ with the integer $n$ labelling the Calogero--Moser space $\cC_n(X,\mathcal{I})$ that contains the isomorphism class of the left ideal $\cl(L)$
under the classification of ideals of $\cD$ given in [Berest--Chalykh, 2012]. 
The proof is based on a dual version of the Cannings--Holland correspondence, which describes left ideals of $\cD$ in terms of finite-dimensional spaces of algebraic distributions on $X$, and an
explicit formula for the Calogero--Moser index given in terms of these distributions. We apply these results to
annihilators of holomorphic functions on $X$: we call a holonomic function $f$ \emph{strongly holonomic} if the Calogero--Moser index of $\Ann_{\cD}(f)$ is positive, and we compute this index in a number of interesting cases.
\end{abstract}

\maketitle

\section{Introduction}
\label{S1}

Let $X$ be a smooth affine irreducible curve over $\C$, and let $\cD = \cD(X)$ be the ring of
global algebraic differential operators on $X$. Given a holomorphic function $f$ on $X$, it is natural to ask which linear differential equations with regular coefficients $f$ satisfies. The set of all such equations is the \emph{annihilator} of $f$,
\[
  \Ann_{\cD}(f) := \{ P \in \cD \,:\, P(f) = 0 \} ,
\]
which is a left ideal of $\cD$. A function $f$ is called \emph{holonomic} if
$\Ann_{\cD}(f) \ne 0$, i.e.\ if $f$ satisfies at least one nontrivial equation. The notion of
holonomicity (for varieties of an arbitrary dimension) goes back to the fundamental work of J.~Bernstein \cite{Bernstein}; for curves, it is equivalent to the above condition. For example, if $X = \A^1$ and $\cD = A_1(\C)$ is the first Weyl algebra, the
functions $e^x$ and $x^n$ ($n \ge 1$) are holonomic, with
\begin{equation}
\label{E11}
  \Ann(e^x) = A_1\,(\del - 1) , \qquad \Ann(x^n) = A_1\,\del^{n+1} + A_1\,(x\del - n) ,
\end{equation}
while $\Gamma(x)$ and $e^{e^x}$ are not (see \cite[Chapter 5, Section~3]{Coutinho1995}). The two ideals in
\eqref{E11} are different in nature: the first is cyclic (i.e., generated by a single element in 
$\cD$), while the second is not. Being cyclic
is a property of the isomorphism class of an ideal, and this suggests a refinement of the above
question: 
\begin{center}
\emph{What is the class of $\Ann_{\cD}(f)$ as a $\cD$-module, and which
properties of $f$ does it encode?}
\end{center}
For a curve $X$, this question is meaningful because the classes of (left) ideals of
$\cD(X)$ admit a simple geometric classification. Let $\cL(\cD)$ denote the set of isomorphism classes
of (nonzero) left ideals of $\cD$. Since $\cD$ is a hereditary Noetherian domain, every ideal is
a projective module of rank $1$, and taking its class in $K_0(\cD) \cong K_0(X) \cong \Z \oplus
\Pic(X)$ defines a map $\gamma: \cL(\cD) \to \Pic(X)$. By \cite{BerestChalykhIdeals}
(generalizing \cite{BerestWilson00, BerestWilsonNC} in the case $X = \A^1$), each fibre of
$\gamma$ decomposes into a countable union of quotients of smooth affine varieties:
\begin{equation}
\label{E12}
  \omega : \ \bigsqcup_{n \ge 0}\,\cC_n(X, \mathcal{I})/\!\Aut_X({\mathcal I})
           \ \stackrel{\sim}{\longrightarrow}\ \gamma^{-1}[\mathcal{I}] ,
           \qquad [\mathcal{I}] \in \Pic(X) ,
\end{equation}
called the \emph{Calogero--Moser spaces} $\,
\cC_{n}(X, \mathcal{I})$ (see Section~\ref{S2.5}). Thus, to every nonzero left ideal
$M \subseteq \cD$ one can attach a non-negative integer $n(M)$ --- the index $n$ of the stratum
in \eqref{E12} containing the orbit $\omega^{-1}(M)$ --- which we call the \emph{Calogero--Moser index} of
$M$. This invariant has a simple meaning at its lowest value: $n(M) = 0$ if and only if
$M \cong \cD \otimes_{\cO}\mathcal{I}$ for a line bundle $\mathcal{I}$ on $X$ (see
Proposition~\ref{P2.3}); in particular, in the case when $X = \A^1$, we have $\,n(M) = 0$ if and only if $M$ is cyclic. This motivates the following definition.

\begin{definition}
\label{D1}
A (locally) holomorphic function $f$ on $X$ is called \emph{strongly holonomic} if it is holonomic and
$n(\Ann_{\cD}(f)) \ge 1$.
\end{definition}

Thus, in the above examples,  $ x^n$ is strongly holonomic on $\A^1$ (in fact, $n(\Ann(x^n)) = n$), while $ e^x$ is not ($n(\Ann(e^x))=0$). In general, however, the answer is not so obvious: for example, what is $n(\Ann(f))$ for $ f = {}_{2}F_1(a,b,c; x) $? (See Example~\ref{Ex6.9}).

\vspace{1ex}

Our motivation for studying  holonomicity comes from knot theory. To a knot $K \subset S^3$ one  can associate its coloured Jones function $J_K: \Z \to \Z[q^{\pm 1}]$, $k \mapsto J_K(k)$, a fundamental quantum invariant of $K$ generalizing the classical Jones polynomial. 
Viewing $q \in \C^*$ as a generic complex parameter, one can make the space of $\C$-valued functions on 
$\Z$ a module over the algebra $ \cT_q = \C\lan \mathsf{L}^{\pm 1}, \mathsf{M}^{\pm 1} \ran/(\mathsf{L} \mathsf{M} = q^2 \mathsf{M} \mathsf{L})$ called the quantum torus, whose generators $\mathsf{L}$ and $\mathsf{M}$ act by
$$
(\mathsf{L} f)(k) = f(k+1) \quad \mbox{and}\quad (\mathsf{M}f)(k) = q^{2k} f(k)
$$
It is known that $\Ann_{\cT_q}(J_K) \ne 0$ for every knot $K$ (see \cite{GaroufalidisLe2005} and also \cite{Le16}, where this property is
called $q$-holonomicity). The algebra $\cT_q$ shares many properties with rings of differential operators on curves\footnote{Indeed, $\cT_q$ is often referred to as the quantum Weyl algebra and denoted $A_q(\C)$.}: in particular, the left ideals of $\cT_q$ are classified by  multiplicative (`quantum') analogues of the Calogero--Moser spaces (see \cite{BCquantum, Berest2013picard}). It is therefore natural to ask
which stratum of the classifying space the ideal $\Ann_{\cT_q}(J_K)$ belongs to. Preliminary
calculations suggest the following.

\begin{conjecture}[Berest]
\label{Conj}
For any knot $K \subset S^3$, the $($unreduced$)$ coloured Jones function $J_K$ is strongly holonomic;
in fact, $n(\Ann_{\cT_q}(J_K)) \ge 2$.
\end{conjecture}

We hope to address this conjecture elsewhere. In the present paper we study the differential
case, where the question is already interesting: 
\emph{How to compute the Calogero--Moser index of
$\Ann_{\cD}(f)$ in terms of $f$?} The definition of $n(M)$ via the bijection \eqref{E12} is not
well suited for this purpose: the inverse of $\omega$ is constructed by means of noncommutative
projective geometry (see \cite{BerestWilsonNC, BerestChalykhIdeals}), and even the simpler
characterization of $n(M)$ as the codimension of the associated graded module $\gr M$ in its
reflexive hull (see Theorem~\ref{T2.3}) requires  knowledge of a Gr\"obner basis of $M$. We will
show that, for ideals arising as annihilators, the index can be computed explicitly --- in elementary terms --- from a single differential operator.

We begin with an overview of our main results. Recall that the \emph{(Weyl) closure} of an
operator $L \in \cD$ is the left ideal $\cl(L) := \cD(K) L \cap \cD$, where $K = K(X)$ is the
field of rational functions on $X$ and $\cD(K) = K \otimes_{\cO} \cD$ is the ring of rational
differential operators. Every annihilator is closed in this sense: if $f$ is holonomic and $L$ is
a nonzero operator of minimal order in $\Ann_{\cD}(f)$, then $\Ann_{\cD}(f) = \cl(L)$ (see
Lemma~\ref{L2.4}). Our problem thus reduces to computing the Calogero--Moser index of $\cl(L)$
for an arbitrary $L \in \cD$.

\subsection*{The holonomicity index}
Fix a nonzero regular differential operator $L \in \cD$ on $X$. For a point $\lambda \in X$, 
denote by $\cO_\lambda$  the local ring of $X$ at $\lambda$, with maximal ideal $\nl$, and let $t \in \nl$ be a local parameter in $\cO_\lambda$ (i.e., a cyclic generator of $\nl$). Then, there is a
unique (minimal) integer $w = w_{\lambda}(L)$, depending on $\lambda$, such that the operator $\bar L := t^{w} L$ preserves the $\nl$-adic filtration
on $\cO_\lambda$, i.e.\ $\bar L(\nl^j) \subseteq \nl^j$ for all $j \ge 0$
(see Lemma~\ref{L3.1}); we call $\bar L$ the \emph{$\lambda$-normalization} of $L$. The operator
$\bar L$ acts on the finite-dimensional vector spaces $\cO_\lambda/\nl^{j+1}$, and we set
\begin{equation}
\label{E13}
  S^j_\lambda(L) := \Ker\bigl[\, \bar L : \ \cO_\lambda/\nl^{j+1} \to \cO_\lambda/\nl^{j+1} \,\bigr] ,
  \qquad j \ge 0 ,
\end{equation}
with  $S^{-1}_\lambda(L) := 0$. It turns out that the dimensions of these spaces form a monotone sequence, which increases by at most one at each step and eventually stabilizes
(Lemma~\ref{L3.3}). We call the finite set
$J_\lambda(L) = \{ j_0 < j_1 < \ldots < j_{s} \}$ of indices $j \ge 0$ at which
$\dim S^j_\lambda(L) > \dim S^{j-1}_\lambda(L)$ the \emph{cotype} of $L$ at $\lambda$, and we
define the \emph{local holonomicity index} of $L$ at $\lambda$ by
\begin{equation}
\label{E14}
  n_\lambda(L) := \sum_{i=0}^{s} (j_i - i) .
\end{equation}
This number vanishes for all but finitely many $\lambda$: specifically, $n_\lambda(L) = 0$ outside
the finite set $\Lambda(L)$ of \emph{singular points} of $L$ (see Definition~\ref{D3}). The
\emph{holonomicity index} of $L$ is then defined by
\begin{equation}
\label{E15}
  n(L) := \sum_{\lambda \in \Lambda(L)} n_\lambda(L) .
\end{equation}
Note that the definition of $n(L)$ involves only linear algebra in the local rings of $X$.

\subsection*{Spaces of conditions}
The main result of this paper (Theorem~\ref{T1.3} below) says that $n(L)$ is the Calogero--Moser
index of $\cl(L)$. We prove the equality $n(\cl(L)) = n(L) $ in two steps, each of which is of independent interest. First, we
give a new description of left ideals of $\cD$, which is dual to the well-known classification of
Cannings and Holland \cite{CanningsHolland1994}. Recall that the latter identifies the (fat)
right ideals of $\cD$ with primary decomposable subspaces of $\cO = \cO(X)$: to a subspace
$W \subseteq \cO$ it assigns the right ideal $\cD(\cO, W)$ of all operators mapping $\cO$ into $W$
(see Section~\ref{S2.4}). Dually, every left ideal of $\cD$ is isomorphic to
$\cD(W, \cO) = \{ P \in \cD(K) : P(W) \subseteq \cO \}$ for some primary decomposable $W$. 
On the other hand, extending a construction of G. Wilson \cite{Wilson1993} to arbitrary curves,
we can describe an (infinite-dimensional) primary decomposable subspace in $\cO $ by a finite-dimensional space of linear conditions. To this end, we consider,
for each $\lambda \in X$,  the right $\cD$-module $\Dl := \cD/\ml \cD$, where $\ml \subset \cO$ is
the maximal ideal of $\lambda$; this module is naturally paired with $\cO$,
\begin{equation}
\label{E16}
  \lan\,\text{--}\,,\text{--}\,\ran_\lambda : \ \Dl \times \cO \to \cO/\ml \cong \C ,
  \qquad \lan F, f \ran_\lambda := F(f) \ (\mathrm{mod}\ \ml) ,
\end{equation}
so that its elements may be viewed as linear functionals on $\cO$ supported at point $\lambda$. A
\emph{homogeneous space of conditions} is a finite-dimensional subspace
$C = \bigoplus_{\lambda} C_\lambda$ of $\bigoplus_{\lambda \in X} \Dl$, and the corresponding
primary decomposable subspace is
\begin{equation}
\label{E17}
  W_C := \{ f \in \cO \, : \, \lan F, f \ran_\lambda = 0 \ \text{ for all } \lambda \in X
  \text{ and } F \in C_\lambda \} .
\end{equation}
Our key observation is that a left ideal $I \subseteq \cD$ determines a space of conditions
directly. Choose an operator $L \in I$ of minimal order, and for each $\lambda \in X$ put
\begin{equation}
\label{E18}
  V_\lambda(I; L) := \Ker\bigl[\, \Dl \xrightarrow{\ \cdot L \ } I/\ml I \,\bigr] ,
  \qquad V_\lambda(L) := \Ker\bigl[\, \Dl \xrightarrow{\ \cdot L \ } \Dl \,\bigr] ,
\end{equation}
where the maps are induced by right multiplication by $L$. These spaces are finite-dimensional,
and they vanish for all but finitely many $\lambda$. Moreover, $V_\lambda(I;L) \subseteq V_\lambda(L)$
with equality when $I = \cl(L)$.

\begin{theorem}
\label{T1.1}
Let $I$ be a nonzero left ideal of $\cD$, and let $L \in I$ be an operator of minimal order. Let
$W = W_C$ be the primary decomposable subspace of $\cO$ associated to the space of conditions
$C = \bigoplus_{\lambda \in X} V_\lambda(I;L)$. Then
\[
  I = \cD(W, \cO) \cdot L .
\]
In particular, $I \cong \cD(W,\cO)$ as left $\cD$-modules, and the dual right ideal of
$I L^{-1}$ is the annihilator of $C$ in $\cD$.
\end{theorem}

Theorem~\ref{T1.1} shows that $W$ is the primary decomposable subspace that the
Cannings--Holland correspondence assigns to the right ideal dual to $I$. As a consequence,
isomorphism classes of left ideals of $\cD$ are parametrized by equivalence classes of
finite-dimensional homogeneous spaces of conditions, which may be thought of as points of a
`dual adelic Grassmannian' (see Corollary~\ref{C4.3}). Along the way, we obtain an explicit
presentation of the closure of an operator on an arbitrary curve.

\begin{corollary}
\label{C1.1}
For any nonzero $L \in \cD$,
\[
  \cl(L) = \cD L + \sum_{\lambda \in \Lambda(L)}\, \cD\, \ml^{-1}\, \widetilde{V}_\lambda(L)\, L ,
\]
where $\widetilde{V}_\lambda(L) = \{F \in \cD : F L \in \ml \cD\}$ is the preimage of
$V_\lambda(L)$ in $\cD$ and $\ml^{-1} \subset K$ is the fractional ideal inverse to $\ml$.
\end{corollary}

For $X = \A^1$, Corollary~\ref{C1.1} is a result of H. Tsai (see \cite[Algorithm~3.4]{Tsai}).

Our second step is a formula for the Calogero--Moser index in terms of spaces of conditions. The
order filtration on $\cD$ induces a filtration on each $\Dl$, so that we may speak of the degree
of a nonzero element of $\Dl$ (if $t$ is a local parameter at $\lambda$ and $\del_t$ is the
derivation with $\del_t(t) = 1$, then $\Dl \cong \C[\del_t]$, and $\deg$ is the usual notion of degree of
polynomials).

\begin{theorem}
\label{T1.2}
Let $C = \bigoplus_{\lambda} C_\lambda$ be a finite-dimensional homogeneous space of conditions.
For each $\lambda$ with $C_\lambda \ne 0$, choose a basis $\{F_{\lambda,1}, \ldots,
F_{\lambda, g_\lambda}\}$ of $C_\lambda$ consisting of elements of pairwise distinct degrees.
Then the Calogero--Moser index of the left ideal $I = \cD(W_C, \cO)$ is given by
\begin{equation}
\label{E19}
  n(I) = \sum_{\lambda}\, \Bigl[\, \sum_{i=1}^{g_\lambda} \deg(F_{\lambda,i})\, -\,
  \frac{g_\lambda\,(g_\lambda - 1)}{2} \,\Bigr] .
\end{equation}
\end{theorem}

Formula \eqref{E19} is dual to the formula for $n$ in terms of the `gaps' of a primary
decomposable subspace given in \cite[(8.3)]{BerestWilson04} (see the Remark at the end of
Section~\ref{S5}).

\subsection*{Main theorem}
The link between the two descriptions is provided by the following duality
(Theorem~\ref{T3.1}). Recall that the holonomicity index of $L$ is defined in terms of the
solution spaces $S^j_\lambda(L)$ of the $\lambda$-normalization $\bar L$ of $L$ in the truncated
local rings $\cO_\lambda/\nl^{j+1}$ (see \eqref{E13}), while the spaces of conditions
$V_\lambda(L) \subseteq \Dl$ are defined by right multiplication by $L$ (see \eqref{E18}). Now, if
$L$ is $\lambda$-normalized (so that $\bar L = L$), then $\dim S^j_\lambda(L)$ equals the
dimension of the subspace of elements of degree $\le j$ in $V_\lambda(L)$. Hence the cotype of $L$ at $\lambda$ is nothing
but the set of degrees of elements of $V_\lambda(\bar L)$, and the local holonomicity index
\eqref{E14} is the $\lambda$-term of \eqref{E19}. Combining this with Theorems~\ref{T1.1}
and~\ref{T1.2}, we arrive at the main result of the paper.

\begin{theorem}
\label{T1.3}
Let $L$ be a nonzero differential operator on $X$. For each $\lambda \in X$, let
$J_\lambda(L) = \{ j_0 < j_1 < \ldots < j_{s_\lambda} \}$ be the cotype of $L$ at $\lambda$, i.e.\
the set of indices $j \ge 0$ at which the dimension of the solution spaces \eqref{E13} of the
$\lambda$-normalization of $L$ jumps, and let
\[
  n(L) = \sum_{\lambda \in \Lambda(L)}\, \sum_{i=0}^{s_\lambda} (j_i - i)
\]
be the holonomicity index of $L$ defined by \eqref{E14} and \eqref{E15}. Then $n(L)$ equals the
Calogero--Moser index of the closure $\cl(L) = \cD(K) L \cap \cD$ of $L$:
\[
  n(L) = n(\cl(L)) .
\]
\end{theorem}

\begin{corollary}
\label{C1.2}
Let $f$ be a holonomic function on $X$, and let $L$ be a nonzero operator of minimal order in
$\Ann_{\cD}(f)$. Then $n(\Ann_{\cD}(f)) = n(L)$. In particular, $f$ is strongly holonomic if and
only if $n_\lambda(L) > 0$ for at least one singular point $\lambda$ of $L$.
\end{corollary}

We illustrate Theorem~\ref{T1.3} by a number of examples in Section~\ref{S6}. Here we mention
only one result, which gives a transparent geometric meaning to the Calogero--Moser index in the
simplest case.

\begin{theorem}
\label{T1.4}
Let $f$ be a nonzero rational function on $X$ or, more generally, a nonzero holomorphic function
on $X$ annihilated by a differential operator of order one. Then $n(\Ann_{\cD}(f))$ equals the
number of zeros of $f$ in $X$ counted with multiplicities. In particular, such a function is
strongly holonomic if and only if it vanishes somewhere on $X$.
\end{theorem}

For example, if $f = p/q$ is a rational function on $\A^1$ written in lowest terms, then
$n(\Ann(f)) = \deg(p)$; if $f$ is a regular function on $X$ defining a finite map
$X \to \A^1$, then $n(\Ann_{\cD}(f))$ is the degree of this map. For a finite set
of functions $\{f_1, \ldots, f_k\}$, the role of $f$ in the above theorem is 
played by the Wronskian $\Wr(f_1, \ldots, f_k)$ (see Corollary~\ref{C6.1}).

\subsection*{Relation to earlier work} The Weyl closure was introduced for the Weyl algebra $A_1$ in the Ph.D. dissertation \cite{Tsai} of H. Tsai, who studied this operation from algorithmic point of view; in particular, he used the spaces $V_\lambda(L)$ to give an algorithm for computing generators of $\cl(L)$. In the end of \cite{Tsai} Tsai observed that there is a relation between his presentation of $\cl(L)$ and the classification of Cannings and Holland \cite{CanningsHolland1994}. Our Theorem~\ref{T1.1} makes this
relation precise and extends it to arbitrary ideals and arbitrary curves. Tsai's work was published in the same year as the paper \cite{BerestWilson00}, where the classification of ideals of $A_1$ in terms of Calogero--Moser spaces was originally established. 
So the invariant $n$ does not appear in \cite{Tsai}. On the other hand, the result of Theorem~\ref{T3.1} does appear in \cite{Tsai} for $A_1$ in a different  form: it is proved by an explicit computation of ranks of two matrices. Our proof is based on the formal continuity of differential operators and works for any curve.
Spaces of conditions were introduced by G. Wilson in his study of bispectral differential operators
\cite{Wilson1993}; the adelic Grassmannian of \cite{Wilson1993, Wilson1998} parametrizes
equivalence classes of primary decomposable subspaces of $\C[x]$, and our `dual adelic
Grassmannian' is its natural counterpart for left ideals. For $X = \A^1$, the Calogero--Moser
index can be read off from the Baker function of the corresponding point of the adelic
Grassmannian (see \cite{Wilson1998}); we use this fact, together with Wilson's Wronskian formula
\cite{Wilson1993}, in our first proof of Theorem~\ref{T1.2}. A second proof, given in the Appendix,
is purely algebraic: it is based on the canonical stairs of an ideal in the sense of Brian\c{c}on
and Maisonobe \cite{BM84}. The codimension characterization of the Calogero--Moser index, on
which both proofs rely, is due to \cite{BerestWilsonNC, BerestWilson09}. Besides \cite{BerestChalykhIdeals}, there is another generalization of the Berest-Wilson construction \cite{BerestWilson00, BerestWilsonNC} to arbitrary curves, 
due to D. Ben-Zvi and T. Nevins  \cite{BenZviNevins}.

\subsection*{Organization of the paper}
Section~\ref{S2} is preliminary: it fixes notation and collects the results on rings of
differential operators, closed ideals, primary decomposable subspaces and Calogero--Moser spaces
needed in the rest of the paper. While most results in this section are known, some are
(apparently) new: in particular, Lemmas~\ref{L2.2} and~\ref{L2.3} on the local duality between
$\Dl$ and the completed local ring, and Proposition~\ref{P2.2} on the invariance of the
Calogero--Moser index under duality. In Section~\ref{S3}, we define the holonomicity index of a
differential operator, establish its basic properties and relate it to the spaces
$V_\lambda(L)$ (Theorem~\ref{T3.1}); at the end of that section we state our main result
(Theorem~\ref{T3.2}, which is a more precise version of Theorem~\ref{T1.3}). Section~\ref{S4}
contains the dual classification of ideals (Theorem~\ref{T4.2}) and its consequences. In
Section~\ref{S5}, we prove the formula for the Calogero--Moser index (Theorem~\ref{T5.1}) and
deduce the main theorem; we refer the reader to the beginning of that section for a summary of
the proof. In Section~\ref{S6}, we consider explicit examples: the case of the affine line, which
plays a special role in this paper, is worked out in detail in Section~\ref{S6.1}, including a
description of the Calogero--Moser data of an ideal in terms of its space of conditions. Finally,
the Appendix contains an alternative, self-contained proof of the local index formula.

The results of this paper are partly based on the second author's Ph.D.\ dissertation
\cite{GoldbergThesis}.

\subsection*{Notation and conventions}
Throughout this paper, we work over the base field $\C$. Unless otherwise specified, $X$ denotes
a smooth affine irreducible curve, $\cO = \cO(X)$ its coordinate ring, $K = K(X)$ its field of
rational functions, and $\cD = \cD(X)$ its ring of differential operators. For $X = \A^1$, we
identify $\cD(\A^1)$ with the first Weyl algebra $A_1 = \C\lan x, \del\ran/(\del x - x \del - 1)$.
An \emph{ideal} of $\cD$ always means a \emph{nonzero} one-sided ideal; we say explicitly whether
it is left or right. We write $P(f)$ for the action of an operator $P$ on a function $f$, as
opposed to the product $P f = P \cdot f$ in the ring $\cD$. For $s \in \C$ and an integer
$k \ge 0$, we write $\ff{s}{k} := s(s-1)\cdots(s-k+1)$, and for a finite set $S$ of non-negative
integers we put
\begin{equation}
\label{E110}
  \nu(S) := \sum_{s \in S} s \, - \binom{|S|}{2} .
\end{equation}
If $S = \{s_0 < s_1 < \ldots < s_{g-1}\}$, then $\nu(S) = \sum_{i} (s_i - i) \ge 0$, with
$\nu(S) = 0$ if and only if $S = \{0, 1, \ldots, g-1\}$.

\section{Preliminaries}
\label{S2}

In this section we fix notation and review the material needed for the rest of the paper.

\subsection{Rings of differential operators}
\label{S2.1}
Let $X$ be a smooth affine irreducible curve over $\C$, with coordinate ring $\cO = \cO(X)$. The
ring $\cD = \cD(X)$ of differential operators on $X$ is a filtered algebra
$\cD = \bigcup_{k \ge 0} \cD_k$, with filtration components defined inductively by
$\cD_{-1} := 0$ and
$\cD_k := \{ P \in \End_{\C}(\cO) : [P, f] \in \cD_{k-1} \text{ for all } f \in \cO\}$. The
elements of $\cD_k$ are differential operators of order $\le k$; we write $\ord(P)$ for the order
of $P$. In particular, $\cD_0 = \cO$ and $\cD_1 = \cO \oplus \Der(\cO)$. As $X$ is smooth, $\cD$ is
generated by $\cO$ and $\Der(\cO)$, and it is a simple Noetherian hereditary domain (see
\cite[Chapter~15]{McConnell}). The associated graded ring
$\bD := \gr \cD = \bigoplus_{k \ge 0} \cD_k/\cD_{k-1}$ is commutative: it is isomorphic to the
coordinate ring of the cotangent bundle $T^*X$, which is a regular domain of dimension $2$.

Let $K = K(X)$ be the field of rational functions on $X$. The ring of rational differential
operators $\cD(K) = K \otimes_{\cO} \cD = \cD \otimes_{\cO} K$ contains $\cD$ as a subring; if
$\delta$ is any nonzero derivation of $K$, then $\cD(K) = K[\delta]$ is an Ore extension of $K$
and hence a (left and right) principal ideal domain. We write $Q$ for the quotient skew-field of
$\cD$. By a \emph{fractional left} (resp., \emph{right}) \emph{ideal} we mean a nonzero finitely
generated left (resp., right) $\cD$-submodule of $Q$. For a fractional left ideal $M$, its
\emph{dual} is the fractional right ideal
\begin{equation}
\label{E21}
  M^* := \{ q \in Q \, : \, M q \subseteq \cD \} ,
\end{equation}
and similarly for right ideals. Since $\cD$ is hereditary, every fractional ideal $M$ is a
projective module of rank $1$; we have $M^* \cong \Hom_{\cD}(M, \cD)$ and $M^{**} = M$ as
submodules of $Q$. We write $\cL(\cD)$ (resp., $\cR(\cD)$) for the set of isomorphism classes of
left (resp., right) ideals of $\cD$; each class in $\cL(\cD)$ contains fractional ideals, and the
assignment $M \mapsto M^*$ induces a bijection $\cL(\cD) \cong \cR(\cD)$.

For subspaces $U, V \subseteq K$, we put
\begin{equation}
\label{E22}
  \cD(U, V) := \{ P \in \cD(K) \, : \, P(U) \subseteq V \} .
\end{equation}
For example, $\cD(\cO, \cO) = \cD$, and $\cD(\cO, V)$ is a right ideal of $\cD$ for any
$V \subseteq \cO$.

\subsubsection*{Local structure}
For a point $\lambda \in X$, we write $\ml \subset \cO$ for the corresponding maximal ideal,
$\cO_\lambda$ for the local ring of $X$ at $\lambda$, and $\nl := \ml \cO_\lambda$ for its maximal
ideal. Since $\cO_\lambda$ is a discrete valuation ring, $\nl = t\,\cO_\lambda$ for some
$t \in \cO_\lambda$ (a \emph{local parameter} at $\lambda$), and every nonzero ideal of
$\cO_\lambda$ is of the form $\nl^k$; we write $\ord_\lambda: \cO_{\lambda} \to \Z \cup \{\infty\}\,$ for the corresponding valuation. The
module $\Der(\cO_\lambda)$ is free of rank $1$ over $\cO_{\lambda}$, generated by the derivation $\del = \del_t$ such
that $\del_t(t) = 1$ (see \cite[Theorem~II.8.8]{Hartshorne}). The ring
$\cD_\lambda := \cD(\cO_\lambda) = \cO_\lambda \otimes_{\cO} \cD$ is then generated by $\cO_\lambda$ and
$\del$ (see \cite[Corollary~15.5.6]{McConnell}):
\begin{equation}
\label{E23}
  \cD_\lambda \cong \cO_\lambda\lan \del \ran \big/ \bigl(\del f - f \del = \del(f) \, , \
  f \in \cO_\lambda \bigr) ,
\end{equation}
so that every element of $\cD_\lambda$ can be uniquely written as $\sum_i a_i \del^i$ with
$a_i \in \cO_\lambda$. We identify $\cD$ with a subring of $\cD_\lambda$ via the localization map. We also need the completion of $\cO_\lambda$ at $\nl$:
$$
\hO_\lambda := \varprojlim \cO_\lambda/\nl^{j} \cong \C[[t]]\,,
$$ 
with maximal ideal $\hn_\lambda = t\,\hO_\lambda$, and the ring
$\hD_\lambda := \hO_\lambda \otimes_{\cO_\lambda} \cD_\lambda \cong \C[[t]][\del]$ of {\it continuous} differential
operators on $\hO_\lambda$ with respect to the $\hn_\lambda$-adic topology. The natural multiplication map $\hO_\lambda \otimes_{\cO_\lambda} \cD_\lambda \xrightarrow{\sim} \cD_\lambda
\otimes_{\cO_\lambda} \hO_\lambda$ is an isomorphism of filtered bimodules, so that $\hD_\lambda \cong \cD_\lambda \otimes_{\cO_\lambda} \hO_\lambda $ (see
\cite[Lemma~2.3]{knop06}). Thus we have ring embeddings
$\cD \subset \cD_\lambda \subset \hD_\lambda$, extending $\cO \subset \cO_\lambda \subset
\hO_\lambda$. Finally, we introduce the notation $\theta := t \del_t$ for the Euler operator with respect to the local parameter $t$ and note the identity
\begin{equation}
\label{E24}
 t^k \del^k  = \ff{\theta}{k} :=   \theta(\theta - 1) \cdots (\theta - k + 1) , \qquad \forall\,k \ge 1 \,.
\end{equation}
(The latter holds in $\cD_\lambda$ since both sides are operators of order $k$ with leading coefficient $t^k$ that annihilate $1, t, \ldots, t^{k-1}$.)

\subsubsection*{Formal continuity}
We recall a general property of algebraic differential operators. Let $A$ be a commutative
$\C$-algebra. For $A$-modules $M$ and $N$, the space $\cD_A(M,N) = \bigcup_{i \ge 0} \cD^i_A(M,N)$
of differential operators from $M$ to $N$ is defined inductively by $\cD^{-1}_A(M,N) := 0$ and
$\cD^i_A(M,N) := \{ \theta \in \Hom_{\C}(M,N) : [a, \theta] \in \cD^{i-1}_A(M,N) \text{ for all }
a \in A \}$. The following properties are well known. We give  proofs for completeness.

\begin{proposition}
\label{P2.1}
Let $A$, $M$ and $N$ be as above.
\begin{enumerate}
\item[$(a)$] For any ideal $\mathfrak{a} \subseteq A$ and all $i, r \ge 0$, we have
$\ \cD^i_A(M,N)\,[\mathfrak{a}^{i+r} M] \subseteq \mathfrak{a}^r N$.
\item[$(b)$] Let $\fm \subset A$ be a maximal ideal with $A/\fm \cong \C$. Then, 
the assignment $\theta \mapsto [\, a + \fm^{i+1} \mapsto \theta(a)\,]$ is a linear isomorphism
$\, \cD^i_A(A, A/\fm) \stackrel{\sim}{\to} \Hom_{\C}(A/\fm^{i+1}, A/\fm)\,$ for all $i \ge 0$.
\end{enumerate}
\end{proposition}

\begin{proof}
$(a)$ We use induction on $i$. For $i = 0$, the claim is obvious, since $\cD^0_A(M,N) =
\Hom_A(M,N)$. If the claim holds for $i$ (and all $r$), then, by definition of $\cD^{i+1}_A$,
\[
  \cD^{i+1}_A(M,N)\,[\mathfrak{a}^{i+1+r} M] \subseteq \cD^{i}_A(M,N)\,[\mathfrak{a}^{i+r} M] +
  \mathfrak{a}\, \cD^{i+1}_A(M,N)\,[\mathfrak{a}^{i+r} M] ,
\]
and the right-hand side is contained in $\mathfrak{a}^r N$ by the induction hypothesis and a
second induction on $r$. $(b)$ By $(a)$, every $\theta \in \cD^i_A(A, A/\fm)$ vanishes on
$\fm^{i+1}$, so the map is well defined and injective. Conversely, if $\theta \in \Hom_{\C}(A,
A/\fm)$ vanishes on $\fm^{i+1}$, then, since $A = \C \oplus \fm$, the $(i+1)$-fold commutator
$[a_0, [a_1, \ldots [a_i, \theta] \ldots ]]$ with $a_j \in A$ may be computed with $a_j \in \fm$,
in which case it takes values in $\fm\, \theta(A) + \theta(\fm^{i+1}) = 0$. Hence
$\theta \in \cD^i_A(A, A/\fm)$.
\end{proof}

\begin{lemma}[\cite{SmithStafford}, 1.3(e),(f)]
\label{L2.1}
Let $A$ be a localization of $\cO(X)$. Then, for any ideal $\mathfrak{a} \subseteq A$, the natural
map $\cD(A)/\mathfrak{a}\, \cD(A) \to \cD_A(A, A/\mathfrak{a})$ is an isomorphism. In particular,
$\cD(A, \mathfrak{a}) = \mathfrak{a}\,\cD(A)$.
\end{lemma}

\subsection{Local duality}
\label{S2.2}
Fix $\lambda \in X$. The right $\cD$-module
\begin{equation}
\label{E25}
  \Dl := \cD/\ml \cD
\end{equation}
plays a central role in this paper. It carries a filtration induced by the order filtration on
$\cD$: we write $\Dl^{\le j}$ for the image of $\cD_j$ in $\Dl$, and define the \emph{degree}
$\deg(F)$ of a nonzero element $F \in \Dl$ to be the minimal $j$ such that $F \in \Dl^{\le j}$. For
a subspace $V \subseteq \Dl$, we put $V^{\le j} := V \cap \Dl^{\le j}$. The next lemma describes
$\Dl$ in local coordinates and identifies its linear dual with the completed local ring.

\begin{lemma}
\label{L2.2}
Let $t$ be a local parameter at $\lambda$, and let $\del = \del_t$.
\begin{enumerate}
\item[$(a)$] The natural maps $\ \Dl \xrightarrow{\sim} \cD_\lambda/\nl \cD_\lambda \xrightarrow{\sim} \hD_\lambda/ t\,\hD_\lambda$
are isomorphisms of filtered vector spaces, and the classes of $\,1, \del, \del^2, \ldots$ form a
basis of $\Dl$; thus $\Dl \cong \C[\del]$, with $\deg$ corresponding to the degree of polynomials.
In particular, $\Dl$ is naturally a right $\hD_\lambda$-module.
\item[$(b)$] The pairing $\lan F, f \ran_\lambda := F(f)\ (\mathrm{mod}\ \hn_\lambda)$ is well
defined on $\Dl \times \hO_\lambda$ and extends \eqref{E16}. It satisfies
$\lan \del^i, t^j \ran_\lambda = i!\, \delta_{ij}$, and therefore identifies $\hO_\lambda$ with the
full linear dual $\Hom_{\C}(\Dl, \C)$, and $\hO_\lambda/\hn_\lambda^{j+1} \cong
\cO_\lambda/\nl^{j+1}$ with the dual of $\Dl^{\le j}$, for each $j \ge 0$.
\item[$(c)$] For all $F \in \Dl$, $P \in \hD_\lambda$ and $f \in \hO_\lambda$, we have
$\lan F P, f \ran_\lambda = \lan F, P(f) \ran_\lambda$.
\item[$(d)$] For $P \in \hD_\lambda$, let $R_P : \Dl \to \Dl$, $F \mapsto F P$, be the right
multiplication by $P$. Then
\[
  (\Ker R_P)^{\perp} = P(\hO_\lambda) \qquad \text{and} \qquad (\im R_P)^\perp = \Ker(P;
  \hO_\lambda) ,
\]
where, for $V \subseteq \Dl$, we write 
$$ 
V^\perp := \{ f \in \hO_\lambda : \lan F, f \ran_\lambda =
0 \text{ for all } F \in V\}\,.
$$
Consequently, there are natural isomorphisms
$(\Ker R_P)^* \cong \Coker(P; \hO_\lambda)$ and $(\Coker R_P)^* \cong \Ker(P; \hO_\lambda)$.
\end{enumerate}
\end{lemma}

\begin{proof}
$(a)$ Since $\cD_\lambda = \cO_\lambda \otimes_{\cO} \cD$ and $\hD_\lambda = \hO_\lambda
\otimes_{\cO} \cD$ as left modules, we have $\cD_\lambda/\nl \cD_\lambda \cong (\cO_\lambda/\nl)
\otimes_{\cO} \cD \cong (\cO/\ml) \otimes_{\cO} \cD = \Dl$, and similarly for $\hD_\lambda$. The
rest follows from \eqref{E23}. $(b)$ If $F \in t\, \hD_\lambda$, then $F(f) \in \hn_\lambda$ for
all $f$, so the pairing is well defined; the formula $\del^i(t^j)|_{t=0} = i!\,\delta_{ij}$ is
obvious. $(c)$ is immediate from the definition. $(d)$ By $(b)$ and $(c)$, the linear map
$P: \hO_\lambda \to \hO_\lambda$ is the transpose of $R_P$. The claim follows then from the
standard facts of linear algebra: if $A: U \to U'$ is a linear map of vector spaces and
$A^{\mathsf{T}} : (U')^* \to U^*$ is its transpose, then $\im(A^{\mathsf{T}})$ is the annihilator
of $\Ker(A)$, and $\Ker(A^{\mathsf{T}})$ is the annihilator of $\im(A)$.
\end{proof}

\begin{remark}
In general, the map $R_P$ in Lemma~\ref{L2.2}$(d)$ is neither injective nor surjective. For
example, if $X = \A^1$, $\lambda = 0$ and $P = x \del - n$ with $n \in \Z_{\ge 0}$, then
$\del^j P \equiv (j - n)\, \del^j$ modulo $x A_1$, so that $\Ker(R_P) = \C\,\del^n$ and
$\im(R_P)$ is spanned by $\del^j$ with $j \ne n$. Dually, $P(x^j) = (j - n)\, x^j$, so that both
the kernel and the cokernel of $P$ on $\C[[x]]$ are spanned by (the class of) $x^n$.
\end{remark}

Given a finite-dimensional subspace $V \subseteq \Dl$, we define its \emph{degree set}
\[
  \Sigma(V) := \{ \deg(F) \, : \, 0 \ne F \in V \} \subset \Z_{\ge 0} ,
\]
and put $\nu(V) := \nu(\Sigma(V))$, see \eqref{E110}. Dually, for a subspace
$U \subseteq \hO_\lambda$, we define its \emph{set of gaps}
$\Gamma(U) := \Z_{\ge 0} \setminus \{ \ord_\lambda(f) : 0 \ne f \in U \}$.

\begin{lemma}
\label{L2.3}
Let $V$ be a finite-dimensional subspace of $\Dl$.
\begin{enumerate}
\item[$(a)$] $V$ has a basis consisting of elements of pairwise distinct degrees; for any such
basis $\{F_1, \ldots, F_g\}$, we have $\Sigma(V) = \{\deg F_1, \ldots, \deg F_g\}$. In particular,
$|\Sigma(V)| = \dim V$ and $\dim V^{\le j} = |\Sigma(V) \cap [0, j]|$ for all $j \ge 0$.
\item[$(b)$] $(V^\perp)^\perp = V\,$ and $\ \Gamma(V^\perp) = \Sigma(V)$.
\item[$(c)$] Let $u \in \hO_\lambda$ be an element of order $k \ge 0$, and let $V' \subseteq \Dl$ be
the (finite-dimensional) subspace such that $(V')^\perp = u\, V^\perp$. Then
\[
  \Sigma(V') = \{0, 1, \ldots, k-1\} \, \sqcup\, \bigl(k + \Sigma(V)\bigr) \qquad \text{and} \qquad
  \nu(V') = \nu(V) .
\]
\end{enumerate}
\end{lemma}

\begin{proof}
$(a)$ is standard (use Gaussian elimination with respect to the basis $\{\del^i\}$). $(b)$ The
first equality holds because $\hO_\lambda$ is the full dual of $\Dl$ and $\dim V < \infty$. For the
second, note that $(V^\perp + \hn_\lambda^{j+1})^\perp = V \cap \Dl^{\le j}$ by
Lemma~\ref{L2.2}$(b)$; hence $\dim \hO_\lambda/(V^\perp + \hn^{j+1}_\lambda) = \dim V^{\le j}$
for all $j \ge 0$. On the other hand, by the modular law,
\[
  \frac{V^\perp + \hn_\lambda^{j}}{V^\perp + \hn_\lambda^{j+1}} \, \cong \,
  \frac{\hn^j_\lambda}{(V^\perp \cap \hn^j_\lambda) + \hn^{j+1}_\lambda} ,
\]
and the right-hand side is zero if $V^\perp$ contains an element of order $j$ and is
one-dimensional otherwise. Thus $j \in \Gamma(V^\perp)$ if and only if $\dim V^{\le j} = \dim
V^{\le j-1} + 1$, i.e.\ if and only if $j \in \Sigma(V)$, by $(a)$. $(c)$ The subspace
$u V^\perp$ is closed and of finite codimension in $\hO_\lambda$, so $V' := (u V^\perp)^\perp$ is
finite-dimensional and $(V')^\perp = u V^\perp$. Clearly,
$\Gamma(u\, V^\perp) = \{0, \ldots, k-1\} \sqcup (k + \Gamma(V^\perp))$, which, in view of $(b)$,
gives the first formula. If $|\Sigma(V)| = g$, then
\[
  \nu(V') = \binom{k}{2} + \sum_{s \in \Sigma(V)} (s + k) - \binom{k+g}{2} =
  \sum_{s \in \Sigma(V)} s - \binom{g}{2} = \nu(V) ,
\]
since $\binom{k+g}{2} = \binom{k}{2} + \binom{g}{2} + k g$.
\end{proof}

\subsection{Closed ideals}
\label{S2.3}
The \emph{closure} of a left ideal $I \subseteq \cD$ is the left ideal
\begin{equation}
\label{E26}
  \cl(I) = \cl_{\cD}(I) := \cD(K)\, I \cap \cD ,
\end{equation}
and $I$ is called \emph{closed} if $\cl(I) = I$. For $L \in \cD$, we write $\cl(L) := \cl(\cD L)$.
For $X = \A^1$, this operation was introduced and studied by Tsai \cite{Tsai} under the name
\emph{Weyl closure}. Clearly, $I \subseteq \cl(I)$ and $\cl(\cl(I)) = \cl(I)$. By an operator of
\emph{minimal order} in $I$ we mean any nonzero $L \in I$ with $\ord(L) = \min\{\ord(P) : 0 \ne
P \in I\}$.

\begin{lemma}
\label{L2.4}
\begin{enumerate}
\item[$(a)$] Let $N$ be a left $\cD(K)$-module (for example, the space of meromorphic functions on
$X$ or on an open subset of $X$ in the classical topology). Then, for any finite subset
$\{u_1, \ldots, u_k\} \subset N$, the left ideal
$\Ann_{\cD}(u_1, \ldots, u_k) := \{P \in \cD : P u_i = 0 \text{ for all } i\}$ is closed.
\item[$(b)$] Let $I$ be a nonzero left ideal of $\cD$, and let $L \in I$ be an operator of minimal
order. Then $\cD(K)\, I = \cD(K)\, L$, and therefore
\[
  \cD L \subseteq I \subseteq \cl(I) = \cl(L) .
\]
In particular, if $I$ is closed, then $I = \cl(L)$. Any other operator of minimal order in $I$
has the form $h L$ with $h \in K^\times$.
\item[$(c)$] For any nonzero $L \in \cD$, the quotient $\cl(L)/\cD L$ is the $\cO$-torsion submodule
of the left $\cD$-module $\cD/\cD L$.
\end{enumerate}
\end{lemma}

\begin{proof}
$(a)$ If $P \in \cD(K)\, \Ann_{\cD}(u_1, \ldots, u_k) \cap \cD$, then
$P = \sum_j q_j b_j$ with $q_j \in \cD(K)$ and $b_j u_i = 0$ for all $i$ and $j$; hence
$P u_i = 0$. $(b)$ Since $\cD(K)$ is a left principal ideal domain, $\cD(K)\, I = \cD(K)\, L'$
for some $L' \in \cD(K)$. As $\cD(K)\, I = K \cdot I$, we have $f L' \in I$ for some nonzero
$f \in \cO$. On the other hand, $L = q L'$ with $q \in \cD(K)$, so that
$\ord(L) \ge \ord(L') = \ord(f L')$. By minimality, $\ord(L) = \ord(L')$; hence $q \in K^\times$
and $\cD(K)\, L = \cD(K)\, L'$. The same argument proves the last statement. $(c)$ An operator
$P \in \cD$ belongs to $\cl(L)$ if and only if $P = q L$ for some $q \in \cD(K) = K \otimes_{\cO}
\cD$, i.e.\ if and only if $f P \in \cD L$ for some nonzero $f \in \cO$.
\end{proof}

For a left $\cD$-module $N$ and a point $\lambda \in X$, we write
$H^0_\lambda(N) := \{ u \in N : \ml^k u = 0 \text{ for some } k \ge 1\}$ for the submodule of
elements supported at $\lambda$, and $N^{\ml} := \{u \in N : \ml u = 0\}$. If $N$ is a torsion
$\cO$-module, then $N = \bigoplus_{\lambda \in X} H^0_\lambda(N)$, since $\cO$ is a Dedekind
domain. We will need the following elementary special case of Kashiwara's theorem.

\begin{lemma}
\label{L2.5}
For any left $\cD$-module $N$ and any $\lambda \in X$, the multiplication map
$(\cD/\cD \ml) \otimes_{\C} N^{\ml} \to H^0_\lambda(N)$ is an isomorphism. In particular,
$H^0_\lambda(N) = \cD \cdot N^{\ml}$, and $N^{\ml}$ is finite-dimensional whenever $H^0_\lambda(N)$
is a finitely generated $\cD$-module.
\end{lemma}

\begin{proof}
Every element of $\cO \setminus \ml$ acts invertibly on $H^0_\lambda(N)$, so $H^0_\lambda(N)$ is
naturally a module over $\cD_\lambda$, and $\cD/\cD \ml \cong \cD_\lambda/\cD_\lambda t \cong
\C[\del]$. If $u \in N^{\ml}$, then $t u = 0$, and $t\, \del^j u = - j\, \del^{j-1} u$ for all
$j \ge 1$; hence $t^j \del^j u = (-1)^j j!\, u$. This implies injectivity: if $\sum_{j=0}^r \del^j
u_j = 0$ with $u_j \in N^{\ml}$, then applying $t^r$ gives $u_r = 0$. For surjectivity, we show by
induction on $k$ that $t^k v = 0$ implies $v \in \sum_j \del^j N^{\ml}$. For $k = 1$, there is
nothing to prove. If $k \ge 2$, put $u := t^{k-1} v \in N^{\ml}$ and
$v' := v - \frac{(-1)^{k-1}}{(k-1)!}\, \del^{k-1} u$. Then $t^{k-1} v' = u - u = 0$, and the claim
follows by induction. The last statement holds because $\cD/\cD \ml$ is a simple $\cD$-module.
\end{proof}

\subsection{Primary decomposable subspaces and spaces of conditions}
\label{S2.4}
We recall the classification of ideals of $\cD$ due to Cannings and Holland
\cite{CanningsHolland1994}. A subspace $W \subseteq \cO$ is called \emph{$\ml$-primary} if it
contains a power of the maximal ideal $\ml$, and \emph{primary decomposable} if it is a finite
intersection $W = \bigcap_{\lambda} W_\lambda$ of $\ml$-primary subspaces $W_\lambda$ (the
subspaces $W_\lambda$ are then uniquely determined by $W$, see
\cite[Theorem~2.4]{CanningsHolland1994}). For example, every nonzero ideal of $\cO$ is primary
decomposable. Two primary decomposable subspaces $W$ and $W'$ are called \emph{equivalent} if
$W' = h\,W$ for some $h \in K^\times$. A right ideal $M \subseteq \cD$ is called \emph{fat} if
$M \cap \cO \ne 0$ (here, as usual, we identify $\cO$ with the subring $\cD_0$ of differential
operators of order zero in $\cD$). For a right ideal $M$, we write $M \cdot \cO$ for the subspace of $\cO$ spanned
by all $P(f)$ with $P \in M$ and $f \in \cO$.

\begin{theorem}[\cite{CanningsHolland1994}]
\label{T2.1}
\begin{enumerate}
\item[$(a)$] The assignments $W \mapsto \cD(\cO, W)$ and $M \mapsto M \cdot \cO$ are mutually inverse
bijections between the set of primary decomposable subspaces of $\cO$ and the set of fat right
ideals of $\cD$. These bijections preserve inclusions, sums and intersections.
\item[$(b)$] Every right ideal of $\cD$ is isomorphic to a fat one, and
$\cD(\cO, W) \cong \cD(\cO, W')$ if and only if $W$ and $W'$ are equivalent. Thus
$W \mapsto [\cD(\cO, W)]$ induces a bijection between the set of equivalence classes of primary
decomposable subspaces of $\cO$ and $\cR(\cD)$.
\item[$(c)$] For any primary decomposable subspace $W$, we have $\cD(\cO, W)^* = \cD(W, \cO)$.
Consequently, every left ideal of $\cD$ is isomorphic to $\cD(W, \cO)$ for a primary decomposable
subspace $W \subseteq \cO$, which is uniquely determined up to equivalence.
\end{enumerate}
\end{theorem}

Parts $(a)$ and $(b)$ are proved in \cite{CanningsHolland1994} (see also
\cite[Section~3]{BerestWilson00} for the case $X = \A^1$). Part $(c)$ is an easy consequence of $(a)$ (cf.\
\cite[Lemma~5.3]{BerestWilson00}): if $q \in \cD(W, \cO)$ and $p \in \cD(\cO, W)$, then
$q p \in \cD(\cO, \cO) = \cD$; conversely, if $q \in Q$ and $q\, \cD(\cO, W) \subseteq \cD$, then
$q \in \cD(K)$, because $\cD(\cO, W)$ contains a nonzero function, and
$q(W) = q(\cD(\cO, W) \cdot \cO) \subseteq \cD(\cO) = \cO$.

Following \cite{Wilson1993}, we now describe primary decomposable subspaces in terms of linear
conditions. Consider the right $\cD$-module
\[
  \Delta(X) := \bigoplus_{\lambda \in X} \Dl , \qquad \Dl = \cD/\ml \cD ,
\]
(cf.\ \cite[Section~3]{BK98}). We call a subspace $C \subseteq \Delta(X)$ \emph{homogeneous} if
$C = \bigoplus_{\lambda \in X} C_\lambda$, where $C_\lambda := C \cap \Dl$.

\begin{definition}
\label{D2}
A \emph{space of conditions} on $X$ is a finite-dimensional homogeneous subspace $C$ of
$\Delta(X)$. The finite set $\Supp(C) := \{\lambda \in X : C_\lambda \ne 0\}$ is called the
\emph{support} of $C$. To a space of conditions $C$ we associate the subspace $W_C \subseteq \cO$
defined by \eqref{E17}, and the right ideal
\begin{equation}
\label{E27}
  \Ann_{\cD}(C) := \{ P \in \cD \, : \, F P = 0 \ \text{in}\ \Delta(X) \ \text{for all } F \in C \}.
\end{equation}
\end{definition}

In local coordinates, an element of $C_\lambda$ is a linear functional on $\cO$ of the form
$f \mapsto [(a_k \del_t^k + \ldots + a_1 \del_t + a_0)\, f](\lambda)$ with $a_i \in \C$. In the
case $X = \A^1$, Definition~\ref{D2} agrees with that of \cite[Section~4]{Wilson1993}.

\begin{lemma}
\label{L2.6}
\begin{enumerate}
\item[$(a)$] For any space of conditions $C$, the subspace $W_C = \bigcap_\lambda W_{C_\lambda}$ is
primary decomposable. The closure of $W_C$ in $\hO_\lambda$ is $C_\lambda^\perp$, and the pairing
\eqref{E16} induces isomorphisms $\,\cO/W_{C_\lambda} \cong C_\lambda^*\,$ and
$\,\cO/W_C \cong C^*$; in particular, $\codim_{\cO}(W_C) = \dim C$.
\item[$(b)$] The assignment $C \mapsto W_C$ is a bijection between the set of spaces of conditions on
$X$ and the set of primary decomposable subspaces of $\cO$.
\item[$(c)$] $\Ann_{\cD}(C) = \cD(\cO, W_C)$; in particular, $\Ann_{\cD}(C)$ is a fat right ideal.
\end{enumerate}
\end{lemma}

\begin{proof}
$(a)$ Choose $N$ such that $C_\lambda \subseteq \Dl^{\le N-1}$ for all $\lambda$. By
Lemma~\ref{L2.2}$(b)$, $C_\lambda^\perp \supseteq \hn_\lambda^N$, so that
$W_{C_\lambda} = \cO \cap C_\lambda^\perp$ contains $\cO \cap \hn^N_\lambda = \ml^N$. Since the
natural map $\cO \to \hO_\lambda/\hn_\lambda^N$ is surjective, $\cO/W_{C_\lambda} \cong
\hO_\lambda/C_\lambda^\perp \cong C_\lambda^*$. The rest follows from the Chinese Remainder Theorem.
$(b)$ If $W = \bigcap_\lambda W_\lambda$ with $W_\lambda \supseteq \ml^N$, put
$C_\lambda := \{ F \in \Dl : \lan F, W_\lambda\ran_\lambda = 0\} \subseteq \Dl^{\le N-1}$. Then
$W_{C_\lambda} = W_\lambda$, by the duality between $\Dl^{\le N-1}$ and $\cO/\ml^N$. Conversely,
$C$ is recovered from $W_C$ by $(a)$ and Lemma~\ref{L2.3}$(b)$. $(c)$ Let $P \in \cD$ and
$F \in C_\lambda$. By Lemma~\ref{L2.1}, $F P = 0$ in $\Dl$ if and only if $(FP)(\cO) \subseteq \ml$,
i.e.\ if and only if $\lan F, P(f) \ran_\lambda = 0$ for all $f \in \cO$. Hence
$P \in \Ann_{\cD}(C)$ if and only if $P(\cO) \subseteq W_C$.
\end{proof}

\subsection{Calogero--Moser spaces and the Calogero--Moser index}
\label{S2.5}
We now recall the classification of ideals of $\cD$ in terms of Calogero--Moser spaces, referring
the reader to \cite{BerestChalykhIdeals} for details. As mentioned in the Introduction, there is a
natural map
\begin{equation}
\label{E28}
  \gamma : \ \cL(\cD) \to K_0(\cD) \cong K_0(X) \xrightarrow{\ \det\ } \Pic(X) ,
\end{equation}
whose fibres are the stable isomorphism classes of ideals (see \cite{BerestWilson09} and
\cite[Section~3.2]{BerestChalykhIdeals}). Given a line bundle $\mathcal{I}$ on $X$, we denote by
$ B := \cO[\mathcal{I}]$ the {\it one-point extension} of $ \cO $ by $\mathcal{I}$, the (noncommutative) square-zero extension of the coordinate ring $\tilde{\cO} = \cO(\tilde X) \cong \cO \times \C $ of the curve $ \tilde{X} = X \cup \{\infty\} $ by the $\cO$-module $\mathcal{I}\,$, and let
$$
\Pi^{\boldsymbol{\lambda}}(B)\,:=\, T_B {\rm Der}(B,\,B^{\otimes 2})/\langle \Delta_B - \boldsymbol{\lambda}\rangle
$$
be the (deformed) preprojective algebra of $B$ of weight
$\boldsymbol{\lambda} = (1, -n)$. The \emph{$n$-th Calogero--Moser space of type
$(X, \mathcal{I})$} is then defined to be the moduli space of representations
(the representation variety)
$$
\cC_n(X, \mathcal{I}) := {\rm Rep}_{\boldsymbol{n}}[\Pi^{\boldsymbol{\lambda}}(B)]/\!/{\rm GL}_{\boldsymbol{n}}(\C)
$$
parametrizing the (isomorphism classes of) finite-dimensional $\Pi^{\boldsymbol{\lambda}}(B)$-modules of dimension vector $\boldsymbol{n} = (n,1)$ (see \cite[Section~3.3]{BerestChalykhIdeals}). 
For all pairs $ (X, \mathcal{I}) $ and each $ n\ge 0 $, $\cC_n(X, \mathcal{I})$ is a smooth affine irreducible variety of (complex) dimension $2n$ (see \cite[Theorem 1.1]{BerestChalykhIdeals}).
We write $\overline{\cC}_n(X, \mathcal{I})$ for the
quotient of $\cC_n(X, \mathcal{I})$ by the natural action of $\Aut_X(\mathcal{I})$; this quotient
depends only on the class of $\mathcal{I}$ in $\Pic(X)$. 
\begin{theorem}[\cite{BerestChalykhIdeals}, Theorem 4.2]
\label{T2.2}
For each $n \ge 0$ and $[\mathcal{I}] \in \Pic(X)$, there are canonical injective maps
$\,\omega_n : \overline{\cC}_n(X, \mathcal{I}) \hookrightarrow \gamma^{-1}[\mathcal{I}]\,$ that assemble $($over all $n$$)$ to a bijection:
\[
  \omega : \ \bigsqcup_{n \ge 0}\, \overline{\cC}_n(X, \mathcal{I})
  \stackrel{\sim}{\longrightarrow} \gamma^{-1}[\mathcal{I}] \,.
\]
Thus, we have a natural bijective correspondence $ \,\bigsqcup_{[{\mathcal I}]\in \Pic(X)}\,\bigsqcup_{n \ge 0}\, \overline{\cC}_n(X, \mathcal{I})\,\cong\,\cL(\cD)\,$.
\end{theorem}
Theorem~\ref{T2.2} provides a classification of left ideals of $\cD $ in terms of finite-dimensional representations of a finitely-generated algebra, which, in practice, reduces
to an explicit description of such ideals in terms of linear algebra data.
For $X = \A^1$, the spaces $\cC_n := \cC_n(\A^1, \cO)$ are the classical Calogero--Moser spaces 
introduced in \cite{Wilson1998}: they have the following simple matrix presentation
\begin{equation}
\label{E29}
  \cC_n = \{ (\mathsf{X}, \mathsf{Y}, v, w) \, : \, [\mathsf{X}, \mathsf{Y}] + \mathrm{Id}_n = v\,w
  \} / \GL_n(\C) ,
\end{equation}
where $\mathsf{X}, \mathsf{Y} \in \End(\C^n)$, $v \in \Hom(\C, \C^n)$, $w \in \Hom(\C^n, \C)$, and
$\GL_n(\C)$ acts by $g \cdot (\mathsf{X}, \mathsf{Y}, v, w) = (g \mathsf{X} g^{-1}, g \mathsf{Y}
g^{-1}, g v, w g^{-1})$. In this case, Theorem~\ref{T2.2} is the main result of
\cite{BerestWilson00, BerestWilsonNC}.

\begin{definition}
\label{DCM}
The \emph{Calogero--Moser index} of a left ideal $M \subseteq \cD$ is the integer $n = n(M) \ge 0$
such that the isomorphism class of $M$ belongs to the image of $\omega_n$ in Theorem~\ref{T2.2}.
\end{definition}

By definition, $n(M)$ depends only on the isomorphism class of $M$; in particular, it is defined
for fractional left ideals. The next theorem gives a characterization of $n(M)$ that does not
refer to Calogero--Moser spaces. If $M$ is a left or right ideal of $\cD$, we equip it with the
induced filtration $M_k := M \cap \cD_k$, and write $\bM := \gr M \subseteq \bD$ for the associated
graded ideal. For a finitely generated $\bD$-module $E$, we write
$E^\vee := \Hom_{\bD}(E, \bD)$.

\begin{theorem}[\cite{BerestWilson09}, Proposition~2.1; \cite{BerestChalykhIdeals}, Theorem~3.1]
\label{T2.3}
Let $M$ be a left ideal of $\cD$. Then
\begin{enumerate}
\item[$(a)$] $\bM^{\vee\vee} \cong \bD \otimes_{\cO} \mathcal{I}_M\,$ for a unique (up to
isomorphism) line bundle $\mathcal{I}_M$ on $X$, and $\gamma[M] = [\mathcal{I}_M]$;
\item[$(b)$] the canonical map $\bM \to \bM^{\vee\vee}$ is injective, and its cokernel
$\,\mathsf{V}_M := \bM^{\vee\vee}/\bM$ is a finite-dimensional vector space, with
\begin{equation}
\label{E210}
  \dim_{\C} \mathsf{V}_M = n(M) ;
\end{equation}
\item[$(c)$] more generally, \eqref{E210} holds if $M$ is any rank one projective $\cD$-module
equipped with a good filtration such that $\gr M$ is torsion-free.
\end{enumerate}
\end{theorem}

The statements $(a)$ and $(b)$ of Theorem~\ref{T2.3}, as well as the fact that
$\dim \mathsf{V}_M$ does not depend on the choice of filtration, hold for right ideals as well (the proofs in
\cite{BerestWilson09} are left-right symmetric). For a right ideal $M$, we \emph{define} its
Calogero--Moser index $n(M)$ by formula \eqref{E210}. We will show (see
Proposition~\ref{P2.2} below) that this is consistent with the duality $M \mapsto M^*$. First, we
give a homological interpretation of \eqref{E210}.

\begin{lemma}
\label{L2.7}
For any (left or right) ideal $M$ of $\cD$, there are natural isomorphisms
\[
  \Ext^1_{\bD}(\bM, \bD) \cong \Ext^2_{\bD}(\mathsf{V}_M, \bD) \cong \mathsf{V}_M^* .
\]
Hence $n(M) = \dim_{\C} \Ext^1_{\bD}(\bM, \bD)$.
\end{lemma}

\begin{proof}
Dualizing the exact sequence $0 \to \bM \to \bD \mathcal{I}_M \to \mathsf{V}_M \to 0$, we get
\[
  \Ext^1_{\bD}(\bD \mathcal{I}_M, \bD) \to \Ext^1_{\bD}(\bM, \bD) \to
  \Ext^2_{\bD}(\mathsf{V}_M, \bD) \to \Ext^2_{\bD}(\bD \mathcal{I}_M, \bD) ,
\]
where the first and the last terms vanish, because $\bD \mathcal{I}_M$ is a projective
$\bD$-module. This proves the first isomorphism. For the second, note that $\mathsf{V}_M$ is a
$\bD$-module of finite length, supported at finitely many (closed) points $p_i \in T^*X$; hence
$\mathsf{V}_M = \bigoplus_i \mathsf{V}_{p_i}$ and $\Ext^2_{\bD}(\mathsf{V}_M, \bD) \cong
\bigoplus_i \Ext^2_{\bD_{p_i}}(\mathsf{V}_{p_i}, \bD_{p_i})$. Each local ring $\bD_{p_i}$ is
regular of dimension $2$, and hence Gorenstein (see \cite[Proposition~3.1.20]{BrunsHerzog98}). By
local duality, $\Ext^2_{\bD_{p_i}}(\mathsf{V}_{p_i}, \bD_{p_i}) \cong \mathsf{V}_{p_i}^*$ for
each $i$, and the result follows.
\end{proof}

\begin{proposition}
\label{P2.2}
For any left ideal $M$ of $\cD$, we have $n(M^*) = n(M)$, where $M^*$ is the dual right ideal
\eqref{E21}.
\end{proposition}

\begin{proof}
Let $\widetilde{\cD} := \bigoplus_{k \ge 0} \cD_k \tau^k \subseteq \cD[\tau]$ be the Rees algebra
of $\cD$, and let $\widetilde{M} := \bigoplus_k M_k \tau^k$ be the Rees module of $M$; thus
$\tau$ is a central non-zero-divisor of degree $1$, with $\widetilde{\cD}/\tau \widetilde{\cD}
\cong \bD$ and $\widetilde{M}/\tau \widetilde{M} \cong \bM$. The graded dual
$\widetilde{M}^\vee := \underline{\Hom}_{\widetilde{\cD}}(\widetilde{M}, \widetilde{\cD})$ is the
Rees module of $M^* \cong \Hom_{\cD}(M, \cD)$ with respect to the dual filtration
$(M^*)_k := \{ q \in M^* : M_j\, q \subseteq \cD_{j+k} \text{ for all } j\}$. Put
$E^i := \underline{\Ext}^i_{\widetilde{\cD}}(\widetilde{M}, \widetilde{\cD})$. These are finitely
generated graded $\widetilde{\cD}$-modules, and $E^i[\tau^{-1}] = 0$ for $i \ge 1$, since $M$ is
projective. Applying $\underline{\Hom}_{\widetilde{\cD}}(\widetilde{M}, \,\text{--}\,)$ to the exact
sequence $0 \to \widetilde{\cD} \xrightarrow{\tau} \widetilde{\cD} \to \bD \to 0$ (we suppress
shifts of grading) and using the isomorphisms $\underline{\Ext}^i_{\widetilde{\cD}}(\widetilde{M},
\bD) \cong \Ext^i_{\bD}(\bM, \bD)$, which hold because $\tau$ is a non-zero-divisor on
$\widetilde{M}$, we get the exact sequence
\[
  0 \to \widetilde{M}^\vee \xrightarrow{\tau} \widetilde{M}^\vee \to \bM^\vee \to E^1
  \xrightarrow{\tau} E^1 \to \Ext^1_{\bD}(\bM, \bD) \to E^2 \xrightarrow{\tau} E^2 \to
  \Ext^2_{\bD}(\bM, \bD) .
\]
Since $\bM$ is a submodule of finite codimension in the projective module $\bD\mathcal{I}_M$ over
the regular ring $\bD$ of dimension $2$, its projective dimension is at most $1$; hence
$\Ext^2_{\bD}(\bM, \bD) = 0$. It follows that $E^2 = \tau E^2$, and therefore $E^2 = 0$ by the
graded Nakayama lemma. Next, $E^1$ is annihilated by a power of $\tau$, and $E^1/\tau E^1 \cong
\Ext^1_{\bD}(\bM, \bD)$ is finite-dimensional by Lemma~\ref{L2.7}; since each
$\tau^i E^1/\tau^{i+1} E^1$ is a quotient of $E^1/\tau E^1$, we get $\dim_{\C} E^1 < \infty$,
and
\[
  \dim \Ker(\tau; E^1) = \dim \Coker(\tau; E^1) = \dim \Ext^1_{\bD}(\bM, \bD) = n(M) .
\]
On the other hand, the above exact sequence shows that $\gr(M^*) = \widetilde{M}^\vee/\tau
\widetilde{M}^\vee$ is a submodule of the reflexive module $\bM^\vee$ with quotient
$\Ker(\tau; E^1)$. Thus the dual filtration on $M^*$ is good, $\gr(M^*)$ is torsion-free, its
reflexive hull is $\bM^\vee$, and $\dim \mathsf{V}_{M^*} = \dim \Ker(\tau; E^1) = n(M)$. By
(the right-handed version of) Theorem~\ref{T2.3}$(c)$, we conclude that $n(M^*) = n(M)$.
\end{proof}

\begin{proposition}
\label{P2.3}
Let $M$ be a left ideal of $\cD$. Then $n(M) = 0$ if and only if $M \cong \cD \otimes_{\cO}
\mathcal{I}$ for some line bundle $\mathcal{I}$ on $X$; in this case, $[\mathcal{I}] =
\gamma[M]$. In particular, $M$ is cyclic if and only if $n(M) = 0$ and $\gamma[M]$ is trivial.
\end{proposition}

\begin{proof}
If $M = \cD \otimes_{\cO} \mathcal{I}$, then, for the natural filtration $M_k = \cD_k
\otimes_{\cO} \mathcal{I}$, we have $\gr M = \bD \otimes_{\cO} \mathcal{I}$, which is reflexive;
hence $n(M) = 0$ by Theorem~\ref{T2.3}$(c)$. Conversely, suppose that $n(M) = 0$. Then, by
Theorem~\ref{T2.3}, $\bM \cong \bD \otimes_{\cO} \mathcal{I}_M$ as graded modules (up to a shift
of grading). Let $p$ be the least integer such that $M_p \ne 0$. Then $M_p \cong \mathcal{I}_M$
as $\cO$-modules, and the multiplication map $\cD \otimes_{\cO} M_p \to M$ is a homomorphism of
filtered $\cD$-modules, which induces an isomorphism on the associated graded modules. Hence it
is an isomorphism. The last statement follows from the fact that the class of
$\cD \otimes_{\cO} \mathcal{I}$ in $K_0(\cD) \cong K_0(X)$ corresponds to $[\mathcal{I}]$.
\end{proof}

Finally, we observe that the Calogero--Moser index is a sum of local contributions. Let $M$ be a
right ideal of $\cD$. Since $\mathsf{V}_M$ is a finite-dimensional $\bD$-module, it has finite
length as a module over $\cO \subset \bD$; hence $\mathsf{V}_M = \bigoplus_{\lambda \in X}
(\mathsf{V}_M)_\lambda$, where $(\mathsf{V}_M)_\lambda := \mathsf{V}_M \otimes_{\cO} \cO_\lambda$,
and
\begin{equation}
\label{E211}
  n(M) = \sum_{\lambda \in X} n_\lambda(M) , \qquad n_\lambda(M) := \dim_{\C}
  (\mathsf{V}_M)_\lambda .
\end{equation}
We call $n_\lambda(M)$ the \emph{local Calogero--Moser index} of $M$ at $\lambda$. It can be
computed after localization and completion. Precisely, let $\widehat{M}_\lambda := M\, \hD_\lambda
\cong M \otimes_{\cO} \hO_\lambda$ be the right ideal of $\hD_\lambda$ generated by $M$, with
filtration induced by the order filtration on $\hD_\lambda$. Then $\gr \widehat{M}_\lambda \cong
\bM \otimes_{\cO} \hO_\lambda$ is an ideal of $\gr \hD_\lambda \cong \hO_\lambda[\xi]$, and, since
the formation of the reflexive hull commutes with flat base change, we have
\begin{equation}
\label{E212}
  n_\lambda(M) = \dim_{\C}\bigl[\, (\gr \widehat{M}_\lambda)^{\vee\vee} \big/ \gr
  \widehat{M}_\lambda \,\bigr] .
\end{equation}
Thus $n_\lambda(M)$ depends only on the right ideal $\widehat{M}_\lambda \subseteq \hD_\lambda
\cong \C[[t]][\del]$.

\subsection{The affine line}
\label{S2.6}
In the case $X = \A^1$, the above classification has additional features, which we will use in
Sections~\ref{S5} and~\ref{S6}. We write $x$ for the coordinate on $\A^1$, so that $\cO = \C[x]$,
$K = \C(x)$ and $\cD = A_1$.

\subsubsection*{The adelic Grassmannian}
If $W = \bigcap_\lambda W_\lambda$ is a primary decomposable subspace of $\C[x]$ and $k_\lambda :=
\codim W_\lambda$, the subspace $\prod_\lambda (x - \lambda)^{-k_\lambda}\, W \subset \C(x)$ depends
only on the equivalence class of $W$. The set of all subspaces of $\C(x)$ arising in this way is
Wilson's \emph{adelic Grassmannian} $\Gra$ (see \cite{Wilson1993, Wilson1998}). By
Theorem~\ref{T2.1}, $\,W \mapsto [\cD(\C[x], W)]$ is a bijection $\Gra \to \cR(A_1)$, and
$W \mapsto [\cD(W, \C[x])]$ is a bijection $\Gra \to \cL(A_1)$. To each $W \in \Gra$ Wilson
associates its (stationary) \emph{Baker function}, which has the form
\begin{equation}
\label{E213}
  \psi_W(x, z) = e^{xz} \Bigl( 1 + \sum_{i \ge 1} a_i(z)\, x^{-i} \Bigr) ,
\end{equation}
where the $a_i(z)$ are rational functions of $z$, and which determines $W$ uniquely (see
\cite[Section~4]{Wilson1993}).\footnote{We interchange the roles of the variables $x$ and $z$ as
compared to \cite{Wilson1993, Wilson1998}, since in this paper $x$ is the coordinate on the curve.
Thus, our $\psi_W(x,z)$ is the function denoted by $\psi_W(z,x)$ in \emph{loc.~cit.}; by
\cite[Theorem~2]{Wilson1993}, it is the Baker function, in Wilson's notation, of the image $b(W)$
of $W$ under the bispectral involution.} For a rational function $a(z)$, we write $\res_{\infty}\, a$ for
the coefficient of $z^{-1}$ in the Laurent expansion of $a(z)$ at $z = \infty$; equivalently,
$\res_\infty a = \sum_{p \in \C} \Res_{z = p}\, a(z)$.

\begin{proposition}[\cite{Wilson1993, Wilson1998}]
\label{P2.4}
Let $C$ be a space of conditions on $\A^1$, with basis $\{c_1, \ldots, c_g\}$, where $c_i \in
C_{\lambda_i}$ is represented by $F_i(\del) \in \C[\del]$. Let $W \in \Gra$ be the point
corresponding to $W_C$, and let $\psi_W$ be its Baker function \eqref{E213}.
\begin{enumerate}
\item[$(a)$] The Calogero--Moser index of the right ideal $\cD(\C[x], W) \cong \Ann_{A_1}(C)$ is
given by $\, n = - \res_{\infty}\, a_1(z)$.
\item[$(b)$] Let $f_i(z) := \lan c_i, e^{xz} \ran_{\lambda_i} = F_i(z)\, e^{\lambda_i z}$. Then,
up to an additive constant, $\,a_1(z) = - \del_z \log \Wr(f_1, \ldots, f_g)$, where
$\Wr(f_1, \ldots, f_g) = \det\bigl(\del_z^{\,i-1} f_j\bigr)_{i,j = 1}^{g}$ is the Wronskian.
\item[$(c)$] $\psi_W(x,z) = e^{xz}\, \det \mathsf{M}(z, x)/\det \mathsf{M}(z, \infty)$, where
$\mathsf{M}(z,x)$ is the $g \times g$ matrix with entries
\[
  \mathsf{m}_{ij} = (\del_z + \lambda_i)^{j-1} \bigl[\, F_i(z) - (x - \lambda_i)^{-1} F_i'(z)
  \,\bigr] , \qquad i,j = 1, \ldots, g .
\]
\end{enumerate}
\end{proposition}

Part $(a)$ follows from Wilson's formula for the Baker function of the point $W \in \Gra$
corresponding to $(\mathsf{X}, \mathsf{Y}, v, w) \in \cC_n$ (see \cite[(3.7)]{Wilson1998}): in
our notation, $e^{-xz}\, \psi_W(x,z) = \det\{ \mathrm{Id} - (z + \mathsf{X})^{-1} (x +
\mathsf{Y})^{-1}\}$, whence $a_1(z) = - \mathrm{tr}\,(z + \mathsf{X})^{-1}$ and
$\res_\infty a_1 = -n$. The fact that the stratification of $\Gra \cong \cR(A_1)$ by the spaces
$\cC_n$ agrees with the one given by the index \eqref{E210} is proved in
\cite[Theorems~1.2 and 4.5]{BerestWilsonNC} (see also \cite[Proposition~5.2]{BerestWilsonNC}).
Parts $(b)$ and $(c)$ are \cite[Lemma~5.8]{Wilson1993} and \cite[Proposition~6.5]{Wilson1993},
respectively. (The additive constant in $(b)$ arises because Wilson's lemma refers to the subspace
$x^{-g}\, W_C$, which differs from $W$ by a factor that is a rational function of $x$ equal to $1$
at infinity; it does not affect $\res_\infty a_1$.)

\subsubsection*{Distinguished representatives}
Unlike the general case, every ideal class of $A_1$ has two distinguished representatives (see
\cite[Section~5.1]{BerestWilsonNC}, where the case of right ideals is treated; the case of left
ideals is similar). Namely, for any left ideal $M \subseteq A_1$, there is a
unique fractional ideal $M_x \cong M$ such that $M_x \subset \C(x)[\del]$, $M_x \cap \C[x] \ne 0$,
all leading coefficients of the elements of $M_x$ are polynomials, and $M_x$ contains an element
with leading coefficient $1$. Similarly, there is a unique fractional ideal $M_\del \cong M$ with
the same properties with respect to $\del$ (i.e.\ obtained by interchanging the roles of $x$ and
$\del$). We write $\kappa = \kappa_M$ for the unique element of $Q$ such that
\begin{equation}
\label{E214}
  M_x = M_\del\, \kappa .
\end{equation}
The element $\kappa$ determines the isomorphism class of $M$ and is related to the Baker function
as follows (see \cite[Proposition~9.2]{BerestWilson00} and
\cite[Remark~5.4]{BerestWilsonNC}): if $M \cong \cD(W, \C[x])$ with
$W \in \Gra$, then
\begin{equation}
\label{E215}
  \kappa^{-1}(x, \del)\,[e^{xz}] = \psi_W(x,z) .
\end{equation}
Here $\kappa^{-1}$ is regarded as a formal pseudo-differential operator
$1 + \sum_{i \ge 1} x^{-i} a_i(\del)$ acting on $e^{xz}$ in the obvious way. For example, if
$f \in \C[x]$ is a polynomial of degree $n$, then $\Ann(f) = A_1 \del^{n+1} + A_1 (f \del - f')$,
and one finds $M_\del = M \del^{-1}$, $M_x = M (\del - f'/f)^{-1}$, so that
\begin{equation}
\label{E216}
  \kappa^{-1} = 1 - \frac{f'(x)}{f(x)}\, \del^{-1} .
\end{equation}
Finally, we recall that the group $G = \Aut_{\C}(A_1)$ acts naturally on $\cL(A_1)$ and
$\cR(A_1)$, and the bijection $\omega$ of Theorem~\ref{T2.2} is equivariant with respect to a
natural action of $G$ on $\cC = \bigsqcup_n \cC_n$, which preserves each $\cC_n$ and is transitive
on it (see \cite[Theorems~1.2 and 1.3]{BerestWilson00}). In particular, the Calogero--Moser index
is invariant under the action of $G$ on ideals.

\section{The holonomicity index}
\label{S3}

In this section, we define the holonomicity index of a differential operator and relate it to the
spaces of conditions introduced in Section~\ref{S2.4}. Throughout, $L$ denotes a nonzero element
of $\cD = \cD(X)$.

\subsection{Singular points}
\label{S3.1}
For each $\lambda \in X$, right multiplication by $L$ defines a linear endomorphism of the vector
space $\Dl = \cD/\ml \cD$, and we put
\begin{equation}
\label{E31}
  V_\lambda(L) := \Ker\bigl[\, \Dl \xrightarrow{\ \cdot L\ } \Dl \,\bigr] .
\end{equation}

\begin{definition}
\label{D3}
A point $\lambda \in X$ is called a \emph{singular point} of $L$ if $V_\lambda(L) \ne 0$. We write
$\Lambda(L)$ for the set of all singular points of $L$, and
$V_\Lambda(L) := \bigoplus_{\lambda \in \Lambda(L)} V_\lambda(L) \subseteq \Delta(X)$.
\end{definition}

By Lemma~\ref{L2.2}$(d)$, $\lambda$ is a singular point of $L$ if and only if the equation
$L(u) = f$ fails to be solvable in $\hO_\lambda$ for some $f \in \hO_\lambda$. In particular,
$\lambda \notin \Lambda(L)$ whenever the leading coefficient of $L$ (with respect to a local
parameter at $\lambda$) does not vanish at $\lambda$, and therefore $\Lambda(L)$ is a finite set.
We will see below that each $V_\lambda(L)$ is finite-dimensional (Corollary~\ref{C3.1}), so that
$V_\Lambda(L)$ is a space of conditions in the sense of Definition~\ref{D2}, and that
$\Lambda(L)$ is precisely the support of the $\cO$-module $\cl(L)/\cD L$
(Theorem~\ref{T4.1}). Note that a zero of the leading coefficient of $L$ need not be a singular
point: for example, $L = x\del + 1 = \del x \in A_1$ has no singular points.

\subsection{Normalization}
\label{S3.2}
Fix $\lambda \in X$ and a local parameter $t$ at $\lambda$. Using \eqref{E23}, we can write $L$ as
an element of $\cD_\lambda$:
\begin{equation}
\label{E32}
  L = a_N \del^N + a_{N-1} \del^{N-1} + \ldots + a_0 , \qquad a_i \in \cO_\lambda , \quad
  \del = \del_t .
\end{equation}
For each $i$ with $a_i \ne 0$, let $k_i := \ord_\lambda(a_i)$, so that $a_i = b_i\, t^{k_i}$ with
$b_i \in \cO_\lambda \setminus \nl$. In view of \eqref{E24}, we may rewrite \eqref{E32} in the
form
\begin{equation}
\label{E33}
  L = \sum_{i \,:\, a_i \ne 0} b_i\, t^{k_i - i}\, \ff{\theta}{i}\, , 
\end{equation}
where $\,\theta = t \del\,$ and $\, \ff{\theta}{i} := t^i \partial^i = \theta(\theta-1) \ldots (\theta -i + 1)  $.

\begin{definition}
\label{D4}
We say that $L$ is \emph{$\lambda$-normalized} (or \emph{balanced at $\lambda$}) if
\begin{enumerate}
\item[$(a)$] $k_i - i \ge 0$ for all $i$ with $a_i \ne 0$, and
\item[$(b)$] there is $i \in \{0, \ldots, N\}$ such that $a_i \ne 0$ and $k_i = i$.
\end{enumerate}
In this case, the \emph{indicial polynomial} of $L$ at $\lambda$ is defined by
$\,q_{L,\lambda}(s) := \sum_{i \,:\, k_i = i} b_i(\lambda)\, \ff{s}{i} \in \C[s]$.
\end{definition}

Note that $q_{L,\lambda}$ is a nonzero polynomial of degree $\le N$, since the polynomials
$\ff{s}{i}$ have distinct degrees. The next two lemmas are almost immediate from the definition.

\begin{lemma}
\label{L3.1}
For any nonzero $L \in \cD_\lambda$, there is a unique $w = w_\lambda(L) \in \Z$ such that
$t^{w} L$ is $\lambda$-normalized, namely $w = \max\{ i - k_i : a_i \ne 0\}$. The integer $w$ and
the property of being $\lambda$-normalized do not depend on the choice of local parameter.
\end{lemma}

\begin{proof}
The first claim is clear from \eqref{E32}: multiplying $L$ by $t^w$ replaces each $k_i$ by
$k_i + w$. The independence of the choice of $t$ follows from Lemma~\ref{L3.2} below, which shows
that $-w$ is the largest integer $m$ such that $L(\nl^j) \subseteq \nl^{j+m}$ for all
$j \gg 0$.
\end{proof}

We call $\bar L := t^{w_\lambda(L)} L$ the \emph{$\lambda$-normalization} of $L$; it is determined
by $L$ up to a factor which is a unit in $\cO_\lambda$.

\begin{lemma}
\label{L3.2}
If $L$ is $\lambda$-normalized, then $L(\nl^j) \subseteq \nl^j$ for all $j \ge 0$, and the induced
action of $L$ on the one-dimensional vector space $\nl^j/\nl^{j+1}$ is multiplication by
$q_{L,\lambda}(j)$. In particular, $L(\nl^j) \not\subseteq \nl^{j+1}$ for $j \gg 0$.
\end{lemma}

\begin{proof}
We have $\nl^j = \cO_\lambda\, t^j$, and $\theta(a\, t^j) = (t\,\del(a) + j\,a)\, t^j$ for all
$a \in \cO_\lambda$. Hence $\theta$ preserves $\nl^j$ and acts on $\nl^j/\nl^{j+1}$ as
multiplication by $j$. The claim follows now from \eqref{E33}, since the terms with $k_i - i > 0$
map $\nl^j$ into $\nl^{j+1}$.
\end{proof}

\subsection{Solution spaces and the cotype}
\label{S3.3}
By Lemma~\ref{L3.2}, the $\lambda$-normalization $\bar L$ of $L$ induces linear endomorphisms of
the finite-dimensional vector spaces $\cO_\lambda/\nl^{j+1} \cong \hO_\lambda/\hn_\lambda^{j+1}$
for all $j \ge 0$. Generalizing the construction of \cite[Section~4]{Tsai}, we define the
spaces $S^j_\lambda(L)$ by formula \eqref{E13} of the Introduction:
$\,S^j_\lambda(L) := \Ker[\,\bar L;\ \cO_\lambda/\nl^{j+1}\,]$ for $j \ge 0$, and
$S^{-1}_\lambda(L) := 0$. Clearly, the dimensions of these spaces do not depend on the choice of
$t$.

\begin{lemma}
\label{L3.3}
Let $q = q_{\bar L, \lambda}$ be the indicial polynomial of the $\lambda$-normalization of $L$.
\begin{enumerate}
\item[$(a)$] For all $j \ge 0$, either $\dim S^j_\lambda(L) = \dim S^{j-1}_\lambda(L)$ or
$\dim S^j_\lambda(L) = \dim S^{j-1}_\lambda(L) + 1$. The second possibility may occur only if
$q(j) = 0$.
\item[$(b)$] The natural maps $S^j_\lambda(L) \to S^{j-1}_\lambda(L)$ are isomorphisms for
$j \gg 0$.
\item[$(c)$] For $j \gg 0$, $\ \dim S^j_\lambda(L) = \dim \Ker(L;\, \hO_\lambda) = \dim \Coker(\bar L;
\, \hO_\lambda)$.
\end{enumerate}
\end{lemma}

\begin{proof}
The operator $\bar L$ induces an endomorphism of the short exact sequence
\[
  0 \to \nl^j/\nl^{j+1} \to \cO_\lambda/\nl^{j+1} \to \cO_\lambda/\nl^{j} \to 0 ,
\]
which, by Lemma~\ref{L3.2}, acts on the first term as multiplication by $q(j)$. By the Snake
Lemma, we get an exact sequence
\[
  0 \to \Ker(q(j)) \to S^j_\lambda(L) \to S^{j-1}_\lambda(L) \xrightarrow{\ \delta\ } \Coker(q(j)) .
\]
If $q(j) \ne 0$, then $\Ker(q(j)) = \Coker(q(j)) = 0$, and $S^j_\lambda(L) \cong
S^{j-1}_\lambda(L)$. If $q(j) = 0$, then $\Ker(q(j)) \cong \Coker(q(j)) \cong \C$, and
$\dim S^j_\lambda(L) = \dim S^{j-1}_\lambda(L) + 1 - \mathrm{rk}(\delta)$. This proves $(a)$. Since
$q \ne 0$ has finitely many roots, $(b)$ follows. To prove $(c)$, choose $m$ so large that
$q(j) \ne 0$ for all $j \ge m$. Then $\bar L$ preserves the filtration
$\{\hn^j_\lambda\}_{j \ge m}$ on $\hn^m_\lambda$ and induces isomorphisms on its graded pieces;
since $\hn^m_\lambda$ is complete, $\bar L: \hn_\lambda^m \to \hn_\lambda^m$ is bijective.
Applying the Snake Lemma to the action of $\bar L$ on $0 \to \hn^m_\lambda \to \hO_\lambda \to
\hO_\lambda/\hn^m_\lambda \to 0$, we see that the kernel and the cokernel of $\bar L$ on
$\hO_\lambda$ are isomorphic to those of $\bar L$ on $\hO_\lambda/\hn^m_\lambda$; the latter have
the same dimension, equal to $\dim S^{m-1}_\lambda(L)$. Finally, $\Ker(L; \hO_\lambda) =
\Ker(\bar L; \hO_\lambda)$.
\end{proof}

Using Lemma~\ref{L3.3}, we introduce the following numerical invariants of $L$ at $\lambda$:
\begin{align}
  d_\lambda & := \dim S^j_\lambda(L) \quad (j \gg 0) , \notag \\
  J_\lambda(L) & := \{ j \in \Z_{\ge 0} \, : \, \dim S^j_\lambda(L) \ne \dim S^{j-1}_\lambda(L) \} ,
  \label{E34} \\
  g_\lambda & := \max J_\lambda(L) , \qquad s_\lambda := |J_\lambda(L)| - 1 . \notag
\end{align}
By Lemma~\ref{L3.3}$(c)$, $d_\lambda$ is the number of linearly independent formal power series
solutions of the equation $L(u) = 0$ at $\lambda$. We call the finite set $J_\lambda(L)$ the
\emph{cotype} of $L$ at $\lambda$, and write its elements in increasing order:
$J_\lambda(L) = \{ j_0 < j_1 < \ldots < j_{s_\lambda} \}$. By Lemma~\ref{L3.3}$(a)$, the cotype is
contained in the set of non-negative integer roots of the indicial polynomial of $\bar L$.

\begin{lemma}
\label{L3.4}
We have $\ (a)\ s_\lambda = d_\lambda - 1$, $\ (b)\ g_\lambda \ge d_\lambda - 1$ (provided
$d_\lambda > 0$).
\end{lemma}

\begin{proof}
$(a)$ By Lemma~\ref{L3.3}, $d_\lambda = \sum_{j \ge 0}\, [\dim S^j_\lambda(L) - \dim
S^{j-1}_\lambda(L)] = |J_\lambda(L)|$. $(b)$ Since the $j_i$ are distinct non-negative integers,
$j_i \ge i$ for all $i$; hence $g_\lambda = j_{s_\lambda} \ge s_\lambda = d_\lambda - 1$.
\end{proof}

We can now state our main definition.

\begin{definition}
\label{D5}
The \emph{local holonomicity index} of $L$ at $\lambda \in X$ is defined by
\begin{equation}
\label{E35}
  n_\lambda(L) := \nu(J_\lambda(L)) = \sum_{i=0}^{s_\lambda} (j_i - i) = \sum_{i=0}^{d_\lambda - 1}
  j_i \, - \binom{d_\lambda}{2} ,
\end{equation}
where $\{j_0 < \ldots < j_{s_\lambda}\}$ is the cotype of $L$ at $\lambda$. The (global)
\emph{holonomicity index} of $L$ is
\begin{equation}
\label{E36}
  n(L) := \sum_{\lambda \in \Lambda(L)} n_\lambda(L) .
\end{equation}
\end{definition}

Thus $n_\lambda(L) \ge 0$, and $n_\lambda(L) = 0$ if and only if the cotype of $L$ at $\lambda$ is
an initial segment $\{0, 1, \ldots, d_\lambda - 1\}$ of $\Z_{\ge 0}$. This is the case, for
example, if the leading coefficient $a_N$ of $L$ does not vanish at $\lambda$: then
$\bar L = t^N L$, and it is easy to see that $J_\lambda(L) = \{0, 1, \ldots, N-1\}$. More
generally, we will show that $n_\lambda(L) = 0$ for all $\lambda \notin \Lambda(L)$ (see
Corollary~\ref{C3.1}), so that the sum in \eqref{E36} may be taken over all points of $X$.

\subsection{Comparison with spaces of conditions}
\label{S3.4}
The next theorem relates the solution spaces $S^j_\lambda(L)$ to the space $V_\lambda(L)$ defined
in \eqref{E31}. For $X = \A^1$, this was proved in \cite[Section~4]{Tsai} by comparing the
matrices that compute both spaces.

\begin{theorem}
\label{T3.1}
Assume that $L$ is $\lambda$-normalized. Then right multiplication by $L$ preserves the filtration
$\{\Dl^{\le j}\}$ on $\Dl$, and for all $j \ge 0$,
\[
  \dim S^j_\lambda(L) = \dim V_\lambda(L)^{\le j} ,
\]
where $V_\lambda(L)^{\le j} = V_\lambda(L) \cap \Dl^{\le j}$. Consequently, $V_\lambda(L)$ is a
finite-dimensional vector space of dimension $d_\lambda$, and its degree set coincides with the
cotype of $L$ at $\lambda$: $\ \Sigma(V_\lambda(L)) = J_\lambda(L)$.
\end{theorem}

\begin{proof}
By Lemma~\ref{L2.2}$(b)$, an element $F \in \Dl$ belongs to $\Dl^{\le j}$ if and only if
$\lan F, \hn^{j+1}_\lambda \ran_\lambda = 0$, and the pairing identifies $\Dl^{\le j}$ with the
dual of $\hO_\lambda/\hn_\lambda^{j+1} \cong \cO_\lambda/\nl^{j+1}$. Since $L$ preserves
$\hn^{j+1}_\lambda$ (Lemma~\ref{L3.2}), it follows from Lemma~\ref{L2.2}$(c)$ that
$\lan F L, \hn^{j+1}_\lambda\ran_\lambda = \lan F, L(\hn^{j+1}_\lambda) \ran_\lambda = 0$ for all
$F \in \Dl^{\le j}$, i.e.\ $\Dl^{\le j} \cdot L \subseteq \Dl^{\le j}$, and that the restriction
of $R_L$ to $\Dl^{\le j}$ is the transpose of the endomorphism $L$ of $\cO_\lambda/\nl^{j+1}$.
Hence
\begin{align*}
  \dim V_\lambda(L)^{\le j} & = \dim \Ker\bigl(R_L |_{\Dl^{\le j}}\bigr) = \dim \Coker(L;\,
  \cO_\lambda/\nl^{j+1}) \\
  & = \dim \Ker(L;\, \cO_\lambda/\nl^{j+1}) = \dim S^j_\lambda(L) ,
\end{align*}
where the third equality holds because $\cO_\lambda/\nl^{j+1}$ is finite-dimensional. The last
statement follows from Lemma~\ref{L3.3} and Lemma~\ref{L2.3}$(a)$.
\end{proof}

\begin{remark}
In terms of Proposition~\ref{P2.1}$(b)$ and Lemma~\ref{L2.1}, the identification of $\Dl^{\le j}$
with the dual of $\cO_\lambda/\nl^{j+1}$ used in the above proof is the natural isomorphism
$\cD^j_\lambda/\nl \cD^j_\lambda \cong \cD^j_{\cO_\lambda}(\cO_\lambda, \cO_\lambda/\nl) \cong
\Hom_{\C}(\cO_\lambda/\nl^{j+1}, \C)$, where $\cD^j_\lambda$ denotes the space of operators of
order $\le j$ in $\cD_\lambda$. Thus Theorem~\ref{T3.1} is a manifestation of the formal
continuity of differential operators.
\end{remark}

Theorem~\ref{T3.1} applies to $\lambda$-normalized operators. For an arbitrary $L \in \cD$, the
space $V_\lambda(L)$ differs from $V_\lambda(\bar L)$, but in a controlled way.

\begin{corollary}
\label{C3.1}
Let $L$ be a nonzero element of $\cD$, $\lambda \in X$, and let $\bar L = t^w L$ be the
$\lambda$-normalization of $L$.
\begin{enumerate}
\item[$(a)$] $V_\lambda(L)$ is finite-dimensional, and $V_\lambda(L)^\perp = L(\hO_\lambda)$.
\item[$(b)$] If $w \ge 0$, then $\,J_\lambda(L) = \Sigma(V_\lambda(\bar L)) = \{0, \ldots, w-1\}
\sqcup (w + \Sigma(V_\lambda(L)))$; if $w < 0$, then $\,\Sigma(V_\lambda(L)) = \{0, \ldots, -w-1\}
\sqcup (-w + J_\lambda(L))$.
\item[$(c)$] In either case, $\ n_\lambda(L) = \nu(V_\lambda(L))$. Explicitly, if
$\{F_1, \ldots, F_g\}$ is a basis of $V_\lambda(L)$ consisting of elements of distinct degrees,
then
\begin{equation}
\label{E37}
  n_\lambda(L) = \sum_{i=1}^g \deg(F_i) - \binom{g}{2} .
\end{equation}
\item[$(d)$] $n_\lambda(L) = 0$ for all $\lambda \notin \Lambda(L)$.
\end{enumerate}
\end{corollary}

\begin{proof}
By Lemma~\ref{L2.2}$(d)$, $V_\lambda(P)^\perp = P(\hO_\lambda)$ for any $P \in \hD_\lambda$; hence
$V_\lambda(\bar L)^\perp = t^w\, V_\lambda(L)^\perp$. Since $V_\lambda(\bar L)$ is
finite-dimensional by Theorem~\ref{T3.1}, the subspace $V_\lambda(L)^\perp$ is closed and has
finite codimension in $\hO_\lambda$, which implies $(a)$. Part $(b)$ follows then from
Lemma~\ref{L2.3}$(c)$ and Theorem~\ref{T3.1}, and $(c)$ follows from $(b)$ and the last formula of
Lemma~\ref{L2.3}$(c)$. Finally, $(d)$ is a special case of $(c)$.
\end{proof}

\begin{remark}
At the cost of using cokernels instead of kernels, one can avoid normalization altogether. If
$L \in \cD$ is an operator of order $N$, then $L(\nl^{j+N+1}) \subseteq \nl^{j+1}$ for all
$j \ge 0$ by Proposition~\ref{P2.1}$(a)$, and the argument of Theorem~\ref{T3.1} shows that
\[
  \dim V_\lambda(L)^{\le j} = \dim \Coker\bigl[\, L : \ \cO_\lambda/\nl^{j+N+1} \to
  \cO_\lambda/\nl^{j+1} \,\bigr] , \qquad j \ge 0 .
\]
Thus the degree set of $V_\lambda(L)$, and hence $n_\lambda(L)$, can be read off from the jumps of
the dimensions on the right.
\end{remark}

\subsection{The main theorem}
\label{S3.5}
We can now give a precise formulation of our main result. Recall that $V_\Lambda(L) =
\bigoplus_{\lambda} V_\lambda(L)$ is a space of conditions on $X$ (see Definition~\ref{D3} and
Corollary~\ref{C3.1}$(a)$); we denote by
\begin{equation}
\label{E38}
  M(L) := \Ann_{\cD}(V_\Lambda(L)) = \{ P \in \cD \, : \, F P = 0 \ \text{in}\ \Delta(X) \
  \text{for all}\ F \in V_\Lambda(L) \}
\end{equation}
the corresponding fat right ideal of $\cD$ (see \eqref{E27}).

\begin{theorem}
\label{T3.2}
Let $L$ be a nonzero differential operator on $X$.
\begin{enumerate}
\item[$(a)$] The right ideal $M(L)$ is dual to the closure of $L$; precisely,
$\,M(L) = (\cl(L)\, L^{-1})^*$.
\item[$(b)$] For every $\lambda \in X$, the local holonomicity index of $L$ equals the local
Calogero--Moser index of $M(L)$: $\ n_\lambda(L) = n_\lambda(M(L))$.
\item[$(c)$] $\,n(L) = n(M(L)) = n(\cl(L))$.
\end{enumerate}
\end{theorem}

Part $(a)$ of Theorem~\ref{T3.2} will be proved in Section~\ref{S4} (see Corollary~\ref{C4.1}),
and parts $(b)$ and $(c)$ in Section~\ref{S5.3}. We conclude this section with two immediate
consequences.

\begin{corollary}
\label{C3.2}
Let $f_1, \ldots, f_k$ be holomorphic functions on $X$ such that the left ideal
$I = \Ann_{\cD}(f_1, \ldots, f_k)$ is nonzero, and let $L$ be an operator of minimal order in $I$.
Then $I = \cl(L)$ and $n(I) = n(L)$. In particular, a holonomic function $f$ on $X$ is strongly
holonomic if and only if $n_\lambda(L) > 0$ for some $\lambda \in \Lambda(L)$, where $L$ is an
operator of minimal order in $\Ann_{\cD}(f)$.
\end{corollary}

\begin{proof}
By Lemma~\ref{L2.4}, the ideal $I$ is closed, and $I = \cl(L)$. The result follows from
Theorem~\ref{T3.2}$(c)$.
\end{proof}

\begin{corollary}
\label{C3.3}
For a nonzero $L \in \cD$, the following are equivalent:
\begin{enumerate}
\item[(i)] $n(L) = 0$;
\item[(ii)] the cotype of $L$ at every point $\lambda \in X$ is $\{0, 1, \ldots, d_\lambda - 1\}$;
\item[(iii)] $\cl(L) \cong \cD \otimes_{\cO} \mathcal{I}$ for a line bundle $\mathcal{I}$ on $X$.
\end{enumerate}
\end{corollary}

\begin{proof}
The equivalence $(i) \Leftrightarrow (ii)$ is clear from \eqref{E35}, and
$(i) \Leftrightarrow (iii)$ follows from Theorem~\ref{T3.2}$(c)$ and Proposition~\ref{P2.3}.
\end{proof}

\section{Dual classification of ideals}
\label{S4}

In this section we describe left ideals of $\cD$ in terms of spaces of conditions. As a
consequence, we obtain a presentation of the closure of a differential operator and prove part
$(a)$ of Theorem~\ref{T3.2}.

\subsection{A presentation of ideals}
\label{S4.1}
Let $I$ be a (nonzero) left ideal of $\cD$, and let $L \in I$ be an operator of minimal order.
Recall that $\cD L \subseteq I \subseteq \cl(L)$ (see Lemma~\ref{L2.4}$(b)$). For each
$\lambda \in X$, we put
\begin{equation}
\label{E41}
  \widetilde{V}_\lambda(I;L) := \{ F \in \cD \, : \, F L \in \ml I \} , \qquad
  V_\lambda(I;L) := \widetilde{V}_\lambda(I;L)/\ml \cD \subseteq \Dl .
\end{equation}
Thus $V_\lambda(I;L)$ is the kernel of the map $\Dl \to I/\ml I$ induced by right multiplication
by $L$, as in \eqref{E18}. When $I = \cl(L)$, we simply write $\widetilde{V}_\lambda(L) := \{F \in
\cD : FL \in \ml \cD\}$ for the preimage of $V_\lambda(L)$ in $\cD$. As usual, $\ml^{-1} := \{ h
\in K : h\, \ml \subseteq \cO\}$ denotes the fractional ideal inverse to $\ml$; recall that
$\ml^{-1} \ml = \cO$ and that $\ml^{-1}/\cO$ is a one-dimensional vector space.

\begin{lemma}
\label{L4.1}
Let $N := I/\cD L$. Then, for every $\lambda \in X$,
\begin{enumerate}
\item[$(a)$] $V_\lambda(I;L) \subseteq V_\lambda(L)$, with equality if $I = \cl(L)$;
\item[$(b)$] $N^{\ml} = \bigl(\ml^{-1}\, \widetilde{V}_\lambda(I;L)\, L + \cD L\bigr)/\cD L\,$,
and there is an isomorphism $N^{\ml} \cong (\ml^{-1}/\cO) \otimes_{\C} V_\lambda(I;L)$. In
particular, $\dim N^{\ml} = \dim V_\lambda(I;L)$.
\end{enumerate}
\end{lemma}

\begin{proof}
$(a)$ The inclusion is obvious, since $\ml I \subseteq \ml \cD$. If $I = \cl(L)$ and
$F L \in \ml \cD$, then $\ml^{-1} F L \subseteq \cD \cap \cD(K) L = \cl(L)$, and hence
$F L \in \ml \ml^{-1} F L \subseteq \ml\, \cl(L)$. $(b)$ By definition, $N^{\ml}$ consists of the
classes of elements $P \in I$ such that $\ml P \subseteq \cD L$. Writing $P = q L$ with
$q \in \cD(K)$, we see that this condition is equivalent to $\ml q \subseteq \cD$. Then
$(\ml q) L \subseteq \ml I$, i.e.\ $\ml q \subseteq \widetilde{V}_\lambda(I;L)$, and therefore
$q \in \ml^{-1} \ml q \subseteq \ml^{-1} \widetilde{V}_\lambda(I;L)$. Conversely, if
$q \in \ml^{-1} \widetilde{V}_\lambda(I;L)$, then $q L \in \ml^{-1} \ml I = I$ and $\ml q L
\subseteq \widetilde{V}_\lambda(I;L)\, L \subseteq \cD L$. This proves the first claim. Since
$q L \in \cD L$ if and only if $q \in \cD = \ml^{-1} (\ml \cD)$, we get
$N^{\ml} \cong \ml^{-1} \widetilde{V}_\lambda(I;L)/\ml^{-1}(\ml \cD) \cong \ml^{-1} \otimes_{\cO}
V_\lambda(I;L)$, which is isomorphic to $(\ml^{-1}/\cO) \otimes_{\C} V_\lambda(I;L)$, because
$\ml$ annihilates $V_\lambda(I;L)$.
\end{proof}

For $X = \A^1$, the isomorphism of Lemma~\ref{L4.1}$(b)$ is given explicitly by $F \mapsto
(x - \lambda)^{-1} F L$ (cf.\ \cite[Theorem 2.5]{Tsai}).

\begin{theorem}
\label{T4.1}
Let $I$ be a left ideal of $\cD$, and let $L \in I$ be an operator of minimal order. Then
$V_\lambda(I;L) = 0$ for all $\lambda \notin \Lambda(L)$, and
\begin{equation}
\label{E42}
  I = \cD L + \sum_{\lambda \in \Lambda(L)} \cD\, \ml^{-1}\, \widetilde{V}_\lambda(I;L)\, L .
\end{equation}
In particular, for any nonzero $L \in \cD$,
\begin{equation}
\label{E43}
  \cl(L) = \cD L + \sum_{\lambda \in \Lambda(L)} \cD\, \ml^{-1}\, \widetilde{V}_\lambda(L)\, L ,
\end{equation}
and $\Lambda(L)$ is the support of the $\cO$-module $\cl(L)/\cD L$.
\end{theorem}

\begin{proof}
By Lemma~\ref{L2.4}$(c)$, $N = I/\cD L$ is a torsion $\cO$-module; hence $N =
\bigoplus_{\lambda \in X} H^0_\lambda(N)$. By Lemma~\ref{L2.5}, $H^0_\lambda(N) = \cD \cdot
N^{\ml}$ for each $\lambda$, and by Lemma~\ref{L4.1}, $N^{\ml}$ is the image of $\ml^{-1}
\widetilde{V}_\lambda(I;L)\, L$ in $N$; moreover, $N^{\ml} = 0$ unless $V_\lambda(L) \ne 0$.
Formula \eqref{E42} follows. The last statement is the special case $I = \cl(L)$, where
$N^{\ml} \ne 0$ if and only if $V_\lambda(L) \ne 0$.
\end{proof}

Formula \eqref{E43} is Corollary~\ref{C1.1} of the Introduction. For $X = \A^1$, we have
$\ml^{-1} = (x - \lambda)^{-1}\, \C[x]$ and $\widetilde{V}_\lambda(L) = \sum_i \C\, F_{\lambda,i} +
\ml A_1$, where $\{F_{\lambda,i}\}$ is a basis of $V_\lambda(L) \subset \C[\del]$, so that
\eqref{E43} becomes
\begin{equation}
\label{E44}
  \cl(L) = A_1 L + \sum_{\lambda \in \Lambda(L)} \sum_{i} A_1\, (x - \lambda)^{-1} F_{\lambda,i}\, L ,
\end{equation}
which is \cite[Algorithm~3.4]{Tsai}.

\subsection{Ideals defined by spaces of conditions}
\label{S4.2}
Let $C = \bigoplus_\lambda C_\lambda$ be a space of conditions on $X$ (see
Definition~\ref{D2}). For each $\lambda \in \Supp(C)$, let $\widetilde{C}_\lambda \subseteq \cD$
denote the preimage of $C_\lambda$ under the projection $\cD \to \Dl$. We associate to $C$ the
left $\cD$-submodule of $\cD(K)$
\begin{equation}
\label{E45}
  J_C := \cD + \sum_{\lambda \in \Supp(C)} \cD\, \ml^{-1}\, \widetilde{C}_\lambda .
\end{equation}
Clearly, $J_C$ is finitely generated, i.e.\ it is a fractional left ideal of $\cD$.

\begin{lemma}
\label{L4.2}
For any space of conditions $C$, we have
\[
  (J_C)^* = \Ann_{\cD}(C) = \cD(\cO, W_C) \qquad \text{and} \qquad J_C = \cD(W_C, \cO) .
\]
\end{lemma}

\begin{proof}
Since $J_C \supseteq \cD$, we have $(J_C)^* \subseteq \cD$. For $P \in \cD$ and $F \in
\widetilde{C}_\lambda$, the inclusion $\ml^{-1} F P \subseteq \cD$ holds if and only if
$F P \in \ml \cD$ (multiply by $\ml$, resp.\ by $\ml^{-1}$). Hence
\[
  (J_C)^* = \{ P \in \cD \, : \, \widetilde{C}_\lambda\, P \subseteq \ml \cD \ \text{ for all }
  \lambda \} = \Ann_{\cD}(C) ,
\]
and $\Ann_{\cD}(C) = \cD(\cO, W_C)$ by Lemma~\ref{L2.6}$(c)$. Since $J_C$ is reflexive, it follows
from Theorem~\ref{T2.1}$(c)$ that $J_C = (J_C)^{**} = \cD(\cO, W_C)^* = \cD(W_C, \cO)$.
\end{proof}

We can now state and prove the main result of this section, which is a more precise version of
Theorem~\ref{T1.1}.

\begin{theorem}
\label{T4.2}
\begin{enumerate}
\item[$(a)$] Let $I$ be a left ideal of $\cD$, $L \in I$ an operator of minimal order, and let
$C := \bigoplus_{\lambda} V_\lambda(I;L)$. Then $C$ is a space of conditions on $X$, and
\[
  I\, L^{-1} = J_C = \cD(W_C, \cO) , \qquad (I\, L^{-1})^* = \Ann_{\cD}(C) = \cD(\cO, W_C) .
\]
In particular, $I \cong \cD(W_C, \cO)$ as left $\cD$-modules.
\item[$(b)$] Conversely, let $C$ be any space of conditions on $X$, and let $L$ be a nonzero
element of $\cD$ such that $I := J_C\, L \subseteq \cD$. Then $L$ is an operator of minimal order
in $I$, and $V_\lambda(I; L) = C_\lambda$ for all $\lambda \in X$.
\item[$(c)$] If $L$ and $L' = h L$ $(h \in K^\times)$ are two operators of minimal order in $I$,
with associated spaces of conditions $C$ and $C'$, then $W_{C'} = h\, W_C$.
\end{enumerate}
\end{theorem}

\begin{proof}
$(a)$ By Theorem~\ref{T4.1} and Corollary~\ref{C3.1}$(a)$, $C$ is a finite-dimensional
homogeneous subspace of $\Delta(X)$. Multiplying \eqref{E42} by $L^{-1}$ on the right, we get
$I L^{-1} = J_C$, and the result follows from Lemma~\ref{L4.2}. $(b)$ Since $1 \in J_C \subset
\cD(K)$, we have $L \in I$ and $\ord(qL) \ge \ord(L)$ for all nonzero $q \in J_C$; hence $L$ has
minimal order in $I$. Let $C' := \bigoplus_\lambda V_\lambda(I;L)$. By part $(a)$,
$\cD(W_{C'}, \cO) = I L^{-1} = J_C = \cD(W_C, \cO)$. Passing to the dual ideals and applying
Theorem~\ref{T2.1}$(a)$, we get $W_{C'} = W_C$, and therefore $C' = C$ by Lemma~\ref{L2.6}$(b)$.
$(c)$ By part $(a)$, $\cD(W_{C'}, \cO) = I L^{-1} h^{-1} = \cD(W_C, \cO)\, h^{-1}$. Dualizing, we
get $\cD(\cO, W_{C'}) = h\, \cD(\cO, W_C)$, whence $W_{C'} = \cD(\cO, W_{C'}) \cdot \cO = h\,
W_C$, by Theorem~\ref{T2.1}$(a)$.
\end{proof}

Applying Theorem~\ref{T4.2} to closed ideals, we obtain the following description of the closure
of a differential operator, which includes Theorem~\ref{T3.2}$(a)$.

\begin{corollary}
\label{C4.1}
Let $L$ be a nonzero element of $\cD$. Let $W(L) := W_{V_\Lambda(L)}$ be the primary decomposable
subspace of $\cO$ associated to the space of conditions $V_\Lambda(L)$, and let
$M(L) = \Ann_{\cD}(V_\Lambda(L))$, see \eqref{E38}. Then
\begin{enumerate}
\item[$(a)$] $\cl(L) = \cD(W(L), \cO) \cdot L\,$ and $\,M(L) = (\cl(L)\, L^{-1})^* = \cD(\cO, W(L))$;
\item[$(b)$] $W(L) = \{ f \in \cO \, : \, f \in L(\hO_\lambda) \text{ for all } \lambda \in X \}$;
in other words, $W(L)$ consists of all regular functions $f$ on $X$ such that the inhomogeneous
equation $L(u) = f$ has a formal power series solution at every point of $X$. In particular,
$L(\cO) \subseteq W(L)$, and $\codim_{\cO} W(L) = \sum_{\lambda} \dim V_\lambda(L)$.
\end{enumerate}
\end{corollary}

\begin{proof}
$(a)$ follows from Theorem~\ref{T4.2}$(a)$ and Lemma~\ref{L4.1}$(a)$. $(b)$ By definition,
$W(L)$ is the set of all $f \in \cO$ such that $f \in V_\lambda(L)^\perp$ for all $\lambda$, and
$V_\lambda(L)^\perp = L(\hO_\lambda)$ by Corollary~\ref{C3.1}$(a)$. The last statement follows
from Lemma~\ref{L2.6}$(a)$.
\end{proof}

\begin{corollary}
\label{C4.2}
Let $L \in \cD$ and $h \in K^\times$ be such that $h L \in \cD$. Then $\cl(hL) = \cl(L)$ and
$W(hL) = h \cdot W(L)$. Moreover, for every $\lambda \in X$, we have $V_\lambda(hL)^\perp = h\,
V_\lambda(L)^\perp$ in $\hO_\lambda$; in particular, if $h$ is regular at $\lambda$ and
$k = \ord_\lambda(h)$, then
$\ \Sigma(V_\lambda(hL)) = \{0, 1, \ldots, k-1\} \sqcup (k + \Sigma(V_\lambda(L)))$.
\end{corollary}

\begin{proof}
The first two statements follow from Theorem~\ref{T4.2}$(c)$ applied to $I = \cl(L)$. The local
statement is a consequence of Lemma~\ref{L2.2}$(d)$ and Lemma~\ref{L2.3}$(c)$.
\end{proof}

\begin{remarks}
1. Corollary~\ref{C4.2} shows how the space of conditions attached to an ideal depends on the
choice of an operator of minimal order: for example, if $X = \A^1$, $a \in \Lambda(L)$ and
$L' = (x - a)^k L$, then the degree set of $V_a(L')$ is obtained from that of $V_a(L)$ by shifting
by $k$ and adding $\{0, 1, \ldots, k-1\}$, while $W(L') = (x-a)^k\, W(L)$.

2. By \cite[Section~5]{BerestWilson09}, the class of the fat right ideal $\cD(\cO, W)$ in
$K_0(\cD) \cong K_0(X)$ is determined by the line bundle $\prod_\lambda \ml^{k_\lambda}$, where
$W = \bigcap_\lambda W_\lambda$ and $k_\lambda = \codim_{\cO} W_\lambda$. Hence, by
Corollary~\ref{C4.1}, the right ideal $M(L)$ dual to $\cl(L)$ corresponds to the ideal
$\prod_{\lambda \in \Lambda(L)} \ml^{d_\lambda(L)}$ of $\cO$, where $d_\lambda(L) = \dim
V_\lambda(L)$, and $\gamma[\cl(L)]$ is the class of the inverse of this ideal in $\Pic(X)$.
\end{remarks}

Combining Theorem~\ref{T4.2} with Theorem~\ref{T2.1}, we may restate our classification in the
following form, which is dual to the description of $\cR(\cD)$ in terms of the adelic
Grassmannian (see Section~\ref{S2.6}).

\begin{corollary}
\label{C4.3}
Let $\mathrm{Gr}^{\mathrm{ad}}_*(X)$ denote the set of all spaces of conditions on $X$ modulo the
equivalence relation: $C \sim C'$ if and only if $W_{C'} = h\, W_C$ for some $h \in K^\times$.
Then the assignment $C \mapsto J_C$, see \eqref{E45}, induces a bijection
\[
  \mathrm{Gr}^{\mathrm{ad}}_*(X) \stackrel{\sim}{\longrightarrow} \cL(\cD) ,
\]
whose inverse is given by $[I] \mapsto [\,\bigoplus_\lambda V_\lambda(I;L)\,]$, where $L$ is any
operator of minimal order in $I$.
\end{corollary}

\begin{proof}
By Lemma~\ref{L4.2}, $J_C = \cD(W_C, \cO)$, so that the map is well defined and injective by
Theorem~\ref{T2.1} and Lemma~\ref{L2.6}$(b)$. It is surjective by Theorem~\ref{T4.2}$(a)$, and the
description of the inverse follows from Theorem~\ref{T4.2}.
\end{proof}

We call $\mathrm{Gr}^{\mathrm{ad}}_*(X)$ the \emph{dual adelic Grassmannian} of $X$. By
Lemma~\ref{L2.6}, it is in natural bijection with the set of equivalence classes of primary
decomposable subspaces of $\cO$, i.e.\ with the adelic Grassmannian of $X$ in the sense of
\cite{CanningsHolland1994, BerestWilson09}; the point of Corollary~\ref{C4.3} is that the
description in terms of conditions is the one adapted to left ideals: both the ideal $J_C$ and
the space of conditions of a given ideal are given by explicit formulas.

\section{The Calogero--Moser index of an ideal}
\label{S5}

In this section we prove Theorem~\ref{T1.2} and deduce our main theorem. We will prove
Theorem~\ref{T1.2} in the following, more precise form. Recall that, for a finite-dimensional
subspace $V \subseteq \Dl$, we write $\Sigma(V)$ for its degree set and $\nu(V) = \sum_{s \in
\Sigma(V)} s - \binom{\dim V}{2}$ (see Section~\ref{S2.2}).

\begin{theorem}
\label{T5.1}
Let $C = \bigoplus_\lambda C_\lambda$ be a space of conditions on $X$. Then the local
Calogero--Moser indices of the right ideal $\Ann_{\cD}(C) = \cD(\cO, W_C)$ are given by
\[
  n_\lambda(\Ann_{\cD}(C)) = \nu(C_\lambda) , \qquad \lambda \in X .
\]
Consequently,
\[
  n(\cD(W_C, \cO)) = n(\cD(\cO, W_C)) = \sum_{\lambda \in \Supp(C)} \nu(C_\lambda) .
\]
\end{theorem}

The proof proceeds in three steps. First, using the additivity \eqref{E211} of the
Calogero--Moser index and its compatibility with completion, we show that $n_\lambda(\Ann_{\cD}(C))$
depends only on the subspace $C_\lambda$ of $\Dl \cong \C[\del]$ (Lemma~\ref{L5.1}). This reduces
the problem to the case when $X$ is the affine line and $C$ is supported at one point. Second, we
treat this special case using Wilson's description of the adelic Grassmannian
(Proposition~\ref{P5.1}). Third, we pass from right ideals to left ideals with the help of
Proposition~\ref{P2.2}. An alternative proof of the second step, which is purely algebraic and
does not use the reduction to the affine line, is given in the Appendix.

\subsection{Localization}
\label{S5.1}
Let $\hD := \C[[t]][\del]$, with $\del\, t - t\, \del = 1$, and let $C_0$ be a finite-dimensional
subspace of $\C[\del] \cong \hD/t\, \hD$. Consider the right ideal
\[
  \Ann_{\hD}(C_0) := \{ P \in \hD \, : \, F P \in t\, \hD \ \text{ for all } F \in C_0 \} ,
\]
which contains $t^m \hD$ for $m \gg 0$, and define
\begin{equation}
\label{E51}
  n(C_0) := \dim_{\C}\bigl[\, (\gr \Ann_{\hD}(C_0))^{\vee\vee} \big/ \gr \Ann_{\hD}(C_0) \,\bigr] ,
\end{equation}
where $\gr$ is taken with respect to the order filtration on $\hD$ and
$(\,\text{--}\,)^\vee$ stands for the dual over $\gr \hD = \C[[t]][\xi]$.

\begin{lemma}
\label{L5.1}
Let $C$ be a space of conditions on $X$, and let $M := \Ann_{\cD}(C)$. Then, for every
$\lambda \in X$, the right ideal of $\hD_\lambda$ generated by $M$ is
\[
  \widehat{M}_\lambda = \{ P \in \hD_\lambda \, : \, F P = 0 \ \text{in}\ \Dl \ \text{for all } F
  \in C_\lambda \} .
\]
Consequently, $n_\lambda(M) = n(C_\lambda)$, where $C_\lambda$ is regarded as a subspace of
$\C[\del]$ via the isomorphisms $\hD_\lambda \cong \C[[t]][\del]$ and $\Dl \cong \C[\del]$
determined by a local parameter $t$ at $\lambda$. In particular, $n_\lambda(M) = 0$ for
$\lambda \notin \Supp(C)$.
\end{lemma}

\begin{proof}
Denote the right-hand side by $M'$, and choose $m$ so that $C_\mu \subseteq \Delta_\mu^{\le m-1}$
for all $\mu$. Then $M = \bigcap_\mu \Ann_{\cD}(C_\mu)$, and $\Ann_{\cD}(C_\mu) \supseteq \fm_\mu^m\, \cD$
for each $\mu$, by Proposition~\ref{P2.1}$(a)$ and Lemma~\ref{L2.1}. If
$\mu \ne \lambda$, then $\Ann_{\cD}(C_\mu)\, \hD_\lambda \supseteq \fm_\mu^m\, \hD_\lambda =
\hD_\lambda$, since $\fm_\mu^m$ contains functions that do not vanish at $\lambda$. As the
extension of right ideals from $\cD$ to $\hD_\lambda \cong \cD \otimes_{\cO} \hO_\lambda$ is
exact and commutes with finite intersections, we get $\widehat{M}_\lambda =
\Ann_{\cD}(C_\lambda)\, \hD_\lambda$. Clearly, $\Ann_{\cD}(C_\lambda)\, \hD_\lambda \subseteq M'$,
and both ideals contain $t^m \hD_\lambda = \ml^m \hD_\lambda$. Since the natural map
$\cD/\ml^m \cD \to \hD_\lambda/t^m \hD_\lambda$ is bijective and the condition $F P \in t
\hD_\lambda$ depends only on the class of $P$ modulo $t^m \hD_\lambda$, the two ideals have the
same image in $\hD_\lambda/t^m \hD_\lambda$. Hence they coincide. The second statement follows
from \eqref{E212} and \eqref{E51}; if $C_\lambda = 0$, then $\widehat{M}_\lambda = \hD_\lambda$,
and $n_\lambda(M) = 0$.
\end{proof}

\subsection{The one-point case}
\label{S5.2}
In view of Lemma~\ref{L5.1}, Theorem~\ref{T5.1} is equivalent to the following statement.

\begin{proposition}
\label{P5.1}
Let $C_0$ be a finite-dimensional subspace of $\C[\del]$ with degree set $S = \Sigma(C_0)$. Then
$\,n(C_0) = \nu(S)$. Explicitly, if $\{F_1, \ldots, F_g\}$ is a basis of $C_0$ with $\deg F_i =
s_i$ and $s_1 < \ldots < s_g$, then $\ n(C_0) = \sum_{i=1}^g s_i - \binom{g}{2}$.
\end{proposition}

We will use the following elementary fact about Wronskians (cf.\
\cite[Proposition~2]{Hermoso2024}).

\begin{lemma}
\label{L5.2}
Let $F_1, \ldots, F_g \in \C[z]$ be polynomials of pairwise distinct degrees $s_1, \ldots, s_g$.
Then $\Wr(F_1, \ldots, F_g)$ is a nonzero polynomial of degree $\sum_{i} s_i - \binom{g}{2}$.
Similarly, if $u_1, \ldots, u_g \in \C[[t]]$ are power series of pairwise distinct orders
$\rho_1, \ldots, \rho_g$, then $\Wr(u_1, \ldots, u_g)$ is a nonzero power series of order
$\sum_i \rho_i - \binom{g}{2}$.
\end{lemma}

\begin{proof}
Since the Wronskian is multilinear, it suffices to compute the contribution of the leading
terms: $\Wr(z^{s_1}, \ldots, z^{s_g}) = \det\bigl(\ff{s_j}{i-1}\, z^{s_j - i + 1}\bigr) =
z^{\sum_j s_j - \binom{g}{2}}\, \det\bigl(\ff{s_j}{i-1}\bigr)$. As $\ff{s}{i-1}$ is a monic
polynomial in $s$ of degree $i-1$, the last determinant equals the Vandermonde determinant
$\prod_{i < j} (s_j - s_i) \ne 0$. All other terms of $\Wr(F_1, \ldots, F_g)$ have lower degree
(resp., higher order, in the case of power series).
\end{proof}

\begin{proof}[Proof of Proposition~\ref{P5.1}]
We regard $C_0$ as a space of conditions on the affine line $\A^1$ (with coordinate $t$) supported
at the origin, and put $N := \Ann_{A_1}(C_0) = \cD(\C[t], W_{C_0})$. By Lemma~\ref{L5.1} and
\eqref{E211}, $n(N) = n_0(N) = n(C_0)$. On the other hand, by Proposition~\ref{P2.4}$(a)$, $n(N) =
- \res_\infty a_1(z)$, where $a_1$ is the first coefficient of the Baker function of the
corresponding point of $\Gra$. Since $C_0$ is supported at the origin, the functions $f_i$ of
Proposition~\ref{P2.4}$(b)$ are simply $f_i(z) = F_i(\del_t)\, e^{tz}|_{t=0} = F_i(z)$; hence
\[
  n(N) = \res_\infty\, \del_z \log \Wr(F_1, \ldots, F_g) = \deg \Wr(F_1, \ldots, F_g) ,
\]
because, for a polynomial $P(z)$, the coefficient of $z^{-1}$ in the expansion of $P'/P$ at
infinity is $\deg P$. The result follows now from Lemma~\ref{L5.2}.
\end{proof}

\subsection{Proofs of the main results}
\label{S5.3}

\begin{proof}[Proof of Theorem~\ref{T5.1}]
The first formula follows from Lemma~\ref{L5.1} and Proposition~\ref{P5.1}. Summing over
$\lambda$, we get $n(\cD(\cO, W_C)) = \sum_\lambda \nu(C_\lambda)$, by \eqref{E211}. Finally,
$\cD(\cO, W_C)$ is the dual of the fractional left ideal $\cD(W_C, \cO)$ (see
Theorem~\ref{T2.1}$(c)$), so that $n(\cD(W_C, \cO)) = n(\cD(\cO, W_C))$ by
Proposition~\ref{P2.2}.
\end{proof}

Theorem~\ref{T1.2} is a restatement of the last formula of Theorem~\ref{T5.1}, in view of
Lemma~\ref{L2.3}$(a)$.

\begin{proof}[Proof of Theorem~\ref{T3.2}]
Part $(a)$ is Corollary~\ref{C4.1}$(a)$. To prove $(b)$, we apply Theorem~\ref{T5.1} to the space
of conditions $C = V_\Lambda(L)$: for every $\lambda \in X$,
\[
  n_\lambda(M(L)) = \nu(V_\lambda(L)) = n_\lambda(L) ,
\]
where the second equality is Corollary~\ref{C3.1}$(c)$. Summing over $\lambda$ gives $n(M(L)) =
n(L)$. Finally, by part $(a)$, $M(L)$ is the dual of the fractional left ideal
$\cl(L) L^{-1} \cong \cl(L)$; hence $n(\cl(L)) = n(M(L))$ by Proposition~\ref{P2.2}. This proves
$(c)$, and thus completes the proof of Theorem~\ref{T1.3}.
\end{proof}

For an arbitrary (not necessarily closed) ideal, Theorems~\ref{T4.2} and~\ref{T5.1} give

\begin{corollary}
\label{C5.1}
Let $I$ be a left ideal of $\cD$, and let $L \in I$ be an operator of minimal order. Then
\[
  n(I) = \sum_{\lambda \in \Lambda(L)} \nu(V_\lambda(I;L)) .
\]
\end{corollary}

\begin{remark}
Formula \eqref{E19} is dual to the formula for the Calogero--Moser index given in
\cite[(8.3)]{BerestWilson04}. To see this, consider the case when $X = \A^1$ and $C$ is supported
at the origin and spanned by monomials: $C = \spn\{\del^{s_1}, \ldots, \del^{s_g}\}$ with
$s_1 < \ldots < s_g$. Then $W_C$ is spanned by the monomials $x^r$ with $r \notin \{s_1, \ldots,
s_g\}$. Let $r_0 < r_1 < r_2 < \ldots$ be the exponents occurring in $W_C$; thus $g$ is the number
of `gaps' in $W_C$, and $r_i = g + i$ for $i > s_g - g$. By Theorem~\ref{T5.1},
\[
  n = \sum_{i=1}^g s_i - \sum_{i=0}^{g-1} i = \sum_{k=g}^{s_g} k \, - \sum_{i=0}^{s_g - g} r_i =
  \sum_{i \ge 0} (g + i - r_i) ,
\]
which is precisely the formula of \cite{BerestWilson04}.
\end{remark}

\section{Explicit constructions and examples}
\label{S6}

In this section we illustrate our results. We begin with the case of the affine line, which is
both the simplest example and (via Proposition~\ref{P5.1}) an ingredient of the proof of the
general index formula.

\subsection{The affine line}
\label{S6.1}
Let $X = \A^1$, so that $\cO = \C[x]$ and $\cD = A_1$. Every line bundle on $X$ is trivial; hence,
by Proposition~\ref{P2.3}, a left ideal $M \subseteq A_1$ is cyclic if and only if $n(M) = 0$, and
a holonomic function on $\C$ is strongly holonomic if and only if its annihilator in $A_1$ is not
cyclic. For $\lambda \in \C$, we take $t = x - \lambda$ as a local parameter, so that $\del_t =
\del$, $\,\Dl = A_1/(x-\lambda) A_1 \cong \C[\del]$, and the pairing \eqref{E16} is
$\lan F(\del), f \ran_\lambda = [F(\del) f](\lambda)$. We write $\theta_\lambda := (x - \lambda)
\del$.

\subsubsection{Computing $V_\lambda(L)$}
The spaces $V_\lambda(L)$ and the cotype of an operator can be computed by linear algebra, as in
\cite{Tsai}. The following lemma is valid for any curve.

\begin{lemma}
\label{L6.1}
Let $L \in \cD$, $\lambda \in X$, and suppose that $w_\lambda(L) \le 0$, so that $L$ can be
written in $\hD_\lambda$ in the form
\begin{equation}
\label{E61}
  L = \sum_{i \ge 0} t^i\, q_i(\theta) , \qquad q_i \in \C[s] , \quad \deg q_i \le \ord(L) .
\end{equation}
Then, for all $j \ge 0$,
\begin{equation}
\label{E62}
  \del^j L \, \equiv \, \sum_{i=0}^{j} \ff{j}{i}\, q_i(j - i)\, \del^{j-i} \quad (\mathrm{mod}\ t\,
  \hD_\lambda) .
\end{equation}
Thus, in the basis $\{\del^j\}$ of $\Dl$, right multiplication by $L$ is given by a triangular
matrix $R_\lambda(L)$ with diagonal entries $q_0(j)$, and $V_\lambda(L) = \Ker R_\lambda(L)$. If
$L$ is $\lambda$-normalized, then $q_0 = q_{L,\lambda}$ is its indicial polynomial, and every
element of the cotype $J_\lambda(L)$ is a non-negative integer root of $q_0$.
\end{lemma}

\begin{proof}
The expansion \eqref{E61} follows from \eqref{E33}. Since $[\del^j, t] = j\, \del^{j-1}$, we have
$\del^j t^i \equiv \ff{j}{i}\, \del^{j-i}$ modulo $t\, \hD_\lambda$ (the right-hand side being
zero for $i > j$). On the other hand, $q(\theta)\, \del^m = \del^m q(\theta - m)$ for any
polynomial $q$; hence $\del^m q(\theta) = q(\theta + m)\, \del^m \equiv q(m)\, \del^m$, because
$\theta\, \del^m \in t\, \hD_\lambda$. Combining these two formulas gives \eqref{E62}. The last
statement follows from Theorem~\ref{T3.1} (cf.\ also Lemma~\ref{L3.3}$(a)$).
\end{proof}

For $X = \A^1$ and $L \in A_1$, the sum in \eqref{E61} is finite: it is the decomposition of
$L$ into homogeneous components with respect to the grading of $A_1$ defined by the adjoint
action of $\theta_\lambda$. The operator $(x - \lambda)^{w_\lambda(L)} L$ is called in
\cite{Tsai} the \emph{balanced form} of $L$ at $\lambda$. By \eqref{E44}, the kernels of the
matrices $R_\lambda(L)$, $\lambda \in \Lambda(L)$, give a finite set of generators of $\cl(L)$,
and, by Theorem~\ref{T3.2} and Corollary~\ref{C3.1}$(c)$, they determine the Calogero--Moser index
of $\cl(L)$: if $\{F_{\lambda, 1}, \ldots, F_{\lambda, g_\lambda}\}$ is a basis of $V_\lambda(L)$
consisting of polynomials of distinct degrees, then
\begin{equation}
\label{E63}
  n(\cl(L)) = \sum_{\lambda \in \Lambda(L)} \Bigl[\, \sum_{i=1}^{g_\lambda} \deg F_{\lambda, i} -
  \binom{g_\lambda}{2} \Bigr] .
\end{equation}
If the right-hand side of \eqref{E63} is zero, then $\cl(L)$ is cyclic; in fact, in this case
$\cl(L) = A_1\, h^{-1} L$, where $h \in \C[x]$ is the greatest common divisor of the coefficients
of $L$.

\begin{example}
\label{Ex6.1}
Let $L = x \del - n$ with $n \in \Z_{\ge 0}$. Then $L = \theta_0 - n$ is normalized at $\lambda =
0$, with indicial polynomial $q_0(s) = s - n$, and $\del^j L \equiv (j - n)\, \del^j$ modulo
$x A_1$. Hence $\Lambda(L) = \{0\}$, $V_0(L) = \C\, \del^n$, the cotype of $L$ at $0$ is $\{n\}$,
and $n(L) = n$. Formula \eqref{E44} gives
\[
  \cl(L) = A_1 (x \del - n) + A_1\, x^{-1} \del^n (x \del - n) = A_1 (x\del - n) + A_1\,
  \del^{n+1} = \Ann(x^n) ,
\]
in agreement with \eqref{E11}, and we recover the fact that $n(\Ann(x^n)) = n$. The primary
decomposable subspace of Corollary~\ref{C4.1} is $W(L) = \{ f \in \C[x] : f^{(n)}(0) = 0 \}$; it
has exactly one gap, in degree $n$, and $M(L) = \{P \in A_1 : \del^n P \in x A_1\}$. Dually, $L$
acts on $\C[[x]]$ by $x^j \mapsto (j - n)\, x^j$, and $L(\C[[x]])$ is the annihilator of $\del^n$,
as it should be by Lemma~\ref{L2.2}$(d)$.
\end{example}

\begin{example}
\label{Ex6.2}
Let $f(x) = x^2 (x-1)^3$. An operator of minimal order in $\Ann(f)$ is $L = x(x-1)\del - 5x +
2$, with $\Lambda(L) = \{0, 1\}$. At $\lambda = 0$, we have $L = -(\theta_0 - 2) + x\,(\theta_0
- 5)$, and \eqref{E62} gives $\del^j L \equiv (2 - j)\, \del^j + j (j - 6)\, \del^{j-1}$; hence
$V_0(L) = \C\,(\del^2 + 8 \del + 20)$. At $\lambda = 1$, we have $L = (\theta_1 - 3) + (x-1)
(\theta_1 - 5)$, so that $\del^j L \equiv (j - 3)\, \del^j + j (j-6)\, \del^{j-1}$ modulo
$(x-1) A_1$, and $V_1(L) = \C\,(\del^3 - 9 \del^2 + 36 \del - 60)$. Thus $n(\Ann(f)) = 2 + 3 = 5 =
\deg f$.
\end{example}

\subsubsection{Calogero--Moser data}
By Corollary~\ref{C4.3}, on any curve $X$ the space of conditions of a left ideal determines
not only its Calogero--Moser index but the ideal class itself. For $X = \A^1$, this can be made
completely explicit: from the space of conditions one can compute the Calogero--Moser data of the
ideal, namely the element $\kappa$ defined in \eqref{E214}, which is a complete isomorphism
invariant of the ideal (see Section~\ref{S2.6}). Let $M$ be a left ideal of $A_1$, $L \in M$ an
operator of minimal order, and let $C = \bigoplus_\lambda V_\lambda(M; L)$. By
Theorem~\ref{T4.2}, $M \cong \cD(W_C, \C[x])$, so that the point $W \in \Gra$ corresponding to
$M$ is $W = \prod_\lambda (x - \lambda)^{-\dim C_\lambda}\, W_C$. Choose a basis $\{c_1, \ldots,
c_g\}$ of $C$ with $c_i \in C_{\lambda_i}$ represented by $F_i(\del) \in \C[\del]$. Then, by
\eqref{E215} and Proposition~\ref{P2.4}$(c)$, the element $\kappa = \kappa_M$ relating the two
distinguished representatives of $M$ (see \eqref{E214}) is given by
\begin{equation}
\label{E64}
  \kappa^{-1}(x, \del)\,[e^{xz}] = \kappa^{-1}(x,z)\, e^{xz} , \qquad \kappa^{-1}(x, z) =
  \frac{\det \mathsf{M}(z, x)}{\det \mathsf{M}(z, \infty)} ,
\end{equation}
where $\mathsf{M}(z,x)$ is the $g \times g$ matrix with entries $\mathsf{m}_{ij} = (\del_z +
\lambda_i)^{j-1} [F_i(z) - (x - \lambda_i)^{-1} F_i'(z)]$. In other words, to find
$\kappa^{-1}(x, \del)$ one expands the rational function $\kappa^{-1}(x,z)$ in the form
$1 + \sum_k \phi_k(x)\, \chi_k(z)$ and replaces $z$ by $\del$, writing the functions of $x$ on
the left.

\begin{example}
\label{Ex6.3}
Suppose that $C$ consists of a single condition $F(\del)$ at the point $\lambda$, so that
$M \cong J_C = A_1 + A_1 (x - \lambda)^{-1} F(\del)$. Then \eqref{E64} gives
\begin{equation}
\label{E65}
  \kappa^{-1}(x, \del) = 1 - \frac{1}{x - \lambda}\, \frac{F'(\del)}{F(\del)} .
\end{equation}
In this case the distinguished representatives of $M$ are easy to find directly: $M_x =
J_C\, (x - \lambda) = A_1 (x - \lambda) + A_1\,[F(\del) + (x-\lambda)^{-1} F'(\del)]$, and
\[
  M_\del = M_x\, \kappa^{-1} = A_1 \Bigl( x - \lambda - \frac{F'(\del)}{F(\del)} \Bigr) + A_1\,
  F(\del) .
\]
For $F = \del^n$ and $\lambda = 0$, formula \eqref{E65} gives $\kappa^{-1} = 1 - n\, x^{-1}
\del^{-1}$, in agreement with \eqref{E216} for $f = x^n$. Note that the coefficient of $x^{-1}
z^{-1}$ in the expansion of $\kappa^{-1}(x,z)$ at infinity is $- \deg F$, in agreement with
Theorem~\ref{T5.1} and Proposition~\ref{P2.4}$(a)$.
\end{example}

\begin{example}
\label{Ex6.4}
Suppose that $C = \spn\{F_1(\del), F_2(\del)\}$ is supported at one point $\lambda$. A
straightforward calculation with \eqref{E64} gives
\[
  \kappa^{-1}(x, \del) = 1 - \frac{1}{x - \lambda}\, \frac{\Wr'}{\Wr}(\del) + \frac{1}{(x -
  \lambda)^2}\, \frac{F_1' F_2'' - F_1'' F_2'}{\Wr}(\del) ,
\]
where $\Wr = \Wr(F_1, F_2) = F_1 F_2' - F_1' F_2$. Thus the coefficients of $\kappa^{-1}$ at the
powers of $(x - \lambda)^{-1}$ coincide with the coefficients of the monic differential operator
$\del_z^2 - (\Wr'/\Wr)\, \del_z + (F_1' F_2'' - F_1'' F_2')/\Wr$ of order $2$ in the variable $z$
whose kernel is spanned by $F_1(z)$ and $F_2(z)$. This is a manifestation of bispectral duality
(cf.\ Proposition~\ref{P2.4}$(b)$).
\end{example}

\begin{example}
\label{Ex6.5}
Let $f = x^2 (x-1)^3$, as in Example~\ref{Ex6.2}. Then $C$ is spanned by $F_1 = \del^2 + 8 \del +
20$ at $\lambda_1 = 0$ and $F_2 = \del^3 - 9 \del^2 + 36 \del - 60$ at $\lambda_2 = 1$. Computing
the determinant in \eqref{E64}, we find $\det \mathsf{M}(z,x) = z^5 - (2/x + 3/(x-1))\, z^4$ and
$\det \mathsf{M}(z, \infty) = z^5$, whence
\[
  \kappa^{-1} = 1 - \Bigl( \frac{2}{x} + \frac{3}{x-1} \Bigr)\, \del^{-1} = 1 - \frac{f'}{f}\,
  \del^{-1} ,
\]
in agreement with \eqref{E216}.
\end{example}

\begin{example}
\label{Ex6.6}
Let $p \in \C[x]$, and consider the function $f = e^{p(x)}\, x^n$, $n \ge 1$. It is annihilated by
$L = x \del - n - x\, p'(x)$, which is normalized at $0$, with $q_0(s) = s - n$ and $q_i(s) = -
p_{i}$ for $i \ge 1$, where $p'(x) = \sum_{i \ge 1} p_i\, x^{i-1}$. Hence $V_0(L)$ is spanned by a
single polynomial $F(\del) = \del^n + \ldots$ of degree $n$, and $n(\Ann(f)) = n$ for all $p$.
However, the isomorphism class of $\Ann(f)$ depends on $p$: by Corollary~\ref{C4.1}, $\Ann(f)
\cong \cD(W, \C[x])$, where $W = \{ g \in \C[x] : [F(\del) g](0) = 0\} \supseteq L(\C[x])$, and
two such ideals are isomorphic if and only if the corresponding polynomials $F$ are proportional
(by Theorem~\ref{T2.1}, since two subspaces of codimension one in $\C[x]$, both primary at the
origin, can be equivalent only if they coincide). Now, $F$ is proportional to $\del^n$ if and only if the coefficient of
$x^n$ vanishes in $L(x^j) = (j - n)\, x^j - x^{j+1} p'(x)$ for all $j < n$, i.e.\ if and only if
$p_1 = \ldots = p_n = 0$. Thus
\[
  \Ann(e^{p(x)} x^n) \cong \Ann(x^n) \quad \Longleftrightarrow \quad p \in \C + x^{n+1}\, \C[x] .
\]
In other words, the stabilizer of the class of $\Ann(x^n)$ in the group $\{e^{\mathrm{ad}\, p(x)}
: p \in \C[x]\} \subset \Aut(A_1)$ consists of the automorphisms with $p \in \C + x^{n+1} \C[x]$.
The same conclusion can be read off from the element $\kappa$, which determines the
isomorphism class of an ideal (see Section~\ref{S2.6}). For example, for $n = 2$, one finds
$F = (\del - p_1)^2 - p_2$, and \eqref{E65} gives
\[
  \kappa^{-1} = 1 - \frac{2}{x}\, \frac{\del - p_1}{(\del - p_1)^2 - p_2} .
\]
Comparing this with $\kappa^{-1} = 1 - 2\, x^{-1} \del^{-1}$ for $\Ann(x^2)$,
we see again that $\Ann(e^{p} x^2) \cong \Ann(x^2)$ if and only if $p_1 = p_2 = 0$, i.e.\ if and
only if $e^{\mathrm{ad}\, p}$ stabilizes the class of $\Ann(x^2)$.
\end{example}

\subsection{Operators of order one}
\label{S6.2}
We now return to an arbitrary curve $X$ and prove Theorem~\ref{T1.4}.

\begin{theorem}
\label{T6.1}
Let $f$ be a nonzero rational function on $X$, or a nonzero holomorphic function on $X$
annihilated by an operator of order one in $\cD$. Then $\Ann_{\cD}(f) = \cl(L)$ for any nonzero
operator $L \in \cD$ of order one such that $L(f) = 0$; for every $\lambda \in X$,
\[
  n_\lambda(L) = \max\{\ord_\lambda(f), 0\} ,
\]
and therefore
\[
  n(\Ann_{\cD}(f)) = \sum_{\lambda \in X} \max\{\ord_\lambda(f), 0\} .
\]
\end{theorem}

\begin{proof}
If $f \in K^\times$, choose a nonzero derivation $\delta \in \Der(\cO)$ and a nonzero $g \in \cO$
such that $g\, \delta(f)/f \in \cO$; then $L = g\, \delta - g\, \delta(f)/f$ is an operator of
order one annihilating $f$. In either case, $\Ann_{\cD}(f)$ is closed by Lemma~\ref{L2.4}$(a)$,
and it contains no nonzero operators of order zero; hence $\Ann_{\cD}(f) = \cl(L)$ by
Lemma~\ref{L2.4}$(b)$. Fix $\lambda \in X$ and a local parameter $t$ at $\lambda$. We can write
$f = t^m u$, where $m = \ord_\lambda(f) \in \Z$ and $u$ is an invertible element of
$\hO_\lambda$. If $L = a\, \del + b$ with $a, b \in \cO_\lambda$, then $b = - a f'/f$, so that in
$\hD_\lambda$
\[
  L = a\, t^{-1}\, (\theta - m - t\, c) , \qquad c := u'/u \in \hO_\lambda .
\]
Hence the $\lambda$-normalization of $L$ is $\bar L = \theta - m - t\,c$ (up to a unit), with
indicial polynomial $q(s) = s - m$. By Lemma~\ref{L3.2}, the matrix of $\bar L$ acting on
$\cO_\lambda/\nl^{j+1}$ in the basis $\{1, t, \ldots, t^j\}$ is triangular, with diagonal entries
$-m, 1 - m, \ldots, j - m$. If $m < 0$, this matrix is invertible for all $j$, so the cotype of
$L$ at $\lambda$ is empty. If $m \ge 0$, it is invertible for $j < m$ and has a one-dimensional
kernel for $j \ge m$; hence $J_\lambda(L) = \{m\}$. In both cases, $n_\lambda(L) = \max\{m,
0\}$. The last formula follows now from Theorem~\ref{T3.2}$(c)$.
\end{proof}

Thus, a function $f$ as in Theorem~\ref{T6.1} is strongly holonomic if and only if it has a zero
on $X$, and $n(\Ann_{\cD}(f))$ is the degree of the divisor of zeros of $f$ on $X$. In particular,
if $f \in \cO$ defines a finite morphism $f: X \to \A^1$, then $n(\Ann_{\cD}(f)) = \deg(f)$.

\begin{remark}
The assumption that $f$ is annihilated by an operator of order one cannot be dropped: in general,
the Calogero--Moser index of $\Ann_{\cD}(f)$ is not determined by the zeros of $f$. For example,
the function $\cos x$ has infinitely many zeros in $\C$, but $\Ann(\cos x) = A_1 (\del^2 + 1)$ is
a cyclic ideal, so that $n(\Ann(\cos x)) = 0$ (see Example~\ref{Ex6.8}$(2)$ below). Similarly, for $f = (e^x - 1)/x$, an operator of minimal order in
$\Ann(f)$ is
\[
  L = \del\, (\del - 1)\, x = x\, \del^2 + (2 - x)\, \del - 1 .
\]
Its normalization at the origin is $x L = \theta_0 (\theta_0 + 1) - x\, (\theta_0 + 1)$, with
indicial polynomial $s (s+1)$; hence, by Lemma~\ref{L6.1}, the cotype of $L$ at $0$ is $\{0\}$, and
$L$ has no singular points. Thus $n(L) = 0$, and $\Ann(f) = A_1 L$ is cyclic, although $f$ has
infinitely many zeros. The correct generalization of Theorem~\ref{T6.1} to operators of higher
order involves the Wronskian of a basis of solutions (see Corollary~\ref{C6.1} below): for
example, $\Wr(\cos x, \sin x) = 1$ has no zeros.
\end{remark}

\begin{example}
\label{Ex6.7}
Let $X = \A^1$.

$(1)$ If $f = p/q$ is a rational function with $p, q \in \C[x]$ coprime, then $n(\Ann(f)) =
\deg(p)$. In particular, $n(\Ann(p)) = \deg(p)$ for any polynomial $p$, while $\Ann(1/q)$ is a
cyclic ideal for any $q$: explicitly, $\Ann(1/q) = A_1\, h^{-1} \del\, q$, where $h$ is the
greatest common divisor of $q$ and $q'$ (in particular, $\Ann(1/q) = A_1\, \del\, q$ if and only
if $q$ has simple roots).

$(2)$ If $f = p\, e^{q}$ is a quasi-polynomial, with $p, q \in \C[x]$, then $n(\Ann(f)) =
\deg(p)$. This also follows from the fact that $\Ann(p\, e^{q})$ is the image of $\Ann(p)$ under
the automorphism $e^{\mathrm{ad}\, q}$ of $A_1$ (see Section~\ref{S2.6}).

$(3)$ The arguments of Theorem~\ref{T6.1} apply to any (multivalued) function annihilated by an
operator of order one. For example, the operator $L = x \del - \alpha$ annihilating $x^\alpha$ has
indicial polynomial $s - \alpha$ at the origin; hence, if $\alpha \notin \Z_{\ge 0}$, then
$n(L) = 0$ and $\cl(L) = A_1 (x \del - \alpha)$ is cyclic.
\end{example}

\subsection{Apparent singularities and Wronskians}
\label{S6.3}
Theorem~\ref{T6.1} has the following generalization to operators of higher order.

\begin{proposition}
\label{P6.1}
Let $L \in \cD$ be an operator of order $k \ge 1$, and let $\lambda \in X$. Assume that the
equation $L(u) = 0$ has $k$ linearly independent solutions in $\hO_\lambda$. Choose a basis
$u_1, \ldots, u_k$ of the space of solutions with $\rho_i := \ord_\lambda(u_i)$ pairwise
distinct. Then $J_\lambda(L) = \{\rho_1, \ldots, \rho_k\}$, and
\[
  n_\lambda(L) = \sum_{i=1}^k \rho_i - \binom{k}{2} = \ord_\lambda \Wr(u_1, \ldots, u_k) ,
\]
where the Wronskian is taken with respect to the derivation $\del_t$ (its order at $\lambda$ does
not depend on the choice of local parameter).
\end{proposition}

\begin{proof}
Let $\bar L$ be the $\lambda$-normalization of $L$, and let $q$ be its indicial polynomial. If
$u = c\, t^\rho + \ldots$ with $c \ne 0$, then $\bar L(u) = c\, q(\rho)\, t^\rho + \ldots$ by
Lemma~\ref{L3.2}; hence $q(\rho_i) = 0$ for all $i$. Since $\deg q \le k$, the numbers $\rho_1,
\ldots, \rho_k$ are all the roots of $q$. On the other hand, by Lemma~\ref{L3.3}, the cotype
$J_\lambda(L)$ is a subset of the set of roots of $q$ consisting of $d_\lambda = \dim \Ker(L;
\hO_\lambda) = k$ elements. Hence $J_\lambda(L) = \{\rho_1, \ldots, \rho_k\}$. The second equality
is Lemma~\ref{L5.2}.
\end{proof}

A point $\lambda \in \Lambda(L)$ satisfying the assumption of Proposition~\ref{P6.1} is called an
\emph{apparent singularity} of $L$. In classical terms, Proposition~\ref{P6.1} says that, at an
apparent singularity, the cotype of $L$ is the set of local exponents of $L$, and the local
holonomicity index measures the deviation of these exponents from the generic values
$0, 1, \ldots, k-1$.

\begin{corollary}
\label{C6.1}
Let $f_1, \ldots, f_k$ be linearly independent holomorphic functions on $X$ such that the ideal
$I = \Ann_{\cD}(f_1, \ldots, f_k)$ contains an operator of order $k$. Then
\[
  n(I) = \sum_{\lambda \in X} \ord_\lambda \Wr(f_1, \ldots, f_k) ,
\]
i.e.\ $n(I)$ equals the number of zeros of the Wronskian of $f_1, \ldots, f_k$ on $X$ counted with
multiplicities. For $X = \A^1$, if $L = a_k \del^k + a_{k-1} \del^{k-1} + \ldots \in I$ is an
operator of order $k$, then $n(I) = - \res_\infty (a_{k-1}/a_k)$; in particular, if $f_1, \ldots,
f_k$ are polynomials, then $n(I) = \deg \Wr(f_1, \ldots, f_k)$.
\end{corollary}

\begin{proof}
Since $f_1, \ldots, f_k$ are linearly independent, $I$ contains no nonzero operators of order less
than $k$; hence $I = \cl(L)$ for any $L \in I$ of order $k$ (Lemma~\ref{L2.4}), and $n(I) = n(L)$
by Theorem~\ref{T3.2}. At each point $\lambda \in X$, the Taylor expansions of $f_1, \ldots, f_k$
form a basis of solutions of $L(u) = 0$ in $\hO_\lambda$. Replacing this basis by one with distinct
orders does not change the order of the Wronskian, and the first formula follows from
Proposition~\ref{P6.1}. For $X = \A^1$, Abel's formula $\Wr' = - (a_{k-1}/a_k)\, \Wr$ shows that
the number of zeros of the entire function $\Wr = \Wr(f_1, \ldots, f_k)$ is the sum of the
residues of the rational function $- a_{k-1}/a_k$, i.e.\ $- \res_\infty (a_{k-1}/a_k)$.
\end{proof}

\begin{remark}
The second formula of Corollary~\ref{C6.1} may be viewed as a consequence of the classical Fuchs
relation: if $\lambda$ is a regular singular point of the operator $\del^k + p_1(x) \del^{k-1} +
\ldots + p_k(x)$, with local exponents $\rho_1(\lambda), \ldots, \rho_k(\lambda)$, then $\sum_i
\rho_i(\lambda) = \binom{k}{2} - \Res_{x = \lambda}\, p_1(x)$. It is instructive to compare this
formula with Proposition~\ref{P2.4}: in both cases the Calogero--Moser index appears as (minus)
the residue at infinity of the logarithmic derivative of a Wronskian.
\end{remark}

\begin{example}
\label{Ex6.8}
Let $X = \A^1$.

$(1)$ Let $0 \le k_1 < \ldots < k_s$ be integers. Then $L = (x \del - k_1) \cdots (x \del -
k_s)$ is an operator of minimal order in $I = \Ann(x^{k_1}, \ldots, x^{k_s})$; it is normalized
at $0$, and $V_0(L) = \spn\{\del^{k_1}, \ldots, \del^{k_s}\}$. Hence $n(I) = \sum_i k_i -
\binom{s}{2}$. For example, $\Ann(x^2, x^4) = A_1 \del^5 + A_1 (x \del - 2)(x \del - 4) + A_1 (x
\del^4 - \del^3)$ has Calogero--Moser index $5$. More generally, if $f_1, \ldots, f_s$ are
polynomials of pairwise distinct degrees, then $n(\Ann(f_1, \ldots, f_s)) = \sum_i \deg(f_i) -
\binom{s}{2}$, by Corollary~\ref{C6.1} and Lemma~\ref{L5.2}.

$(2)$ Let $f = \sum_{i=1}^k c_i\, e^{a_i x}$, with distinct $a_i \in \C$ and $c_i \ne 0$. Then
$L = (\del - a_1) \cdots (\del - a_k)$ is an operator of minimal order annihilating $f$. Since
$L$ has constant leading coefficient, $\Lambda(L) = \varnothing$, and $\Ann(f) = A_1 L$ is cyclic.
Thus sums of exponentials (for example, $\cos x$) are holonomic but not strongly holonomic. 
The same conclusion holds whenever a {\it minimal} order annihilator of $f$ is monic, i.e. has constant nonzero leading coefficient (for example, for the classical Airy function $f = {\rm Ai}(x)$ satisfying $ f'' - x f = 0$).

$(3)$ Let $f = x + x e^x$. An operator of minimal order in $\Ann(f)$ is $L = x^2 \del^2 - (x^2 +
2x) \del + (x + 2)$, whose space of solutions is spanned by $x$ and $x e^x$. Since $\Wr(x, x e^x)
= x^2 e^x$, Corollary~\ref{C6.1} gives $n(\Ann(f)) = 2$.

$(4)$ The Calogero--Moser index is not additive. For instance, let $f = e^{p} + q$ with $p, q \in
\C[x]$ non-constant, and let $h := q' - p' q$. Then $L = (h \del - h')(\del - p')$ is an operator
of minimal order in $\Ann(f)$, with solutions $e^p$ and $q$, and $\Wr(e^p, q) = h\, e^p$. Hence
$n(\Ann(f)) = \deg(h) = \deg(p) + \deg(q) - 1$, while $n(\Ann(e^p)) = 0$ and $n(\Ann(q)) =
\deg(q)$.
\end{example}

When the singular points of $L$ are not apparent, the cotype has to be computed directly, as in
Section~\ref{S6.1}. The following example is a corrected version of \cite[Example~4.5]{Tsai}.

\begin{example}
\label{Ex6.9}
Consider the Gauss hypergeometric operator $L = x(1-x) \del^2 + (c - (a+b+1) x) \del - ab$, with
$a, b, c \in \C$. Its singular points are contained in $\{0, 1\}$. At $\lambda = 0$, the
normalization of $L$ is
\[
  x L = \theta_0\, (\theta_0 + c - 1) - x\, (\theta_0 + a)(\theta_0 + b) ,
\]
with indicial polynomial $q_0(s) = s (s + c - 1)$. By Lemma~\ref{L6.1}, $\del^j (xL) \equiv
q_0(j)\, \del^j - j\, (j - 1 + a)(j - 1 + b)\, \del^{j-1}$. It follows that $0$ always belongs to
the cotype of $L$ at $0$, and the second root $m_0 := 1 - c$ of $q_0$ belongs to the cotype if
and only if $m_0$ is a positive integer and $\prod_{j=0}^{m_0-1} (j + a)(j + b) = 0$. Hence
\[
  n_0(L) = \begin{cases} m_0 - 1 & \text{if } m_0 \in \{2, 3, \ldots\} \text{ and } \{a, b\} \cap
  \{0, -1, \ldots, 1 - m_0\} \ne \varnothing , \\ 0 & \text{otherwise} . \end{cases}
\]
Similarly, at $\lambda = 1$, we have $-(x-1) L = \theta_1 (\theta_1 + a + b - c) + (x - 1)
(\theta_1 + a)(\theta_1 + b)$; the indicial polynomial is $s (s - m_1)$ with $m_1 := c - a - b$,
and the same argument gives
\[
  n_1(L) = \begin{cases} m_1 - 1 & \text{if } m_1 \in \{2, 3, \ldots\} \text{ and } \{a, b\} \cap
  \{0, -1, \ldots, 1 - m_1\} \ne \varnothing , \\ 0 & \text{otherwise} . \end{cases}
\]
Thus,  $\,n(L) = n_0(L) + n_1(L)\,$, where $n_0(L)$ and $n_1(L)$ are computed by the above formulas with $m_0 =  1 -c $ and $ m_1 = c - a - b $. This answers the question about  $\,n(\Ann(f))\,$ for the hypergeometric series $ f = {}_{2}F_1(a,b,c; x) $ posed in the Introduction,
provided the Gauss hypergeometric operator $L$ is a minimal order annihilator of $f$. In the resonant cases, the answer may differ. For a concrete example, take 
$ f = {}_{2}F_1(-1,1,3; x) = 1-x/3 $, in which case $\,n(\Ann(f)) = 1 $ by Theorem~\ref{T6.1}, while $n(L) = n_0(L) + n_1(L) = 0 + 2 = 2 $.
\end{example}

\subsection{Ideals that are not closed}
\label{S6.4}
The next example illustrates Corollary~\ref{C5.1}.

\begin{example}
\label{Ex6.10}
Let $I := A_1\, x^{n+1} + A_1 (x \del + n + 1)$, $n \ge 1$. This is the image of $\Ann(x^n)$ under
the Fourier automorphism $x \mapsto \del$, $\del \mapsto -x$ of $A_1$; hence $n(I) = n$ (see
Section~\ref{S2.6}). An operator of minimal order in $I$ is $L = x^{n+1}$, so that $\cl(L) = A_1$,
and we have a chain of left ideals $A_1 x^{n+1} \subset I \subset A_1$, whose first and last terms
are cyclic. Since $\del^j x^{n+1} \in x A_1$ if and only if $j \le n$, we have $V_0(L) =
\spn\{1, \del, \ldots, \del^n\}$ and $\nu(V_0(L)) = 0$, as it should be. On the other hand, for
$1 \le j \le n$,
\[
  \del^j x^{n+1} = x^{n+1-j} \prod_{i = n+2-j}^{n+1} (x \del + i) \, \in \, x\, A_1\, (x \del +
  n + 1) \subseteq x I ,
\]
while $x^{n+1} \notin x I$, since $x^n \notin I$. Hence $V_0(I; L) = \spn\{\del, \del^2, \ldots,
\del^n\}$, and Corollary~\ref{C5.1} gives $n(I) = \sum_{j=1}^n j - \binom{n}{2} = n$. Thus the three
ideals $A_1 L \subset I \subset \cl(L)$ correspond to the subspaces $0 \subset V_0(I;L) \subset
V_0(L)$, with indices $0$, $n$ and $0$, respectively.
\end{example}

\subsection{Other curves}
\label{S6.5}

\subsubsection{Open subsets}
The holonomicity index is local in nature: if $U \subseteq X$ is an affine open subset, and
$L \in \cD(X)$ is regarded as an element of $\cD(U)$, then $n(L|_U) = \sum_{\lambda \in U}
n_\lambda(L)$. For example, let $X = \A^1 \setminus \{0\}$, so that $\cD(X) = \C[x, x^{-1}][\del]$.
Then $\Ann_{\cD(X)}(x^n) = \cD(X)\, (\del - n/x)$ is cyclic, and $n = 0$, while $n(\Ann_{A_1}(x^n))
= n$: removing the only singular point of $x \del - n$ kills the index. Similarly, by
Theorem~\ref{T6.1}, a rational function $f$ regarded as a function on the complement of its
zeros in $X$ is never strongly holonomic.

\begin{remark}
The local nature of the holonomicity index suggests that it should be of interest for curves
which are not necessarily affine. Indeed, it is immediate from the definition that $n_\lambda(hL)
= n_\lambda(L)$ for any $h \in K^\times$; hence, for a holonomic function $f$, the number
$n_\lambda(f) := n_\lambda(L)$, where $L$ is an operator of minimal order annihilating $f$, is
well defined at every point $\lambda$ of any smooth curve $Y$ on which $f$ is defined, and we may
put $n_Y(f) := \sum_{\lambda \in Y} n_\lambda(f)$. When $Y$ is affine, $n_Y(f)$ is the
Calogero--Moser index of $\Ann_{\cD(Y)}(f)$, by Theorem~\ref{T3.2}. When $Y$ is projective, the
ring of global differential operators on $Y$ is too small, and the classification of ideals used
in this paper is not available; nevertheless, the invariant $n_Y(f)$ still makes sense. For
example, the function $x^{-n}$ $(n \ge 1)$ has a zero of order $n$ at the point $\infty \in
\mathbb{P}^1$, so that $n_{\mathbb{P}^1}(x^{-n}) = n$, while $n_{\A^1}(x^{-n}) = 0$. More
generally, the proof of Theorem~\ref{T6.1} shows that $n_Y(f) = \deg(f)$ for any non-constant
rational function $f$ on a smooth projective curve $Y$. It would be interesting to relate these
invariants to the classification of torsion-free $\cD$-modules on projective curves given in \cite{BenZviNevins}. 
\end{remark}

\subsubsection{Plane curves}
Let $X \subset \A^2$ be a smooth plane curve defined by an equation $P(x,y) = 0$, so that
$\cO = \C[x,y]/(P)$. In this case the ring $\cD(X)$ has a simple description (cf.\
\cite[Section~6.2.3]{BerestChalykhIdeals}).

\begin{lemma}
\label{L6.2}
The $\cO$-module $\Der(\cO)$ is free of rank one, generated by the derivation $\del_\theta$
defined by $\del_\theta(x) = P_y$, $\del_\theta(y) = - P_x$, where $P_x$ and $P_y$ are the partial
derivatives of $P$. Consequently,
\[
  \cD(X) \cong \cO\lan \del_\theta \ran \big/ \bigl( [\del_\theta, x] = P_y , \ [\del_\theta, y] =
  - P_x \bigr) .
\]
\end{lemma}

\begin{proof}
A derivation $\delta$ of $\cO$ is determined by the pair $(u, v) = (\delta(x), \delta(y)) \in
\cO^2$, subject to the relation $u P_x + v P_y = 0$ in $\cO$. Since $X$ is smooth, the functions
$P_x$ and $P_y$ have no common zeros on $X$; hence $\alpha P_x + \beta P_y = 1$ for some $\alpha,
\beta \in \cO$. Put $h := \beta u - \alpha v$. Then $u = u\,(\alpha P_x + \beta P_y) = - \alpha v
P_y + \beta u P_y = h\, P_y$, and similarly $v = - h\, P_x$. Thus $\delta = h\, \del_\theta$. The
second statement follows from the fact that $\cD(X)$ is generated by $\cO$ and $\Der(\cO)$.
\end{proof}

\begin{example}
\label{Ex6.11}
Let $E$ be the affine elliptic curve $y^2 = x^3 + x$. By Lemma~\ref{L6.2}, $\cD(E)$ is generated
by $\cO = \cO(E)$ and the derivation $\del_\theta$ with $\del_\theta(x) = 2y$ and $\del_\theta(y)
= 3x^2 + 1$. Consider the function $x \in \cO$. It defines a finite map $E \to \A^1$ of degree
$2$; its only zero on $E$ is the point $o = (0,0)$, where $y$ is a local parameter and $x$ has
order $2$. By Theorem~\ref{T6.1},
\[
  n(\Ann_{\cD(E)}(x)) = 2 .
\]
Let us verify this directly, and find the corresponding ideal. An operator of minimal order
annihilating $x$ is $L = x\, \del_\theta - 2y$, with $\Lambda(L) = \{o\}$. Since $\del_\theta(y)$
does not vanish at $o$, the classes of $\del_\theta^j$, $j \ge 0$, form a basis of $\Delta_o =
\cD/\fm_o \cD$. A direct calculation gives
\begin{align*}
  \del_\theta L & = x\, \del_\theta^2 - 6x^2 - 2 , \qquad \del_\theta^2 L = x\, \del_\theta^3 + 2y\,
  \del_\theta^2 - (6x^2 + 2)\, \del_\theta - 24 xy , \\
  \del_\theta^3 L & = x\, \del_\theta^4 + 4y\, \del_\theta^3 - 48 xy\, \del_\theta - 120 x^3 - 72 x .
\end{align*}
Hence $\del_\theta L \equiv -2$, $\,\del_\theta^2 L \equiv - 2\, \del_\theta$ and $L \equiv
\del_\theta^3 L \equiv 0$ modulo $\fm_o \cD$, so that $1$ and $\del_\theta^3$ belong to
$V_o(L)$. In fact,
\[
  V_o(L) = \spn\{1, \del_\theta^3\} .
\]
Indeed, $L = y\, \bar L$ up to a unit of $\cO_o$, where $\bar L$ is the normalization of $L$ at $o$,
and $J_o(L) = \{2\}$ by Theorem~\ref{T6.1}; hence $\Sigma(V_o(L)) = \{0, 3\}$ by
Corollary~\ref{C3.1}$(b)$. Formula \eqref{E37} gives $n_o(L) = 0 + 3 - 1 = 2$, as expected. Next,
$\fm_o^{-1} = \cO + \cO\,(y/x)$, and Theorem~\ref{T4.1} yields
\[
  \Ann_{\cD(E)}(x) = \cl(L) = \cD L + \cD\, \tfrac{y}{x} L + \cD\, \del_\theta^3 L + \cD\,
  \tfrac{y}{x}\, \del_\theta^3 L ,
\]
where
\[
  \tfrac{y}{x} L = y\, \del_\theta - 2(x^2 + 1) , \qquad \tfrac{y}{x}\, \del_\theta^3 L = y\,
  \del_\theta^4 + 4(x^2 + 1)\, \del_\theta^3 - 48 y^2\, \del_\theta - 120 x^2 y - 72 y .
\]
This ideal contains the operator of order three $\del_\theta^3 - 6x\, \del_\theta - 12 y$, which
equals $-\frac{1}{3}\bigl[\del_\theta^3\, \tfrac{y}{x} L - \tfrac{y}{x}\, \del_\theta^3 L - 3x\,
\del_\theta^2 L\bigr] + 2y\, \del_\theta L + 2x^2 L$. It follows that $\gr \Ann(x)$ contains the
ideal $(x \xi, y \xi, \xi^3)$ of $\gr \cD(E) = \cO[\xi]$, which has codimension $2$ in its
reflexive hull $\xi\, \cO[\xi]$; this confirms, via Theorem~\ref{T2.3}, that $n(\Ann(x)) \le 2$.
Finally, by Corollary~\ref{C4.1}, $\Ann_{\cD(E)}(x) = \cD(W, \cO) \cdot L$, where
\[
  W = \{ f \in \cO(E) \, : \, f(o) = 0 , \ (\del_\theta^3 f)(o) = 0 \} .
\]
\end{example}

\subsection*{Acknowledgements}
We are grateful to O. Chalykh and P. Etingof for stimulating discussions, many interesting questions and critical comments. Part of this research was conducted while the first author was visiting the Simons Laufer Mathematical Sciences Institute (formerly known as MSRI) as Research Professor for the program `Noncommutative Algebraic Geometry' in Spring 2024. He is very grateful to the Institute for its hospitality and financial support under the NSF Grant DMS-1928930. 

\appendix
\section{Canonical stairs and the local index formula}
\label{SA}

In this Appendix we give an alternative proof of Proposition~\ref{P5.1}, and hence of
Theorem~\ref{T5.1}. In fact, we give two different proofs: one combinatorial in Section~\ref{SA.2} and one homological in Section~\ref{SA.3}.
Unlike the proof given in Section~\ref{S5.2}, these two proofs are purely algebraic and do
not use the classification of ideals of the Weyl algebra in terms of the adelic Grassmannian. They are based on the notion of canonical stairs of an ideal introduced in \cite{BM84} (see also \cite{Maisonobe}).

\subsection{The stairs of a fat ideal}
\label{SA.1}
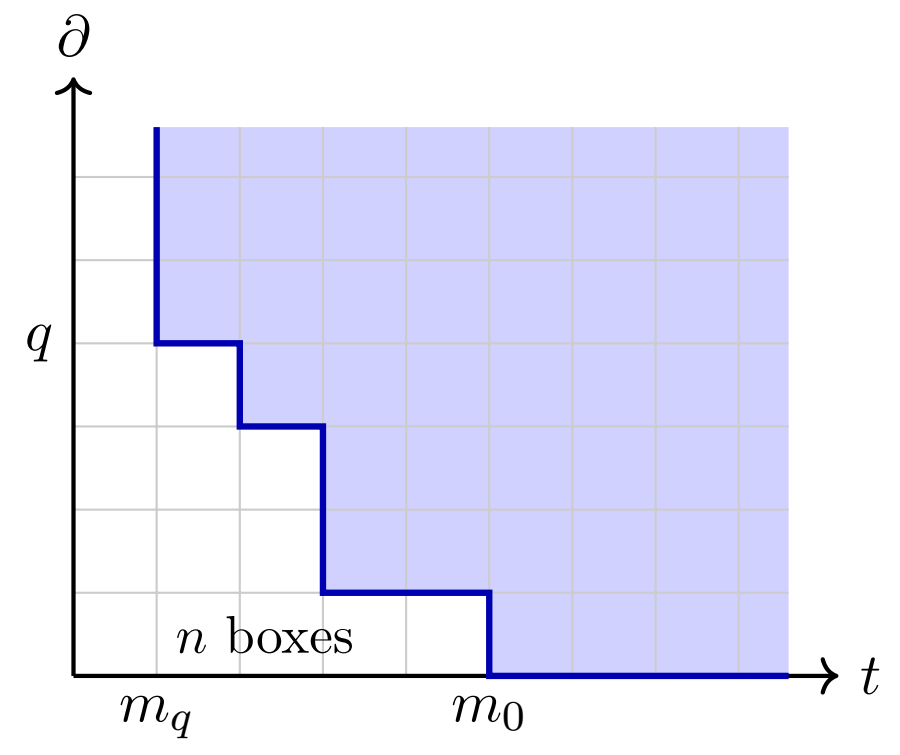
\begin{figure}[h]
\begin{center}
\begin{tikzpicture}[scale=0.55]
  \fill[blue!18] (5,0) -- (8.6,0) -- (8.6,6.6) -- (1,6.6) -- (1,4) -- (2,4) -- (2,3) -- (3,3) --
  (3,1) -- (5,1) -- cycle;
  \draw[gray!40] (0,0) grid (8.6,6.6);
  \draw[thick,->] (0,0) -- (9.2,0) node[right] {$t$};
  \draw[thick,->] (0,0) -- (0,7.2) node[above] {$\del$};
  \draw[very thick,blue!70!black] (8.6,0) -- (5,0) -- (5,1) -- (3,1) -- (3,3) -- (2,3) -- (2,4) --
  (1,4) -- (1,6.6);
  \node[below] at (5,0) {$m_0$};
  \node[below] at (1,0) {$m_q$};
  \node[left] at (0,4) {$q$};
  \node at (2.3,0.5) {\small $n$ boxes};
\end{tikzpicture}
\end{center}
\caption{The stairs of a fat ideal: the shaded region represents $\gr N$, and the number of boxes
under the stairs to the right of the line $t = m_q$ is the local index.}
\label{Fig}
\end{figure}

As in Section~\ref{S5.1}, write $\hD = \C[[t]][\del]$ for the ring of differential operators with coefficients in $ \C[[t]]$ and equip this ring with its order filtration $\{\hD_{\le k}\}$,
so that $\gr \hD \cong \C[[t]][\xi]$. Let $N$ be a right ideal of $\hD$, which is \emph{fat} in the sense that  $\,N \cap \C[[t]] \ne 0$. For each $k \ge 0$, denote by
\[
  \mathfrak{I}_k(N) := \{ a \in \C[[t]] \, : \, a\, \del^k + (\text{terms of order} < k) \in N \}
\]
the set of leading coefficients of the elements of $N$ of order $\le k$. These are nonzero ideals of $\C[[t]]$ such that 
$$
\mathfrak{I}_0(N) \subseteq \mathfrak{I}_{1}(N) \subseteq \ldots \subseteq \mathfrak{I}_k(N) \subseteq \mathfrak{I}_{k+1}(N) \subseteq \ldots
$$
Since every nonzero ideal of $\C[[t]]$ is
generated by a power of $t$, we have $\mathfrak{I}_k(N) = (t^{m_k})$ for a unique integer
$m_k \ge 0$, and
\[
  m_0 \ge m_1 \ge m_2 \ge \ldots \ge 0 .
\]
Let $q$ be the smallest integer such that $m_k = m_q$ for all $k \ge q$. Following \cite{BM84}, the finite set
\[
  \Stairs(N) := \{ (m_0, 0), (m_1, 1), \ldots, (m_q, q) \}
\]
is called the \emph{stairs} of $N$ (see Figure~\ref{Fig}). Clearly, this set is determined by
the partition $ (m_0, m_1, \ldots, m_q) $, which we will refer to as stairs as well. By construction,
\begin{equation}
\label{EA1}
  \gr N = \bigoplus_{k \ge 0} t^{m_k}\, \C[[t]]\, \xi^k \, \subseteq \, \C[[t]][\xi] ,
\end{equation}
which shows that the reflexive hull of $\gr N$ is the principal ideal generated by $t^{m_q}$,
and that
\begin{equation}
\label{EA2}
  \dim_{\C}\bigl[\, (\gr N)^{\vee\vee}\! / \gr N \,\bigr] = \sum_{k=0}^{q-1} (m_k - m_q) , \qquad
  \dim_{\C}\bigl[\, \hD_{\le k}/(N \cap \hD_{\le k}) \,\bigr] = \sum_{j=0}^{k} m_j .
\end{equation}
The first number in \eqref{EA2} is the `number of boxes under the stairs'. If $M$ is a fat right
ideal of $\cD(X)$ and $\lambda \in X$, then, by \eqref{E212}, the local Calogero--Moser index
$n_\lambda(M)$ is given by the first formula of \eqref{EA2} applied to $N = \widehat{M}_\lambda$.

\subsection{The stairs of an annihilator}
\label{SA.2}
Let $C_0$ be a finite-dimensional subspace of $\C[\del] \cong \hD/t\,\hD$, with degree set
$S = \Sigma(C_0) = \{s_1 < \ldots < s_g\}$, and let $N := \Ann_{\hD}(C_0)$. Recall that $N$ is a
fat right ideal of $\hD$, and $n(C_0)$ is defined by \eqref{E51}. The next theorem computes the
stairs of this ideal explicitly.
\begin{theorem}
\label{TA.1}
Let $m_0 \ge m_1 \ge \ldots$ be the stairs of the ideal $N = \Ann_{\hD}(C_0)$. Then
\begin{equation}
\label{EAstairs}
  m_k = \begin{cases}
    s_{g-k} + k + 1 , & 0 \le k < g , \\
    g , & k \ge g .
  \end{cases}
\end{equation}
In particular, we have $q \le g$ and $m_q = g$.
\end{theorem}

If $g = 0$, then $N = \hD$ and every $m_k = 0$, so  the theorem is immediate. We may therefore
assume that $g > 0$. Choose a basis $\{F_1, \ldots, F_g\}$ in $C_0$, with each
$F_i = \sum_j f_{ij}\, \del^j$ being monic of degree $s_i$. Put $f_{ij} = 0$ for $j < 0$ or
$j > s_i$. For each $k \ge 0$, consider the linear map
\[
  \Phi_k : \ \hD_{\le k} \to T_k := \bigoplus_{i=1}^g \C[\del]_{\le s_i + k} , \qquad P \mapsto
  (F_1 P, \ldots, F_g P) \ \ (\mathrm{mod}\ t\, \hD) .
\]
By definition, $\Ker \Phi_k = N \cap \hD_{\le k}$; hence, by \eqref{EA2},
\begin{equation}
\label{EA3}
  \mathrm{rk}\, \Phi_k = \sum_{j=0}^{k} m_j .
\end{equation}
To compute the rank of $\Phi_k$ we need a convenient parametrization of $\hD_{\le k}$. Every
$P \in \hD_{\le k}$ can be uniquely written in the form
\[
  P = \sum_{a \ge 0} \frac{t^a}{a!}\, P_a(\del) , \qquad P_a = \sum_{b=0}^k p_{a,b}\, \del^b .
\]
For each $r \ge 0$, we define the polynomial
\[
  \pi_r(s) := \sum_{a} \binom{s}{a}\, p_{a,\, k - r + a} \, \in \C[s] , \qquad \binom{s}{a} :=
  \frac{\ff{s}{a}}{a!} ,
\]
where the sum is taken over all $a$ with $\max\{0, r - k\} \le a \le r$. Since the polynomials
$\binom{s}{a}$ are linearly independent, the assignment $P \mapsto (\pi_r)_{r \ge 0}$ is a
bijection of $\hD_{\le k}$ onto $\prod_{r \ge 0} \Pi_r$, where
\[
  \Pi_r = \C[s]_{\le r} \quad \text{for } r \le k , \qquad \Pi_r = \ff{s}{r-k} \cdot \C[s]_{\le k}
  \quad \text{for } r > k .
\]
Finally, let $\vartheta := \del\, \frac{d}{d \del}$ denote the Euler operator on $\C[\del]$, so
that $\pi(\vartheta)\, \del^j = \pi(j)\, \del^j$ for any $\pi \in \C[s]$.

\begin{lemma}
\label{LA.1}
For any $F \in \C[\del]$ and $P \in \hD_{\le k}$, we have
\[
  F P \, \equiv \, \sum_{r \ge 0} \del^{k-r} \cdot [\pi_r(\vartheta) F](\del) \quad (\mathrm{mod}\ t\,
  \hD) .
\]
\end{lemma}

\begin{proof}
Since $[F(\del), t] = F'(\del)$, we have $F(\del)\, t^a \equiv F^{(a)}(\del)$ modulo $t\,\hD$.
Hence, if $F = \sum_j f_j\, \del^j$, then
\[
  F P \equiv \sum_{a \ge 0} \frac{1}{a!}\, F^{(a)}(\del)\, P_a(\del) = \sum_{a, b, j} \binom{j}{a}\,
  f_j\, p_{a,b}\, \del^{j - a + b} .
\]
Putting $r := k - b + a$, we can rewrite the last sum as $\sum_{r \ge 0} \sum_j f_j\, \pi_r(j)\,
\del^{j + k - r}$, which is the required formula. (Note that $\pi_r(j) = 0$ for integers
$0 \le j < r - k$, so that no negative powers of $\del$ occur.)
\end{proof}
\begin{proof}[Proof of Theorem~\ref{TA.1}]
Let $k \ge 0$. For $c \ge 0$, let $T_k^{(c)} \subseteq T_k$ be the subspace spanned by the
elements $\del^{s_i + k - c}$ of the $i$-th summand, $i = 1, \ldots, g$ (with $s_i + k - c \ge
0$); thus $T_k = \bigoplus_{c \ge 0} T_k^{(c)}$, and $\dim T_k^{(c)} = g$ for $c \le k$. By
Lemma~\ref{LA.1}, the component of $\Phi_k(P)$ in $T_k^{(c)}$ is the vector with coordinates
\begin{equation}
\label{EA4}
  \gamma_{i,c} = \pi_c(s_i) + \sum_{r=0}^{c-1} f_{i,\, s_i - c + r}\ \pi_r(s_i - c + r) .
\end{equation}
Thus $\Phi_k$ is `block triangular' with respect to the decompositions $\hD_{\le k} \cong \prod_r
\Pi_r$ and $T_k = \bigoplus_c T_k^{(c)}$, with diagonal blocks given by the evaluation maps
$\mathrm{ev}_c : \Pi_c \to T_k^{(c)}$, $\,\pi \mapsto (\pi(s_i))_i$. Only the blocks with
$c \le s_g + k$ contribute to the target, and input blocks with larger indices do not affect
these components.

If $c \le k$, then $\mathrm{ev}_c$ is the evaluation of polynomials of degree $\le c$ at $g$
distinct points, so that $\mathrm{rk}(\mathrm{ev}_c) = \min\{c+1, g\}$. If $c > k$, put
$g_c := |\{i : s_i \ge c-k\}|$. In this case, $\mathrm{ev}_c$ is the evaluation of the
polynomials $\ff{s}{c-k}\,\pi(s)$, $\deg \pi \le k$, at the $g_c$ points $s_i \ge c-k$, at
which $\ff{s}{c-k}$ does not vanish. Hence
\[
  \mathrm{rk}(\mathrm{ev}_c) = \min\{k+1, g_c\} \qquad (c > k) .
\]
Put $L:=s_g+k$ and let $M_k$ be the restriction of $\Phi_k$ to
$\bigoplus_{c=0}^{L}\Pi_c$. Then $\mathrm{rk}\,\Phi_k=\mathrm{rk}\,M_k$. Since
$\dim\Pi_c=\min\{c+1,k+1\}$ is nondecreasing and $\dim T_k^{(c)}$ is nonincreasing, with
$\dim T_k^{(L)}=1$, the integer
\[
  \ell := \min\{c:0\le c\le L,\ \dim\Pi_c\ge\dim T_k^{(c)}\}
\]
satisfies
\[
  \Ker(\mathrm{ev}_c)=0 \quad (c<\ell), \qquad
  \Coker(\mathrm{ev}_c)=0 \quad (c\ge\ell).
\]
Splitting the source and target at $\ell$, triangularity gives
\[
  M_k=\begin{pmatrix} A&0\\ B&D\end{pmatrix}, \qquad
  \Ker A=0, \qquad \Coker D=0.
\]
Thus $\Ker M_k=0\oplus\Ker D$, and hence
\[
  \mathrm{rk}\,\Phi_k
  =\sum_{c<\ell}\dim\Pi_c+\sum_{c=\ell}^{L}\dim T_k^{(c)}
  =\sum_{c=0}^{L}\mathrm{rk}(\mathrm{ev}_c).
\]

Write $G_h:=|\{i:s_i\ge h\}|$ for $h\ge1$, and put $R_k:=\mathrm{rk}\,\Phi_k$, with
$R_{-1}:=0$. The preceding rank formulas give
\[
  R_k=\sum_{c=0}^{k}\min\{c+1,g\}+\sum_{h\ge1}\min\{k+1,G_h\}.
\]
By \eqref{EA3},
\[
  m_k =R_k-R_{k-1}\,=\,\min\{k+1,g\}
       +\sum_{h\ge1}\bigl(\min\{k+1,G_h\}-\min\{k,G_h\}\bigr).
\]
For $0\le k<g$, we have $G_h\ge k+1$ if and only if $h\le s_{g-k}$; for $k\ge g$, we have
$G_h\le k$ for all $h$. Therefore
\[
  \sum_{h\ge1}\bigl(\min\{k+1,G_h\}-\min\{k,G_h\}\bigr)
  =\begin{cases}
    \displaystyle\sum_{h=1}^{s_{g-k}}1=s_{g-k},&0\le k<g,\\[4pt]
    0,&k\ge g.
  \end{cases}
\]
This proves \eqref{EAstairs}. The remaining assertions follow since the sequence $m_k$ stabilizes at $g$ for $k\ge g$.
\end{proof}

\begin{remark}
Theorem~\ref{TA.1} shows, in particular, that the reflexive hull of $\gr \Ann_{\hD}(C_0)$ is
generated by $t^{g}$ with $g = \dim C_0$, in agreement with Remark~2 after
Corollary~\ref{C4.2}.
\end{remark}

Proposition~\ref{P5.1} is an immediate consequence of Theorem~\ref{TA.1}:  
\begin{corollary}
\label{CA.1}
For any $\,C_0 \subset \C[\partial]$ with degree set $S$ and stairs $m_0 \ge m_1 \ge \ldots \ge m_q$,
\[
  n(C_0) = \sum_{k=0}^{q-1}(m_k-m_q)
  = \sum_{s\in S}s-\binom{g}{2} = \nu(S)\,.
\]
\end{corollary}

\begin{proof}
We have
$$
\displaystyle n(C_0)=\sum_{k=0}^{q-1}\,(m_k-m_q)\,
=\,\sum_{k=0}^{g-1}\,(s_{g-k}+k+1-g)=\sum_{s\in S}s-\binom{g}{2}\,,
$$
where the first equality follows from \eqref{EA2} and the second from \eqref{EAstairs}.
\end{proof}

\subsection{The canonical resolution}
\label{SA.3}
In this section, we give an alternative -- homological -- proof of Theorem~\ref{TA.1} using
the standard Gr\"obner resolution of an ideal.
We will keep the assumptions of Section~\ref{SA.2}. 
To simplify the notation we set $ A:=\C[[t]]$, and let $\,W := C_0^\perp \subseteq A\,$ be the annihilator of $C_0$ in $A$ (with respect to the pairing of
Lemma~\ref{L2.2}). We equip $W$ with the natural adic filtration
$\,W = F^0W \supset F^1W \supset \ldots\,$, where
$F^rW := W \cap t^r A$ for $r \ge 0$, and denote by
$\overline{W} := \bigoplus_{r\ge0} F^rW/F^{r+1}W$ its associated graded vector space.
By Lemma~\ref{L2.3}$(b)$, the set of gaps of $W$ is $S$. Since each subquotient $F^rW/F^{r+1}W$ is a subspace
of the one-dimensional space $t^rA/t^{r+1}A$, we have
\[
  \dim_{\C}(F^rW/F^{r+1}W) = \begin{cases}
    0 , & r \in S , \\
    1 , & r \notin S .
  \end{cases}
\]
The Poincar\'e series of $\overline{W}$ is therefore given by
\begin{equation}
\label{EAPoincare}
  P(\overline{W},z) := \sum_{r\ge0}\dim_{\C}(F^rW/F^{r+1}W)\,z^r
  = \frac{1}{1-z}-\sum_{s\in S}z^s\,.
\end{equation}
Our goal is to compute this Poincar\'e series in terms of the stairs of the ideal $ N = \Ann_{\hD}(C_0)\,$.

To this end, we choose operators $ P_0,\ldots,P_q $ in $N$ with symbols
$t^{m_0},\,t^{m_1}\xi,\,\ldots,\,t^{m_q}\xi^q$, respectively. By \eqref{EA1}, these operators
generate $N$ as a right ideal. The division algorithm in $\hD $ provides then the free resolution of $N$
(see \cite[Proposition 3]{BM84}):
\begin{equation}
\label{EA5}
  0 \to \hD^{\oplus q} \xrightarrow{\ \Psi\ } \hD^{\oplus(q+1)}
  \xrightarrow{\ (P_0,\ldots,P_q)\ } N \to 0 ,
\end{equation}
where the matrix $ \Psi $ is given by
\[
  \Psi=\begin{pmatrix}
  \del+b_0 & u_{01} & \cdots & u_{0,q-1} \\
  -t^{m_0-m_1} & \del+b_1 & \cdots & u_{1,q-1} \\
  0 & -t^{m_1-m_2} & \ddots & \vdots \\
  \vdots & & \ddots & \del+b_{q-1} \\
  0 & \cdots & 0 & -t^{m_{q-1}-m_q}
  \end{pmatrix}, \qquad b_i,u_{ki}\in A .
\]
Indeed, the operator $P_i\del-P_{i+1}t^{m_i-m_{i+1}}$ has order at most $i$; successive
division by $P_i,\ldots,P_0$ gives the $i$-th column of $\Psi$. The leading relations
$\xi e_i-t^{m_i-m_{i+1}}e_{i+1}$ generate the relations among $t^{m_j}\xi^j$, so order
division also proves exactness in the middle of \eqref{EA5}. Injectivity on the left follows
successively from the last row upwards. If $q=0$, the first term is understood to be zero.

We identify $N\otimes_{\hD}A$ with $W$ as follows. The pairing of Lemma~\ref{L2.2} gives
\[
  N=\{P\in\hD:P(A)\subseteq W\}.
\]
Choose $d\ge1$ such that $t^dA\subseteq W$, and put
\[
  T_d:=\sum_{j=0}^{d-1}\frac{(-t)^j}{j!}\,\del^j .
\]
Since $T_d(h)\equiv h(0)\pmod{t^dA}$ for $h\in A$, every $f\in W$ satisfies
$fT_d\in N$ and $(fT_d)(1)=f$. Thus evaluation at $1$ maps $N$ onto $W$.
Moreover, if $Q\del\in N$, then $(FQ)\del=0$ in $\hD/t\hD\cong\C[\del]$ for every
$F\in C_0$. Multiplication by $\del$ is injective on $\C[\del]$, so $Q\in N$. Consequently,
\[
  N\cap\hD\del=N\del, \qquad
  N\otimes_{\hD}A=N/N\del\ \xrightarrow{\ \sim\ }\ W, \quad [P]\mapsto P(1).
\]
Tensoring \eqref{EA5} with $A=\hD/\hD\del$ therefore gives an exact sequence of vector spaces
\begin{equation}
\label{EA6}
  0\to A^{\oplus q}\xrightarrow{\ \Psi\ }A^{\oplus(q+1)}
  \xrightarrow{\ \Phi\ }W\to0, \qquad
  \Phi(f_0,\ldots,f_q)=\sum_{j=0}^qP_j(f_j).
\end{equation}
Here injectivity of the first map again follows from the last row upwards.

To pass to associated graded spaces, we use the weight filtration on $\hD$.
For $0\ne P=\sum_j c_j(t)\del^j$, put
\[
  w(P):=\max_{c_j\ne0}\{j-\ord_t(c_j)\}.
\]
The initial part $\inn(P)\in\C[t][\del]$ is the sum of the terms
$\gamma_jt^{\ord_t(c_j)}\del^j$ with $j-\ord_t(c_j)=w(P)$, where
$\gamma_jt^{\ord_t(c_j)}$ is the leading term of $c_j$. We set $w(0)=-\infty$ and
$\inn(0)=0$. The relation $[\del,t]=1$ gives, for nonzero $P,Q$,
\[
  w(PQ)=w(P)+w(Q), \qquad \inn(PQ)=\inn(P)\inn(Q).
\]
Also, $w(P+Q)\le\max\{w(P),w(Q)\}$, with equality if the two weights differ, and
$\ord_t(Pf)\ge\ord_t(f)-w(P)$ for $f\in A$, with $\ord_t(0)=\infty$.

\begin{lemma}
\label{LA.2}
Put $Q_j:=\inn(P_j)$ and $d_i:=m_i-m_{i+1}$. Then
\[
  w(P_j)=j-m_j \quad (0\le j\le q), \qquad
  Q_i\del=Q_{i+1}t^{d_i} \quad (0\le i<q).
\]
\end{lemma}

\begin{proof}
Since $P_0=t^{m_0}$, the first assertion holds for $j=0$. Suppose it holds for $j\le i$.
The $i$-th column of $\Psi$ gives
\[
  P_i\del+P_ib_i+\sum_{k<i}P_ku_{ki}-P_{i+1}t^{d_i}=0.
\]
Now $w(P_i\del)=i-m_i+1$, while
\[
  w(P_ib_i)\le i-m_i, \qquad
  w(P_ku_{ki})\le k-m_k\le i-1-m_i \quad (k<i).
\]
Hence $w(P_{i+1})-d_i=i-m_i+1$, and taking initial parts gives
$Q_i\del=Q_{i+1}t^{d_i}$. This proves both assertions by induction.
\end{proof}

Write
\[
  a_j:=m_j-j, \qquad \beta_i:=m_i-i-1=a_i-1,
\]
and equip the $j$-th target summand in \eqref{EA6} with the shifted filtration
$F^nA[-a_j]:=t^{\max\{n-a_j,0\}}A$; use the shift $\beta_i$ on the $i$-th source summand.
By Lemma~\ref{LA.2}, $\Phi$ preserves these filtrations. The same holds for $\Psi$, since
\[
  w(\Psi_{ki})\le a_k-\beta_i.
\]
Equality holds on the diagonal and subdiagonal; for $k<i$, the right-hand side is at least
$2$, whereas $w(u_{ki})\le0$. Thus the initial matrix $\Psi_0$ is bidiagonal.
Let $ \bar{A} := \gr(A) =\C[\bar t\,]$, and identify $\bar t$ with $t$ when writing homogeneous
differential operators. Our convention is that the generator of $\bar{A}[-a]$ has degree $a$.

\begin{lemma}
\label{LA.3}
The maps in \eqref{EA6} are strictly compatible with the above filtrations. In particular,
there is an exact sequence of graded vector spaces
\begin{equation}
\label{EA7}
  0\to\bigoplus_{i=0}^{q-1}\bar{A}[-\beta_i]
  \xrightarrow{\ \Psi_0\ }\bigoplus_{j=0}^{q}\bar{A}[-a_j]
  \xrightarrow{\ \Phi_0\ }\overline{W}\to0,
\end{equation}
where $\Phi_0=(Q_0,\ldots,Q_q)$ and
$\Psi_0(fe_i')=f'e_i-t^{d_i}fe_{i+1}$.
\end{lemma}

\begin{proof}
Put $E:=t\del$. The identities of Lemma~\ref{LA.2} and $Et^d=t^d(E+d)$ give, by induction,
\[
  Q_j=t^{a_j}\prod_{h<j}(E+a_j-\beta_h).
\]
This formula is initially read in $\C[t,t^{-1}][\del]$. If $q>0$, minimality of $q$ gives
\[
  \beta_0>\beta_1>\cdots>\beta_{q-1}>a_q.
\]
It follows that
\[
  Q_q(1)=\prod_{h<q}(a_q-\beta_h)\,t^{a_q}\ne0.
\]
But $Q_q\in\C[t][\del]$, so $a_q\ge0$. The same conclusion is immediate for $q=0$.
In particular, all the shifts $a_j$ and $\beta_i$ are nonnegative.

Let $B:=\bigoplus_{j=0}^q\bar{A}[-a_j]$ and $C:=\bigoplus_{i=0}^{q-1}\bar{A}[-\beta_i]$.
The map $\Psi_0:C\to B$ is injective, as is seen from its last row upwards, and
$\Phi_0\Psi_0=0$. Viewing $\Phi_0$ for the moment as a map to $\bar{A}$, we have
\[
  Q_j(t^{n-a_j})=\prod_{h<j}(n-\beta_h)\,t^n \qquad (n\ge a_j).
\]
Since $a_i=\beta_i+1$ for $i<q$, comparison in degree $n$ gives
\[
  \dim B_n-\dim C_n=\mathrm{rk}(\Phi_{0,n})=
  \begin{cases}
    1,& n\ge a_q,\quad n\notin\{\beta_0,\ldots,\beta_{q-1}\},\\
    0,& \text{otherwise}.
  \end{cases}
\]
Indeed, for $n=\beta_i$ only the summands with $j>i$ contribute to $B_n$, and every
corresponding product vanishes. For $n\ge a_q$ outside this set, the coefficient of $Q_q$
is nonzero. Hence $\Ker\Phi_0=\im\Psi_0$.

It remains to check that the graded image is $\overline W$. For $w\in F^nW$, choose
$f=(f_j)\in A^{\oplus(q+1)}$ with $\Phi(f)=w$, using \eqref{EA6}. Put
$v(f):=\min_j\{\ord_t(f_j)+a_j\}$. If $v(f)=r<n$, the leading class of $f$ lies in
$\Ker\Phi_0=\im\Psi_0$. Lift a homogeneous preimage to the source of $\Psi$ and subtract
its image from $f$. This preserves $\Phi(f)$ and increases $v(f)$. After finitely many
steps we obtain a preimage with $v(f)\ge n$. Thus $\Phi$ is strict and its associated
graded map is onto $\overline W$. Finally, injectivity of $\Psi_0$ implies strictness of
$\Psi$, which completes the proof of \eqref{EA7}.
\end{proof}

Taking Poincar\'e series in \eqref{EA7} gives the required formula for \eqref{EAPoincare} in terms of the stairs.

\begin{proposition}
\label{PA.1}
With the above notation, we have
\begin{equation}
\label{EA8}
  P(\overline{W},z) \,= \,
    \sum_{j=0}^{q}\,\frac{z^{m_j-j}}{1-z}\,-\,\sum_{i=0}^{q-1}\,\frac{z^{m_i-i-1}}{1-z}\,.
\end{equation}
\end{proposition}

\begin{proof}
This follows from the exactness of \eqref{EA7} and the identity
$P(\bar{A}[-a],z)=z^a/(1-z)$.
\end{proof}

\begin{corollary}
\label{CA.2}
Proposition~\ref{PA.1} is equivalent to Theorem~\ref{TA.1}.
\end{corollary}

\begin{proof}
Put $r:=m_q-q$. The integers $\beta_i=m_i-i-1$, $0\le i<q$, are strictly decreasing
and greater than $r$. Formula \eqref{EA8} can therefore be rewritten as
\[
  P(\overline W,z)=\frac{z^r}{1-z}-\sum_{i=0}^{q-1}z^{\beta_i}.
\]
Its lowest term is $z^r$, so $r\ge0$. Comparing coefficients with \eqref{EAPoincare} shows that
\begin{equation}
\label{EAgaps}
  S=\{0,\ldots,r-1\}\,\sqcup\,\{\beta_0,\ldots,\beta_{q-1}\}.
\end{equation}
Consequently, $g=r+q=m_q$ and $\beta_i=s_{g-i}$ for $0\le i<q$. Hence
$m_i=s_{g-i}+i+1$ in this range. For $q\le i<g$, we have
$s_{g-i}=g-i-1$ and $m_i=g$, while $m_i=g$ for every $i\ge g$. This proves
\eqref{EAstairs}.

Conversely, assume \eqref{EAstairs}. Then $m_q=g$, $q\le g$ and $r=g-q$. If $r>0$,
the formula at $k=q$ gives $s_r=r-1$, whence $s_j=j-1$ for $1\le j\le r$.
For $i<q$, it gives $\beta_i=s_{g-i}$. Thus \eqref{EAgaps} holds, and substituting these values of $s \in S$ into \eqref{EAPoincare} yields \eqref{EA8}.
\end{proof}

The proof of Proposition~\ref{PA.1} is independent of Theorem~\ref{TA.1}; hence, combining this proposition with Corollary~\ref{CA.2}, we get another proof of Theorem~\ref{TA.1}. 

\bibliographystyle{amsplain}
\bibliography{preprint}

\end{document}